%% file: main_v7_arxiv.tex
\documentclass[hidelinks,onefignum,onetabnum]{siamart220329}
\newsiamthm{assump}{Assumption}
\crefname{subsection}{section}{sections}
\Crefname{subsection}{Section}{Sections}
\crefname{section}{Section}{Sections}
\Crefname{section}{Section}{Sections}
\crefrangeformat{figure}{Figures~#3#1#4 and~#5#2#6}
\Crefrangeformat{figure}{Figures~#3#1#4 through~#5#2#6}
\crefrangeformat{section}{Sections~#3#1#4 and~#5#2#6}
\Crefrangeformat{section}{Sections~#3#1#4 through~#5#2#6}
\crefrangeformat{subsection}{Sections~#3#1#4 and~#5#2#6}
\Crefrangeformat{subsection}{Sections~#3#1#4 through~#5#2#6}
\crefrangeformat{equation}{equations~(#3#1#4) through~(#5#2#6)}
\Crefrangeformat{equation}{Equations~(#3#1#4) through~(#5#2#6)}

\input{shared}
\ifpdf
\hypersetup{
  pdftitle={Fast, High-order Accurate Volume Potentials},
  pdfauthor={T. G. Anderson, M. Bonnet, L. M. Faria, and C. P\'erez-Arancibia}
}
\fi

\newcommand{\bol}{\boldsymbol}

\newcommand{\ney}{\boldsymbol{y}}
\newcommand{\nex}{\boldsymbol{x}}

\newcommand{\nez}{\boldsymbol{z}}

\newcommand{\de}{\,\mathrm{d}}
\newcommand{\e}{\operatorname{e}}

\newcommand{\inc}{\mathrm{inc}}

\newcommand{\andtext}{\quad\mbox{and}\quad}

\newcommand{\p}{\partial}

\newcommand{\lf}{\left}
\newcommand{\rg}{\right}

\newcommand{\R}{\mathbb{R}}
\newcommand{\C}{\mathbb{C}}
\newcommand{\N}{\mathbb{N}}

\newcommand{\vv}[1]{\boldsymbol{#1}}
\usepackage{amsmath}
\DeclareMathOperator*{\argminB}{argmin}

\usepackage{paralist}
\usepackage{multirow}
\usepackage{amsmath}
\usepackage{amsfonts}
\usepackage{amsmath}
\usepackage{amssymb}
\usepackage{graphicx}
\usepackage{mathrsfs}
\usepackage{upgreek}
\usepackage{booktabs}
\usepackage{cite}
\usepackage{hyperref}
\usepackage{mathtools}
\usepackage{enumitem}
\usepackage[inkscapearea=page]{svg}
\usepackage{adjustbox}

\makeatletter
\def\widebreve{\mathpalette\wide@breve}
\def\wide@breve#1#2{\sbox\z@{$#1#2$}%
     \mathop{\vbox{\m@th\ialign{##\crcr
\kern0.08em\brevefill#1{0.8\wd\z@}\crcr\noalign{\nointerlineskip}%
                    $\hss#1#2\hss$\crcr}}}\limits}
\def\brevefill#1#2{$\m@th\sbox\tw@{$#1($}%
  \hss\resizebox{#2}{\wd\tw@}{\rotatebox[origin=c]{90}{\upshape(}}\hss$}
\makeatletter

\makeatletter
\let\save@mathaccent\mathaccent
\newcommand*\if@single[3]{%
  \setbox0\hbox{${\mathaccent"0362{#1}}^H$}%
  \setbox2\hbox{${\mathaccent"0362{\kern0pt#1}}^H$}%
  \ifdim\ht0=\ht2 #3\else #2\fi
  }
\newcommand*\rel@kern[1]{\kern#1\dimexpr\macc@kerna}
\newcommand*\widebar[1]{\@ifnextchar^{{\wide@bar{#1}{0}}}{\wide@bar{#1}{1}}}
\newcommand*\wide@bar[2]{\if@single{#1}{\wide@bar@{#1}{#2}{1}}{\wide@bar@{#1}{#2}{2}}}
\newcommand*\wide@bar@[3]{%
  \begingroup
  \def\mathaccent##1##2{%
    \let\mathaccent\save@mathaccent
    \if#32 \let\macc@nucleus\first@char \fi
    \setbox\z@\hbox{$\macc@style{\macc@nucleus}_{}$}%
    \setbox\tw@\hbox{$\macc@style{\macc@nucleus}{}_{}$}%
    \dimen@\wd\tw@
    \advance\dimen@-\wd\z@
    \divide\dimen@ 3
    \@tempdima\wd\tw@
    \advance\@tempdima-\scriptspace
    \divide\@tempdima 10
    \advance\dimen@-\@tempdima
    \ifdim\dimen@>\z@ \dimen@0pt\fi
    \rel@kern{0.6}\kern-\dimen@
    \if#31
      \overline{\rel@kern{-0.6}\kern\dimen@\macc@nucleus\rel@kern{0.4}\kern\dimen@}%
      \advance\dimen@0.4\dimexpr\macc@kerna
      \let\final@kern#2%
      \ifdim\dimen@<\z@ \let\final@kern1\fi
      \if\final@kern1 \kern-\dimen@\fi
    \else
      \overline{\rel@kern{-0.6}\kern\dimen@#1}%
    \fi
  }%
  \macc@depth\@ne
  \let\math@bgroup\@empty \let\math@egroup\macc@set@skewchar
  \mathsurround\z@ \frozen@everymath{\mathgroup\macc@group\relax}%
  \macc@set@skewchar\relax
  \let\mathaccentV\macc@nested@a
  \if#31
    \macc@nested@a\relax111{#1}%
  \else
    \def\gobble@till@marker##1\endmarker{}%
    \futurelet\first@char\gobble@till@marker#1\endmarker
    \ifcat\noexpand\first@char A\else
      \def\first@char{}%
    \fi
    \macc@nested@a\relax111{\first@char}%
  \fi
  \endgroup
}
\makeatother

\usepackage[normalem]{ulem}

\definecolor{burntumber}{rgb}{0.54, 0.2, 0.14}

\makeatletter
\newcommand{\proofstep}[1]{%
  \par
  \addvspace{\smallskipamount}
  \textit{#1\@addpunct{.}}\enspace\ignorespaces
}
\makeatother
\usepackage[format=hang,font=small,labelfont={bf,sl},textfont=sl,width=0.95\textwidth]{caption}
\usepackage{subcaption}

\begin{document}
\maketitle

\begin{abstract}
  This article presents a high-order accurate numerical method for the evaluation of singular volume integral operators, with attention focused on operators associated with the Poisson and Helmholtz equations in two dimensions. Following the ideas of the density interpolation method for boundary integral operators, the proposed methodology leverages Green's third identity and a local polynomial interpolant of the density function to recast the volume potential as a sum of single- and double-layer potentials and a volume integral with a regularized (bounded or smoother) integrand. The layer potentials can be accurately and efficiently evaluated everywhere in the plane by means of existing methods (e.g.\ the density interpolation method), while the regularized volume integral can be accurately evaluated by applying elementary quadrature rules. Compared to straightforwardly computing corrections for every singular and nearly-singular volume target, the method significantly reduces the amount of required specialized quadrature by pushing all singular and near-singular corrections to near-singular layer-potential evaluations at target points in a small neighborhood of the domain boundary. Error estimates for the regularization and quadrature approximations are provided. The method is compatible with well-established fast algorithms, being both efficient not only in the online phase but also to set-up. Numerical examples demonstrate the high-order accuracy and efficiency of the proposed methodology; applications to inhomogeneous scattering are presented.
\end{abstract}



\maketitle 

\section{Introduction}\label{sec:intro}
This paper considers the numerical evaluation of  volume integral operators of the form
\begin{equation}\label{eq:vol_pot}
\mathcal V_k [f](\nex) \coloneqq \int_{\Omega} G_k(\nex,\ney) f(\ney)\de \ney,\quad \nex\in\R^2,
\end{equation}
over bounded domains $\Omega \subset \R^2$ with piecewise-smooth boundary $\Gamma$, where $f \in C^s(\overline\Omega)$, $s \in\N_0:=\N\cup\{0\}$, is a given function (termed the source density), and where the operator kernel $G_k$ is the free-space Green's function
\begin{equation}
G_k(\nex,\ney) \coloneqq \begin{cases}
\displaystyle-\frac{1}{2\pi} \log\left|\nex-\ney\right|,&k=0,\medskip\\
\displaystyle\frac{i}{4}H_0^{(1)}\left(k\left|\nex-\ney\right|\right), & k\neq0,
\end{cases}\label{eq:Green_Function}
\end{equation}
for the Laplace ($k = 0$) or Helmholtz ($k \neq 0$) equation. In fact, the method we propose is considerably more general than the examples presented, extending to (i) volume operators corresponding to these partial differential operators that possess stronger singularities {such as the gradient of $\mathcal{V}_k$ as appears e.g.\ in~\cite{martin2003acoustic},
\begin{equation}\label{eq:deriv_volpot}
 \mathcal W_k [\mathbf{f}](\nex) \coloneqq \int_\Omega \nabla_{\ney} G_k(\nex, \ney) \cdot {\bf f}(\ney) \de \ney, \quad \nex \in \R^2,
\end{equation}
where $\nabla_{\ney}$ denotes a gradient with respect to the $\ney$ variable} and $\bold f\in [C^s(\overline\Omega)]^2$, (ii) higher-dimensional analogues of these operators, and (iii) a variety of elliptic partial differential operators of interest in mathematical physics, e.g.\ in elasticity and fluids.  Despite the broad utility of these operators, for example, to solve a nonlinear or an inhomogeneous linear PDE, until very recently the computation of volume potentials has been relatively neglected in the context of complex geometries; in what follows we outline some difficulties of this problem and present some alternatives to the volume potential.

There are several well-known challenges to be met for the efficient and accurate evaluation of operators such as~\cref{eq:vol_pot} or~\cref{eq:deriv_volpot} in complex geometries---some mirroring those of the more well-studied problem of layer potential evaluation and some unique to volume problems. As is well known, the free-space Green's function $G_k$ given in~\cref{eq:Green_Function} features a logarithmic singularity (whose location is, of course, dependent on the target evaluation point $\nex$) that hinders the accuracy of standard quadrature rules for the numerical evaluation of the operator~\cref{eq:vol_pot} or~\cref{eq:deriv_volpot}. Despite the requirement for some degree of \emph{local} numerical treatment, the potential is a volumetric quantity depending \emph{globally} on the source density, and as such it is especially important to couple its computation to global fast algorithms---of which a variety have been developed including $\mathcal{H}$-matrix compression and their directional counterparts~\cite{Borm:03,Borm:17}, the fast multipole method (FMM)~\cite{Greengard:87}, and, more recently, interpolated-factored Green function methods~\cite{Bauinger:21}. Generally, it is also desirable to construct numerical discretizations of the operators which are efficient for repeated application (e.g.\ in schemes involving iterative linear solvers, Newton-type iterations for nonlinearities, or time-stepping) or with values that can be efficiently accessed (e.g.\ for the construction of direct solvers).

The prototypical, but by no means exclusive, use case for volume operators of the form~\cref{eq:vol_pot} is to tackle the slightly more specialized problem of producing a particular solution to an inhomogeneous constant-coefficient elliptic PDE. Restricting attention to work in this vein where the resulting homogeneous problem is solved with integral equation techniques, there have been a variety of methods proposed for producing a valid particular solution, typically on uniform grids using finite difference, finite element, or Fourier methods---see~\cite{Anderson:22a} for more details. With the proviso that methods based on finite elements or finite differences are often of limited order of accuracy, they nonetheless do apply to piecewise-smooth domains. Fourier methods, in turn, relying on various concepts of Fourier extension to generate a highly-accurate Fourier series expansion of the density function (after which the production of a particular solution is straightforward) in particular have received significant attention as they allow the use of highly-efficient FFT algorithms. Such Fourier-based methods are either fundamentally limited to globally smooth domains~\cite{Stein:16,Stein:22} or have only been demonstrated on such domains~\cite{Fryklund:18,Bruno:22}. Recently a scheme has been proposed~\cite{Fryklund:22} that allows extension on piecewise-smooth domains and, by coupling to highly-efficient volume Laplace FMMs~\cite{Ethridge:01} achieves an efficient and high-order accurate algorithm for producing a particular solution to Poisson's equation. Another significant application of the volume operator~\cref{eq:vol_pot} and~\cref{eq:deriv_volpot} is in the numerical solution of Lippmann-Schwinger integral equations that arise in some formulations of inhomogeneous scattering problems, see e.g.~\cite{saranen2001periodic,vainikko2006multidimensional,Bruno:04,Greengard:16}; for completeness we mention also the recent work~\cite{Bruno:19} for the same problem that proceeds along a different path by performing a direct PDE discretization rather than using volume potentials.

In this work, the volume potential is computed at each target point by utilizing Green's third identity to relate it to a linear combination of regularized domain integrals (over a standard triangularization), boundary integrals, and a local solution to the PDE that can be explicitly computed in terms of monomial sources. The method relies on local interpolants of the density function and explicit solutions to the underlying inhomogeneous PDE associated with those interpolated density functions.

Although the proposed approach is in spirit a natural extension of the density interpolation method~\cite{perez2019harmonic,perez2019planewave,perez2018plane,gomez2021regularization,faria2021general}, it shares several elements with previous methods. We first describe prior work in constructing local PDE solutions and then discuss the manner in which such solutions have historically been coupled to yield particular solutions valid on the entirety of the domain of interest. The first work we are aware of in which particular solutions were constructed using monomials is that of reference~\cite{Atkinson:85}; of course, the solutions are also polynomials. Subsequently, a variety of contributions studied combinations of polynomial sources and various elliptic partial differential operators; for monomials there exist solutions for Laplace~\cite{Atkinson:85}, Helmholtz~\cite{Golberg:03monom}, and elasticity operators~\cite{Matthys:96monom}. A particular solution formula that is applicable to second-order constant coefficient partial differential operators featuring a non-vanishing zeroth-order term and a monomial right-hand side was presented in~\cite{Dangal:17monom}. Many of these methods can be considered to be closely related to the dual reciprocity method~\cite{Nardini:83,Partridge}---such schemes rely on the density function being \emph{globally} well-approximated by a collection of simple functions (e.g.\ polynomials, sinusoidal, or radial functions), and implicitly rely on the accurate and stable determination of coefficients in such expansions. This latter task has not generally been successfully executed to even moderate accuracy levels for complex geometries and/or density functions.

Thus, for an arbitrary complex geometry, in order to use the simple basis function expansions of the density function envisioned in reference~\cite{Atkinson:85} for a stable, high-order solver, one is naturally led to the use of subdivisions of the domain, wherein \emph{global} solutions are stitched together from solutions defined with \emph{local} data. The case of a rectangular domain that can be covered entirely by boxes is considered in reference~\cite{GreengardLee:96} to construct a high-order Poisson solver. That work, similarly to the present method, uses local polynomial solutions (corresponding to Chebyshev polynomial right-hand sides) and layer potential corrections to provide a global solution. As described there, the local solution defined with data supported on a given element undergoes jumps across the boundary of that element, and it follows from Green's identities that the particular solution, i.e. the volume potential, can be represented outside of the local element in terms of single- and double-layer potentials with polynomial densities defined explicitly in terms of the local element solution. This approach was also taken up using multi-wavelets on domains covered by a hierarchy of rectangular boxes in reference~\cite{Averbuch:00} as well as, in large part, in the three-dimensional works~\cite{Israeli:02, McCorquodale:07}. Like these other works, the method proposed here uses a piecewise-polynomial approximation of the source.

We discuss next a few relevant distinctions between the present methodology and other more closely-related schemes. The scheme proposed in the present work differs from that introduced in reference~\cite{Anderson:22a} (see also~\cite{Atkinson:85}), where the volume potential is evaluated over irregular domain regions using a direct numerical quadrature of the singular and near-singular integrals. {While the method presented in~\cite{Anderson:22a}} requires only singular quadrature rules appropriate for the singular asymptotic behavior of the function $G_k$, its applicability may not include some of the more severely singular kernels which the present methodology and the density interpolation method naturally treat~\cite{gomez2021regularization,faria2021general}; additionally, the efficiency of that method suffers when gradients become large (e.g.\ $|k|$ large) with a wider near-field volume correction area needed. We turn next to methods that like ours leverage a polynomial solution to the underlying PDE and transform volume potentials to layer potentials. Unlike the previously-mentioned works~\cite{GreengardLee:96,Averbuch:00,Ethridge:01,Israeli:02,McCorquodale:07} that rely on the method of images and/or pre-computed solutions for singular and near-singular evaluation points---since the grid and therefore the target locations have known structure in those settings---the present work applies to unstructured meshes wherein evaluation points may lay in arbitrary locations relative to element boundaries. Our method shares some elements with the work~\cite{ShenSerkh:22} that has contemporaneously appeared: both use Green's third identity for volume potential evaluation over unstructured meshes, but, unlike that work, the class of near-singular layer potential evaluation schemes which we employ extends to three dimensions, treats operators with stronger singularities, and is, in principle, kernel-independent (see e.g.~\cite{faria2021general}). We stress, however, that with regard to choice of methods for layer potential evaluation the scheme is agnostic. Indeed, alternative techniques such as~\cite{Barnett:14,klockner2013quadrature,helsing2008evaluation,barnett2015spectrally,beale2016simple,bao2024singularity,siegel2018local,af2018adaptive,zhu2022high} can be used instead of the general-purpose density interpolation method~\cite{faria2021general}.

\begin{figure}
    \centering
    \includegraphics[width=.35\textwidth]{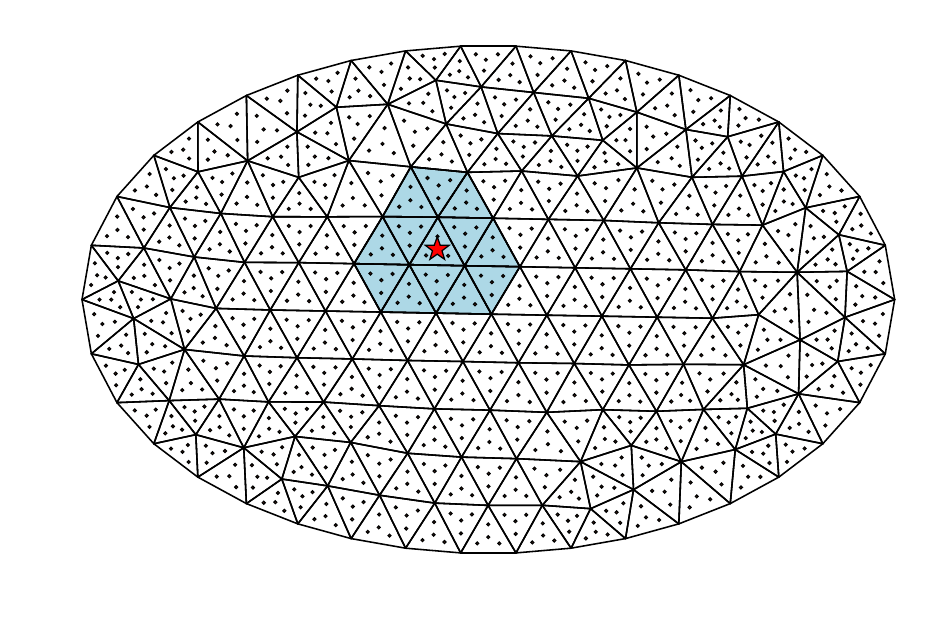}
    \includegraphics[width=.35\textwidth]{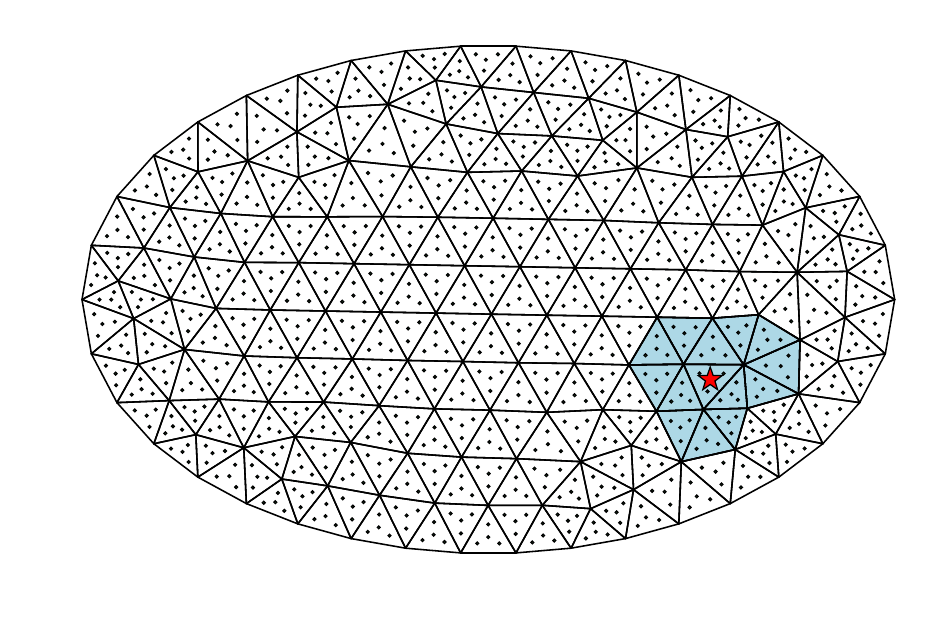}\\
    \includegraphics[width=.35\textwidth]{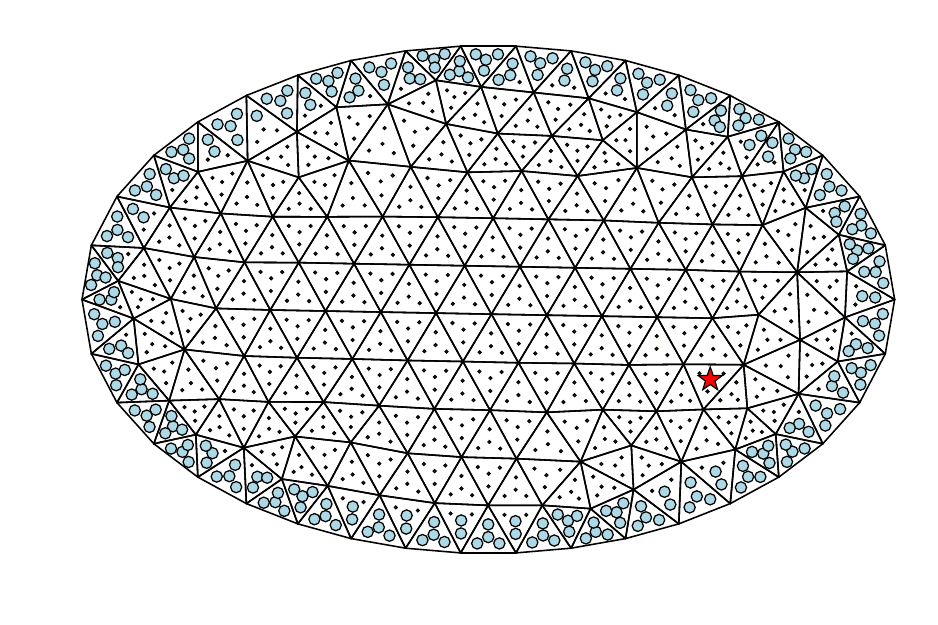}
    \includegraphics[width=.35\textwidth]{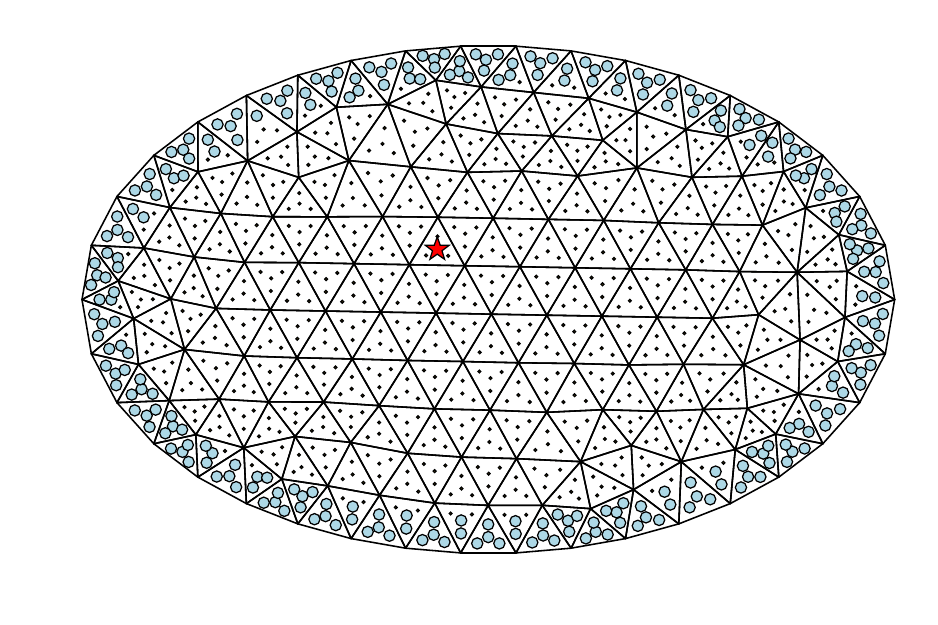}
    \caption{Comparison of the regions for which specialized quadrature is needed in two volume quadrature methods. Top: With the use of standard quadrature rules for $\mathcal{V}_k[f]$, nearby triangles for which the volume integral will be inaccurate for \emph{specific} evaluation points collocated at a quadrature node and marked with a red star. Computation of the potential at all nodes displayed in the volume requires developing many (near-)singular quadrature rules for every triangle. Bottom: In the proposed methodology, all volume evaluation points which require near-singular layer potential evaluation and indeed specialized quadrature of any kind are marked in blue; no specialized quadrature is needed for evaluation points inside the domain.}.
    \label{fig:sing_nearsing_nodes}
\end{figure}

An important \emph{novel} feature of this method is that it results in a significant reduction in the number of \emph{singular} corrections compared to previous volume potential schemes (even as the total number of volume corrections remains the same). This is made possible by the use of Green's identity over a large region $\Gamma = \partial \Omega$ (illustrated in \Cref{fig:sing_nearsing_nodes}) as opposed to individual element-wise boundaries, with the effect that, on the one hand, near-singular layer potential quadrature (e.g.\ with the density interpolation method) is only required in a thin region abutting $\Gamma$ while, on the other hand, volume quadrature for resulting regularized volume integrals can be effected to high-order accuracy with standard quadratures.

To support this last claim the paper presents a detailed error analysis culminating in the main \Cref{thm:tri_error_analysis} with high-order ($h$-refinement) estimates that are observed in the experiments  performed for moderate interpolation degrees $n\in\{0,1,2,3,4\}$ (certain numerical instabilities may be more significant for large interpolation degrees, an issue that is discussed in detail in~\Cref{sec:limitations} with a prescription outlined for managing the instability). For context on the significance of these results, the guarantees established here---most precisely, the main theorem's supporting \Cref{lem:ordquad_convergence_farfield_optimal} that establishes high order accuracy  (an order of accuracy that is, in a sense made precise, optimal) as the mesh size $h \to 0$ when a quadrature rule of fixed order $m$ is used for all
elements not containing the evaluation point---must be understood in contrast to the
state of the theory for methods for nearly-singular integral
evaluation in widespread use. Indeed, in the community it is routine, both for volume and layer potentials, to refer to the `near-field' as a neighborhood of diameter proportional to the element / local patch / panel size $h$; the complement of this neighborhood is called the `far-field'. Common practice is to use specialized corrective quadratures \emph{only} over an
$\mathcal{O}(h)$-sized `near field' in conjunction with a fixed-order (\emph{not} spectral) quadrature rule for the `far-field'. Unfortunately, such methods are often accompanied by unjustified claims of high-order accuracy because the `far-field' encroaches closer to the singularity under $h$-refinement and errors in the `far-field' integrals are not controlled to high-order. \emph{Simply put, as the mesh-size $h\to 0$, either (i) the order of accuracy of far-field
quadratures is expected to eventually degrade, or (ii) if, instead, true high-order accuracy is desired then corrections are made on a growing number of elements or points  (with the need for the inclusion of these corrections impacting on efficiency or even leading to superlinear complexity~\cite{BrunoGarza:20,ying2006high})}. For example, the representative reference~\cite{siegel2018local} even notes ``crucially, this scaling of parameters with $h$
also leads to the desired $\mathcal{O}(1)$ work per target in the local correction step''. A certainly incomplete list of further recent examples in this `locally-corrected' vein include~\cite{greengard2021fast,zhu2022high,sushnikova2023fmm,bremer2015high,bremer2012nystrom,bremer2013numerical} and the volume potential schemes~\cite{Anderson:22a,ShenSerkh:22}.  Nevertheless, the practical effect of this issue may be concealed by the higher-order quadratures frequently used in the far-field so that significant impacts on accuracy are not necessarily seen in practice (reference~\cite[Fig.\ 9]{ShenSerkh:22} does report seriously sub-optimal convergence rates in light of the theory in this paper). In connection with this discussion it is interesting to observe that in the boundary element literature~\cite[\S 5.3]{Sauter2010} similar problems are controlled through use of variable quadrature orders, that is, the use of quadrature orders that increase logarithmically with $h$-refinement.

The remainder of this paper is structured as follows. An introduction to basic identities and background information on quadrature is presented in \Cref{sec:prelim}. \Cref{sec:method_overview} presents a general overview of the method including its compatibility with fast algorithms. The resulting complexity analysis is given in \Cref{sec:complexity}. \Cref{sec:tri_interpolant_constr} focuses on the construction of the interpolation polynomial $f_n$ (the conditioning of the linear system used to recover the interpolant is analyzed in \Cref{sec:conditioning}). A detailed error analysis presented in \Cref{sec:tri_regularization_analysis}, whose outcome is \Cref{thm:tri_error_analysis}, establishes the regularizing effect of the polynomials on the domain integrals and ultimately gives rise to a provably high-order accurate numerical scheme. Finally, numerical examples of the proposed approach are presented in \Cref{sec:numer}, some limitations due to certain instabilities and methods to manage them are given in \Cref{sec:limitations}, and a few comments on future work are given in the conclusion.

\section{Preliminaries}\label{sec:prelim}

In what follows we briefly outline the principles underlying the source approximation and singular integration strategy employed in this paper. The method relies on polynomial approximations of the density function but is somewhat agnostic to the interpolation strategy. We present a specific strategy wherein approximation is performed by Lagrange interpolating polynomials over the entirety of an element $K$ in a triangulation $\mathcal T_h$.
The method locally constructs, on every element $K$, a piece-wise (regularizing) polynomial approximation to the source, but, as will be explained in detail, is represented in a global basis on $\Omega$. Without loss of generality, we assume the domain's centroid lies at the origin, which mitigates the scale of this (polynomial) basis over $\Omega$.

Next, the method builds associated regularization polynomial solutions that are particular solutions to the PDE with the interpolant as a right-hand side, and the interpolants and solutions are then used to regularize the volume integral in~\cref{eq:vol_pot}. Precisely, we first explicitly construct a family of polynomial PDE solutions $\Phi_n(\cdot;K): \R^2  \to \C$, the family being parametrized by $K\in  \mathcal T_h$, that satisfies
\begin{equation}\label{eq:reg_function_ii}
  (\Delta+k^2) \Phi_n(\cdot; K)  = f_n(\cdot; K)\quad\text{in}\quad \R^2,
\end{equation}
where $f_n(\cdot;  K)$ is the polynomial of total degree $n$ that interpolates $f$ on a discrete set of points within the (possibly curvilinear) triangle $ K$ which is closest to the target point $\nex\in\R^2$ (typically, $\nex\in K$); details about the construction of $f_n$ and $\Phi_n$ are given in~\Cref{sec:tri_interpolant_constr} and \cite{anderson2023particular} respectively.

Employing Green's third identity and the solution $\Phi_n$ we arrive at the relation
\begin{equation}\label{eq:GF}
  \mu(\nex)\Phi_n(\nex; K)=-\int_\Omega G_k(\nex,\ney) f_n(\ney; K)\de \ney -\int_{\Gamma} \left\{\frac{\p G_k(\nex,\ney)}{\p \nu(\ney)}   \Phi_n(\ney; K)-G_k(\nex,\ney) \frac{\p \Phi_n(\ney; K)}{\p \nu(\ney)}  \right\}\de s(\ney)
,
\end{equation}
where $$\mu(\nex) \coloneqq \begin{cases}1&\nex\in\Omega,\\
  \gamma(\nex) &\nex\in\Gamma,\\
0&\nex\in\R^2\setminus\overline\Omega,\end{cases}$$ with $\gamma(\nex)$ denoting a function proportional to the interior angle at $\nex\in\Gamma$ that is equal to the familiar $\frac12$ where $\Gamma$ is continuously differentiable~\cite[\S 8.1]{Atkinson}. As usual, normal derivatives in \eqref{eq:GF} are taken with respect to the exterior unit normal to $\Gamma$ that we denote by $\nu$. From adding~\cref{eq:GF} to~\cref{eq:vol_pot} it follows that the volume integral can be recast as
\begin{multline}
  \mathcal V_k [f](\nex) = \mathcal V_k [f-f_n(\cdot;K)](\nex)- \\
  \int_{\Gamma} \left\{\frac{\p G_k(\nex,\ney)}{\p \nu(\ney)}   \Phi_n(\ney; K)-G_k(\nex,\ney) \frac{\p \Phi_n(\ney; K)}{\p \nu(\ney)}  \right\}\de s(\ney)-\mu(\nex)\Phi_n(\nex; K). \label{eq:reg_formula}
 \end{multline}

In addition, the regularization of the operator $\mathcal W_k$ can be achieved similarly by interpolating the vector density~$\mathbf{f}$ component-wise. Indeed, letting $\mathbf f_n(\cdot; K)$ denote the corresponding polynomial interpolant of $\mathbf{f}$ at/around the target point $\nex\in\R^2$,  and letting $\boldsymbol\Phi_n(\cdot; K)$ denote a solution of the vector PDE
\begin{equation}\label{eq:reg_function_W}
  (\Delta+k^2) \boldsymbol\Phi_n(\cdot; K)  = \mathbf{f}_n(\cdot; K)\quad\text{in}\quad \R^2,
\end{equation}
we obtain from the identities
\begin{equation*}\label{eq:GF_W}
  \mu(\nex)\boldsymbol{\Phi}_n(\nex; K)=-\int_\Omega G_k(\nex,\ney) {\bf f}_n(\ney; K )\de \ney -\int_{\Gamma} \left\{\frac{\p G_k(\nex,\ney)}{\p \nu(\ney)} \boldsymbol{\Phi}_n(\ney; K)-G_k(\nex,\ney) \frac{\p \boldsymbol{\Phi}_n(\ney; K)}{\p \nu(\ney)}  \right\}\de s(\ney)
  \end{equation*}
  and $\nabla_{\nex} G_k(\nex,\ney) = -\nabla_{\ney}G_k(\nex,\ney)$,
the following regularized expression for $\mathcal W_k[\bf f]$:
\begin{multline}
  \mathcal W_k [\mathbf{f}](\nex) = \mathcal W_k [\mathbf{f}-\mathbf{f}_n(\cdot;K)](\nex)+ \\
  \int_{\Gamma} \left\{\nabla_{\nex}\frac{\p G_k(\nex,\ney)}{\p \nu(\ney)}   \cdot\boldsymbol\Phi_n(\ney; K)-\nabla_{\nex}G_k(\nex,\ney)\cdot \frac{\p\boldsymbol \Phi_n(\ney; K)}{\p \nu(\ney)}  \right\}\de s(\ney)+\mu(\nex)\operatorname{div}\boldsymbol\Phi_n(\nex; K). \label{eq:reg_formula_W}
\end{multline}

Just like the case of $\mathcal V_k[f]$, the interpolation property of $\mathbf{f}_n$ enables the evaluation of the regularized volume integral in~\cref{eq:reg_formula_W} using a quadrature rule independent of the target point location. Using~\cref{eq:reg_formula_W} in conjunction with existing methods for evaluating the layer potentials we obtain an accurate and efficient algorithm for evaluating $\mathcal W_k[\mathbf{f}]$. We mention in passing that the divergence theorem can be useful in reducing computational effort when both $\mathcal{V}_k$ and $\mathcal{W}_k$ are needed.

\begin{remark}\label{rem:regularity_f}
  In general, in the derivations and  in the analysis that follow, we will assume smoothness of the density functions $f$ and $\mathbf{f}$ over $\Omega$, e.g.\ for a given integer $n\in\N_0$ we assume $f \in C^{s}(\overline\Omega)$ and $\mathbf{f}\in [C^{s}(\overline\Omega)]^2$, for a certain integer $s>n$ (specific $s$ values can be found in \Cref{thm:tri_error_analysis} but we do not claim the required value of $s = n + 3$ is in general tight; if all elements in the mesh are convex the required regularity for the theory is only $s = n + 1$  and should be sharp). In addition, the methodology applies in the event the source densities are only piecewise smooth over~$\Omega$, in which case~$\Omega$ must be considered as a union of disjoint domains where~$f$ (resp. $\bold f$) is smooth. Thus, provided the location of possible discontinuities of a density function is known, such concerns can be addressed using the proposed technique.
\end{remark}

\subsection{Elementary numerical integration}\label{sec:tri_mappings_quadratures}
We will work over a tessellation of the domain that consists of both straight and curved triangles.
Let $\mathcal{T}_h = \left\{K_\ell\right\}_{\ell=1}^{L}$ be a triangulation of $\Omega$, i.e., $\overline\Omega=\cup_{\ell=1}^LK_\ell$ and $\mathring K_\ell\cap \mathring K_{\ell'}=\emptyset$, $\ell\neq\ell'$, which we assume throughout consists of elements with maximum diameter (in the sense of the metric distance) $h>0$ and with a minimum diameter of inscribed circles of $\rho>0$; at times, local element diameters $h_K$ and  inscribed circle sizes $\rho_K$ will be used.

\subsubsection{Mappings for straight and curvilinear triangles}\label{sec:curved_triangles}
All elements are mapped from the closed reference triangle $\widehat{K}$ with vertices
\begin{equation*}\label{eq:vertices}
  \left\{\left(-\frac12,-\frac{\sqrt3}2\right),(1,0),\left(-\frac12,\frac{\sqrt3}2\right)\right\}
\end{equation*}
and the transformation to physical coordinates is denoted by $\vv T: \widehat{K} \to K$. Here, and throughout this paper, superscript-ed hat notation, e.g.\ $\widehat{\vv T}$, denotes quantities in reference space on the reference triangle $\widehat K$; undecorated symbols, e.g. $\vv{T}$ or $K$, denote maps into or quantities in physical space.

When $K$ is a straight triangle, $\vv T$ is an affine transformation that maps the vertices of $\widehat{K}$ to the vertices of $K$. When $K$ is a curved triangle with vertices $\vv v_1 = (v_{11}, v_{12})$, $\vv{v}_2 = (v_{21}, v_{22})$, and $\vv{v}_3 = (v_{31}, v_{32})$ that has exactly one curved edge (i.e.\ on the boundary), say along the edge connecting $\vv{v}_1$ to $\vv{v}_2$, we employ, as recently suggested~\cite{Anderson:22a}, the blending mappings~\cite{Gordon:73a,Gordon:73b} defined on the standard simplex $\Delta$. To do this, let ${\vv T}_\Delta$ denote the affine map ${\vv T}_\Delta: \widehat{K} \to \Delta$ and let the curved element map ${\vv T}: \widehat{K} \to K$ be given by ${\vv T} = {\vv T}_\gamma\circ {\vv T}_\Delta,$ with
\begin{equation*}
    \mkern9mu {\vv T}_\gamma (\eta_1, \eta_2) =
  \begin{pmatrix}
    (1 - \eta_1 - \eta_2) v_{11} + \eta_1 v_{21} + \eta_2 v_{31} + \frac{1 - \eta_1 - \eta_2}{1 - \eta_1}\left(\gamma_1(\eta_1) - (1 - \eta_2) v_{11} - \eta_1 v_{21}\right)\\
    (1 - \eta_1 - \eta_2) v_{12} + \eta_1 v_{22} + \eta_2 v_{32} + \frac{1 - \eta_1 - \eta_2}{1 - \eta_1}\left(\gamma_2(\eta_1) - (1 - \eta_1) v_{12} - \eta_1 v_{22}\right)
  \end{pmatrix}.
\end{equation*}
Here, $\bol\gamma=(\gamma_1,\gamma_2): [0, 1] \to \Gamma\cap K$ is the parametrization of the coordinates of the curved edge connecting $\vv{v}_1$ and $\vv{v}_2$ that satisfy $\bol\gamma(0) = \vv v_{1}$ and $\bol\gamma(1) = \vv v_{2}$. It can easily be seen that $\vv T$ is a $C^1$-invertible map of $\widehat{K}$ onto $K$ for each element $K$; we denote $\left\{\vv T_\ell\right\}_{\ell=1}^{L}$ the set of all such curvilinear and affine mappings for elements in~$\mathcal{T}_h$.\looseness=-1

\subsubsection{Quadrature for boundary integrals and for smooth functions on mapped triangles}
The regularized volume integral operators in~\cref{eq:reg_formula} and~\cref{eq:reg_formula_W} require a quadrature rule for smooth functions defined on the elements $K_\ell\in\mathcal T_h$, $\ell=1,\ldots,L$. For sufficiently regular functions $\phi: \Omega \to \C$, a volume integral can be expressed via the sum, over all $L$ elements of the triangulation $\mathcal{T}_h$,
\begin{equation}\label{eq:tri_smooth_integral_decomp}
  \int_\Omega \phi(\nex)\de\nex \approx \sum_{\ell=1}^L \int_{\widehat{K}} \phi({\vv T}_\ell({\widehat\nex})) \left|\operatorname{det} \mathsf J_\ell({\widehat\nex})\right|\de\widehat\nex,
\end{equation}
where $\mathsf J_\ell$ is the Jacobian matrix of the transformation $\vv T_\ell$ prescribed above. Quadrature for each of the integrals over $\widehat{K}$ in~\cref{eq:tri_smooth_integral_decomp} is performed via the $q_n$-numbered (see~\cref{eq:qn_def}) Vioreanu-Rokhlin quadrature nodes~\cite{Vioreanu:14}
\[
  \widehat{\mathcal{I}}_n = \left\{\widehat\nex_{j}: \widehat\nex_{j} \in \widehat{K}\right\}_{j=1}^{q_n}
\]
with associated weights $\left\{\upomega_j\right\}_{j=1}^{q_n}$. A
Vioreanu-Rokhlin quadrature rule is capable of exactly integrating polynomials of total
degree $m$ (values of $m$ are tabulated~\cite{Vioreanu:14}). Use of other
quadratures is possible but we prefer Vioreanu-Rokhlin nodes because they
doubly function for us as relatively well-conditioned interpolation nodes (see \Cref{sec:conditioning} and \Cref{fig:shiftscale0}). Letting $N = L q_n$, real-space nodes are obtained via the mappings $\{\vv T_\ell\}_{\ell=1}^L$ yielding the node-weight set
\begin{equation}\label{eq:globnodesweights}
\begin{split}
    \scalebox{1.25}{$\chi$} &= \{({\vv \xi}_r, \omega_r)\}_{r=1}^N\\
        &= \bigcup_{\ell = 1}^L\{{\vv T}_\ell(\widehat\nex_j), \upomega_j |\operatorname{det} {\mathsf J}_\ell(\widehat\nex_j)|\}_{j=1}^{q_n}
\end{split}
\end{equation}
of global quadrature nodes and weights,
and the quadrature approximation can thus be written as
\begin{equation}\label{eq:standard_QR}
  \mathcal Q_\Omega[\phi] \coloneqq \sum_{\substack{(\vv{\xi}_j,\omega_j) \in \scalebox{1.0}{$\chi$}\\ \vv{\xi}_j \in \Omega}} \omega_j \phi(\vv{\xi}_j) = \sum_{j=1}^N\omega_j\phi(\vv\xi_j) \approx\int_\Omega \phi(\ney)\de \ney.
\end{equation}
(Note that the notation introduced in~\eqref{eq:standard_QR} allows us to refer to quadrature approximations $\mathcal{Q}_{E}$ applied to integrals over subsets $E\subset\Omega$.)

The layer potentials in~\cref{eq:reg_formula}, on the other hand, can be accurately evaluated by means of existing techniques~\cite{gomez2021regularization,perez2019harmonic,perez2018plane,faria2021general} (or the other layer potential schemes given in the introduction).  For sufficiently regular trace functions $\varphi,\psi:\Gamma\to\C$, such techniques yield high-order approximations
\begin{equation}\label{eq:quadrature_layerpot}
\mathcal B_\Gamma[\varphi,\psi](\nex)\approx \int_{\Gamma} \left\{\frac{\p G_k(\nex,\ney)}{\p \nu(\ney)}  \varphi(\ney)-G_k(\nex,\ney) \psi(\ney) \right\}\de s(\ney),\quad \nex\in\R^2,
\end{equation}
where $\mathcal{B}_\Gamma$ stands for (any) sufficiently accurate approximation of the boundary integrals. When $\nex$ is well-separated from $\Gamma$, of course, standard quadrature rules are sufficient.

\section{Method description and algorithm sketch}\label{sec:method_overview}
In this section, we outline a numerical procedure that exploits~\cref{eq:reg_formula} to develop an efficient, highly accurate method to numerically evaluate~\cref{eq:vol_pot} everywhere in $\overline\Omega$, i.e., here we assume  that  $\nex\in K\in\mathcal T_h$; the method directly extends (component-wise) to allow use of~\cref{eq:reg_formula_W} in evaluating~\cref{eq:deriv_volpot}.  The generalization of the procedure to evaluate volume potentials over narrow regions outside $\overline\Omega$ (where the volume integrals remain nearly singular) is discussed in \Cref{rem:at_boundary}.

  To explain more concretely, we introduce the notation used in the remainder of the paper: we will make use of the standard multi-index notation where, for any $\alpha = (\alpha_1,\alpha_2)\in\N_0^2$, we set $\alpha! = \alpha_1!\alpha_2!$, $|\alpha| = \alpha_1+\alpha_2$, $\ney^\alpha=y_1^{\alpha_1} y_2^{\alpha_2}$ when $\ney=(y_1,y_2)\in\R^2$, and
\[
{\binom{\alpha}{\beta}}=\frac{\alpha!}{(\alpha-\beta)!\beta!}=\frac{\alpha_1!}{(\alpha_1-\beta_1)!\beta_1!}\frac{\alpha_2!}{(\alpha_2-\beta_2)!\beta_2!}
\]
when $\beta=(\beta_1,\beta_2)\in\N_0^2$.

Our goal is to use a well-conditioned nodal basis set to produce a Lagrange interpolation polynomial $f_n(\cdot; K): \Omega \to \mathbb{C}$ that regularizes the singular and near-singular integrals over and near to the element $K \ni\nex$, $K \in \mathcal{T}_h$. (The same local polynomial density interpolant $f_n(\cdot;K)$ will be used for all points $\nex \in K$.) A significant difficulty in using local interpolants with~\cref{eq:reg_formula} is that the volume integrals~\cref{eq:reg_formula} depend on the target point $\nex\in K$ explicitly through the Green's function as well as implicitly through $f_n(\cdot; K)$. Clearly, evaluating the volume integral $\mathcal{V}_k[f - f_n(\cdot; K)](\nex)$ over $\Omega$ at all target points $\nex\in K$ for every $K\in\mathcal T_h$ by means of the quadrature rule~\eqref{eq:standard_QR},  would result in a highly inefficient algorithm even using fast methods. These problems are also present for the layer potentials in~\cref{eq:reg_formula}. To mitigate this cost, we express, using the binomial theorem, every \emph{local} polynomial $f_n(\cdot; K)$ in a \emph{fixed, global} basis of normalized monomials, i.e.,
 \begin{equation}\label{eq:mon_interp}
     f_n(\ney; K)=\sum_{|\alpha|\leq n}c_\alpha[f](K)p_\alpha(\ney),\qquad  p_\alpha(\ney) \coloneqq \frac{\ney^\alpha}{\alpha!}.
\end{equation}
By doing so,  $\mathcal{V}_k[f - f_n(\cdot; K)](\nex)$ can be evaluated by means of the volume quadrature  but via the intermediate application of the volume potential to each monomial $p_\alpha$. The application of $\mathcal Q_\Omega$ to $G_k(\nex, \cdot) p_\alpha$ can be precomputed for a given fixed $\Omega$ and then re-used for each new source $f$, merely with different weights $\{c_\alpha[f](K)\}_{|\alpha|\leq n}$---thus leading to an algorithm with offline ($f$-independent) and online ($f$-dependent) components. The offline computations can be accelerated with fast algorithms, while the $f$-dependent weights themselves are \emph{efficiently computable} in the online phase; details on computational costs are provided in~\cref{sec:complexity}.

The coefficients of the interpolant are obtained in a shifted and scaled coordinate system, because, as we explain in detail in \cref{sec:tri_interpolant_constr}, in Lagrange interpolation it is useful (for stability reasons) to solve the multivariate Vandermonde system in coordinates where the triangle fits inside a unit ball (see \Cref{fig:shiftscale0} in \Cref{sec:conditioning}). The interpolation coefficients  $\{c_\alpha[f](K)\}_{|\alpha|\leq n}$ in~\eqref{eq:mon_interp} are then contained in the vector
\begin{equation}\label{eq:ctildevec_system}
    \bold{c}_K = \bold{C}_K \bold{d}_K = \bold{P}_K\widetilde{\bold{V}}^{-1}_K\bold{d}_K\in\C^{q_n},
\end{equation}
where the vector $\bold{d}_K\in\C^{q_n}$ is filled with samples of $f$ at the $q_n$ quadrature nodes on $K$, and  $\widetilde{\bold{V}}_K\in \R^{q_n\times q_n}$ is the Vardermonde matrix in the shifted and scaled coordinate system. Further details on both the linear system $\widetilde{\bold{V}}_K\widetilde{\bold{c}}_K = \bold{d}_K$ associated with~\cref{eq:ctildevec_system} and the change of coordinates map represented by the matrix  $\bold {P}_K\in\R^{q_n\times q_n}$ are provided in \Cref{sec:tri_interpolant_constr} (more precisely, $\bold{P}_K$ is given by~\cref{eq:trans_pol_tri}). Unfortunately, the improvements to conditioning of Vandermonde systems that are obtained by the employed translation-and-scaling method do come with them a potential instability in the change of coordinates map $\bold{P}_K$. A viable strategy to manage this instability is discussed in \Cref{sec:limitations} by using domain decomposition of the domain such that appropriate norms of $\bold{P}_K$ are controlled. In most practical cases and for interpolation degrees $n\in\{0,1,2,3,4\}$, however, a simple re-centering of the domain suffices to manage this instability.

The square shape of the Vandermonde matrix  $\widetilde{\bold{V}}_K$ in~\eqref{eq:ctildevec_system} displays an advantage of utilizing the  Vioreanu-Rokhlin quadrature rule, which is the fact that the number of quadrature points $q_n$ coincides with the cardinality of the monomial basis $\{p_\alpha\}_{|\alpha|\leq n}$, i.e.,
\begin{equation}\label{eq:qn_def}
q_n :=\frac{(n+1)(n+2)}{2}.
\end{equation}

Once the coefficients $\{c_\alpha[f](K)\}_{|\alpha|\leq n}$ in~\cref{eq:mon_interp} are found, the  polynomial in~\eqref{eq:reg_function_ii} is, by linearity,
\begin{equation}\label{eq:mon_reg}
 \Phi_n(\ney; K)=\sum_{|\alpha|\leq n}c_\alpha[f](K)P_{\alpha}(\ney),
\end{equation}
 where the polynomials $\{P_\alpha\}_{|\alpha|\leq n}$ satisfy the inhomogeneous Laplace/Helmholtz equation:
 \begin{equation}\label{eq:nom_inv}
 (\Delta+k^2)P_\alpha=p_\alpha\quad\text{in}\quad\R^2.
\end{equation}
The construction of the polynomial PDE solutions $P_\alpha$ corresponding to the monomial sources $p_\alpha$ is based on the procedures presented in~\cite{anderson2023particular} and those solutions are used here.

With the polynomial $f_n(\cdot;K)$ and polynomial PDE solution $\Phi_n(\cdot;K)$ in hand, we turn to their use in~\cref{eq:reg_formula} and~\cref{eq:reg_formula_W}. The numerical approximation of the regularized volume integrals in~\cref{eq:reg_formula} and~\cref{eq:reg_formula_W} is performed in this work using general-purpose quadratures for smooth functions, both in the online and offline phases of the algorithm. We perform two approximations: first, the contribution from the element $K$ containing the evaluation point is neglected and then smooth quadrature (i.e.\ a quadrature meant for smooth functions that is oblivious to any possible singularities present) on the remaining region is performed. In the description and analysis of each of these approximations that follows, it will therefore be useful to introduce the operators
\begin{equation}\label{eq:breveop}
  \widebreve{\mathcal{V}}_k[f](\nex) \coloneqq \int_\Omega \widebreve{G}_k(\nex, \ney)f(\ney)\,\mathrm{d}\ney\quad\mbox{and}\quad
  \widebreve{\mathcal{W}}_k[\mathbf{f}](\nex) \coloneqq \int_\Omega \nabla_{\ney} \widebreve{G}_k(\nex, \ney) \cdot \mathbf{f}(\ney)\de\ney,
\end{equation}
where $\widebreve{G}_k$ is the punctured Green's function
\begin{equation}\label{eq:punctured_green_function}
\widebreve{G}_k(\nex,\ney) :=\begin{cases} G_k(\nex,\ney), &\nex\in K, \ney\notin K,
\\0,&\nex\in K,\ney\in K.\end{cases}
\end{equation}

In the sequel we denote by $\widebreve{\mathcal{V}}_k^{h,m}[f](\nex)$ and $\widebreve{\mathcal{W}}^{h,m}_k[\bold f](\nex)$ the approximations of these operators on a triangulation $\mathcal{T}_h$ using reference element quadratures that can exactly integrate polynomials of maximum total degree $m \in \mathbb{N}_0$, see e.g.~\cite{Xiao:10,Vioreanu:14}. In fact, these quadratures are independent of $\nex$, implying they are compatible with fast algorithms, and the quadrature over the whole triangular mesh for $\mathcal{V}_k$ can be written in the \Cref{sec:prelim} notation for generic volume quadratures:
\[
  \widebreve{\mathcal{V}}_k^{h,m}[f](\nex) = \mathcal{Q}_\Omega [\widebreve{G}_k(\nex, \cdot)f(\cdot)],
\]
that, by leveraging the expansion~\cref{eq:mon_interp} of $f_n$ in the global basis $\{p_\alpha\}_{|\alpha| \le n}$,  leads to the following numerical method for computing $\mathcal{V}_k[f - f_n]$:
\begin{equation}\label{eq:reg_vol_smooth_quad}
  \begin{split}
    \mathcal{V}_k[f - f_n](\nex) \approx \widebreve{\mathcal{V}}_k[f - f_n](\nex) &\approx \widebreve{\mathcal{V}}_k^{h,m}\left[f - f_n\right](\nex) \\ &= \mathcal Q_\Omega\big[\widebreve G_k(\nex,\cdot)f(\cdot)\big]-\sum_{|\alpha|\leq n}c_\alpha[f](K)\mathcal Q_\Omega\big[ \widebreve G_k(\nex,\cdot)p_\alpha(\cdot)\big].
  \end{split}
\end{equation}
Here, the first approximation commits regularization error, the second approximation commits quadrature error, and the third equality identifies how the volume potential approximation which is to be analyzed can be efficiently computed, for a fixed geometry, in an online ($f$-dependent) / offline ($f$-independent) setting.

It may be useful at this stage to point out that the approximations in \eqref{eq:reg_vol_smooth_quad} can be expected to be of high-quality in view of the fact that $f - f_n(\cdot; K)$ is small (in a sense made precise in \Cref{taylor_lagrange_interpolation_lemma} and \Cref{thm:tri_error_analysis}) over the region near to $\nex$ where $G_k$ is singular or experiences large gradients. Based on the intuition that since $f - f_n(\cdot; K)$ and its derivatives are small near $\nex$ the integrand should be significantly smoother than the kernel $G_k(\cdot,\nex)$ itself, we simply integrate with standard quadrature rules $\mathcal{Q}_\Omega$~\eqref{eq:standard_QR} intended for smooth functions $\phi$ defined over $\Omega$. Precise analysis of this intuition is provided in \Cref{sec:tri_regularization_analysis}, with convergence rates for the approximations in~\cref{eq:reg_vol_smooth_quad} established. (Analogous statements to~\cref{eq:reg_vol_smooth_quad} for $\mathcal{W}_k$ can be written, and the error analysis for these approximations is also given in \Cref{sec:tri_regularization_analysis}.).

Combining~\cref{eq:reg_vol_smooth_quad} with~\cref{eq:quadrature_layerpot} in~\cref{eq:reg_formula} it then follows that $\mathcal V_k [f](\nex)$ for $\nex\in\Omega$, can be approximated as
\begin{equation}\label{eq:approx_form}
\mathcal V_k[f](\nex)\approx \mathcal Q_\Omega\left[\widebreve G_k(\nex,\cdot)f(\cdot)\right]-\sum_{|\alpha|\leq n}c_\alpha[f](K)\big\{\mathcal Q_\Omega\big[ \widebreve G_k(\nex,\cdot)p_\alpha(\cdot)\big] + \mathcal B_\Gamma[P_\alpha,\p_\nu P_\alpha](\nex)+\mu(\nex)P_\alpha(\nex)\big\}.
\end{equation}
As usual, the symbol $\p_{\nu}$ in~\cref{eq:approx_form} denotes the derivative along the exterior unit normal to the boundary $\Gamma$. Clearly, being that the set $\big\{\mathcal Q_\Omega\big[ \widebreve G_k(\nex,\cdot)p_\alpha(\cdot)\big] + \mathcal B_\Gamma[P_\alpha,\p_\nu P_\alpha](\nex)+\mu(\nex)P_\alpha(\nex)\big\}_{|\alpha|\leq n}$ in~\eqref{eq:approx_form} is independent of $f$, it can be precomputed and used for any suitable input density.

\begin{remark}\label{rem:at_boundary} The regularization procedure outlined above for $\nex\in\overline\Omega$ can be directly extended to arbitrary target point locations. Indeed, for target points $\nex\in\R^2\setminus\overline{\Omega}$ lying close to $\Gamma$, the same numerical regularization procedure can be followed, using the Lagrange interpolation polynomial $f_n(\cdot;K)$ of the density $f$ on the quadrature nodes contained in a triangle $K\ni\vv{x}^\star$ where $\vv{x}^* \coloneqq \displaystyle\argminB_{1\leq j\leq N}|\nex-\vv{\xi}_j|$.
\end{remark}

\section{Complexity analysis}\label{sec:complexity}

\begin{table}
\begin{center}
\begin{tabular}{ cc c c }
\toprule
 & Task & Cost \\\cmidrule{2-4}
  \multirow{2}{*}{$f$-dependent}   &$\mathcal Q_\Omega[\widebreve G_k(\vv\xi_j,\cdot)f(\cdot)]$, $j=1,\ldots,N$& $\mathcal O(N{\log N})$  \\
 &$c_\alpha[f](K_\ell)$, $\ell=1,\ldots,L$, $|\alpha|\leq n$ & $\mathcal O(q_n N)$\\
\cmidrule{2-2}\cmidrule{3-3}\cmidrule{4-4}
  \multirow{4}{*}{$f$-independent}&$ \mathcal Q_\Omega\big[ \widebreve G_k(\vv\xi_j,\cdot)p_\alpha(\cdot)\big],$\ \ \  $j=1,\ldots,N,$ $|\alpha|\leq n$ & $\mathcal O(q_n N {\log N})$\\
  &$\mathcal B_\Gamma[P_\alpha,\p_\nu P_\alpha](\vv\xi_j)$,\quad\ \ \ $j=1,\ldots,N,$ $|\alpha|\leq n$ & $\mathcal  O(q_nN {\log N})$    \\
&$\mu(\vv\xi_j)P_\alpha(\vv\xi_j)$,\quad\quad\quad\ \ \  $j=1,\ldots,N,$ $|\alpha|\leq n$& $\mathcal  O(q_nN)$\\
&$\bold{C}_{K_\ell}$,\quad\quad\quad\quad\quad\quad \ \ $\ell=1,\ldots,L$\quad\quad\quad\quad \  & $\mathcal  O(q_n^2N)$\\\bottomrule
\end{tabular}
\end{center}
  \caption{Computational costs for both the online ($f$-dependent) and offline ($f$-independent) components of the method; both the online and offline components are required. Here $N = L q_n$ is the number of volume evaluation points in a mesh with $L$ elements. The cardinality of the polynomial basis is denoted by $q_n$ per~\cref{eq:qn_def}.}\label{tab:complexity-estimates}
\end{table}

The efficiency of the proposed methodology becomes evident from considering the evaluation of the approximation~\cref{eq:approx_form} to $\mathcal V_k[f]$ at all the $N$-numbered volume quadrature nodes $\{\vv\xi_j\}_{j=1}^N$ in the mesh. Table~\ref{tab:complexity-estimates} provides an overview of the computational costs associated with these evaluations, assuming that the FMM is utilized to accelerate the volume and boundary integral computations. (Clearly, the layer and volume potential pre-computations inherent in~\cref{eq:approx_form}, i.e.\ the computations listed in the $f$-independent rows in \Cref{tab:complexity-estimates}, can be computed in an embarrassingly-parallel fashion across the polynomial degree multi-index~$\alpha$.)

The costs in \Cref{tab:complexity-estimates} for evaluating $\mathcal Q_\Omega[\widebreve G_k(\vv\xi_j,\cdot)f(\cdot)]$, $\mathcal Q_\Omega\big[ \widebreve G_k(\vv\xi_j,\cdot)p_\alpha(\cdot)\big]$, and $\mathcal B_\Gamma[P_\alpha,\p_\nu P_\alpha](\vv\xi_j)$ are the costs of FMM-accelerated summation at the $N$ volume target points; the FMM-accelerated general-purpose DIM~\cite{faria2021general} is utilized to evaluate the boundary integrals $\mathcal{B}_\Gamma$, with the boundary discretized using $N_b\propto N^{1/2}$ points. Turning to the evaluation of the coefficients $c_\alpha[f](K)$, $K\in\mathcal T_h$, the computation of the $q_n\times q_n$ matrix $\mathbf{C}_K=\mathbf{P}_K\widetilde{\bold{V}}^{-1}_K$ in~\cref{eq:ctildevec_system} entails the LU-factorization of the Vandermonde matrix $\widetilde{\bold{V}}_K$ for each of the $L$ elements $K\in\mathcal T_h$, which amounts to a cost of $\mathcal{O}(q_n^3 L) = \mathcal{O}(q_n^2 N)$, while evaluating the coordinate map~$\bold{P}_K$ for all of the $L$ elements costs a total of $\mathcal{O}(q_n N)$ in view of the prescription~\cref{eq:trans_pol_tri} for this map. Precomputing $\bold{C}_K$ for every mesh element thus costs $\mathcal{O}(q_n^2 N)$, while evaluation of $c_\alpha[f](K)$, $K\in\mathcal T_h$, by computing $\bold c_K=\bold{C}_K\bold d_K$ in~\eqref{eq:ctildevec_system} costs $\mathcal{O}(q_n N)$. Overall, the operation count estimate shows that, for a given interpolation degree $n\in\N_0$, the VDIM methodology achieves quasi-optimal complexity which is confirmed by the numerical results in~\Cref{sec:numer}.

\section{Interpolant construction}\label{sec:tri_interpolant_constr}
We construct a Lagrange interpolation polynomial $f_n(\cdot; K)$ of total degree at most $n\in\N_0$ that interpolates data at a prescribed nodal set (with $q_n$ defined in~\cref{eq:qn_def}):
\[
\mathcal{I}_n = \{\nex_j: \nex_j\in K\}_{j=1}^{q_n}.
\]
For detailed discussions concerning multivariate polynomial interpolation see~\cite{Sauer:95,Sauer:00,Olver:06}; we discuss the interpolation problem in more depth later but remark here only that given a nodal set $\mathcal{I}_n$ the solvability of the Lagrange interpolation problem (in which case $\mathcal{I}_n$ is called poised) in multiple dimensions is not assured. It suffices that the associated multivariate Vandermonde matrix is nonsingular, a requirement which is achieved, as we will prove, by the choice of the interpolation set $\mathcal{I}_n$ as the mapped Vioreanu-Rokhlin nodes~\cite{Vioreanu:14}, i.e., $\mathcal{I}_n = \vv T(\,\widehat{\mathcal{I}}_n\,)$, $\vv T:\widehat K\to K$, with the reference space nodes $\widehat{\mathcal{I}}_n=\{\widehat \nex_j \}_{j=1}^{q_n} \subset \widehat K$
 being precisely the Vioreanu-Rokhlin nodes.
(The multivariate Vandermonde matrices associated with these nodes on the reference triangle have favorable condition numbers~\cite [Tab.\ 5.1-2]{Vioreanu:14}, at least on a certain basis, suggesting their use here.)

\begin{figure}
\centering
\begin{minipage}[b]{0.58\linewidth}
    \centering
    \begin{tabular}{ c|c|c }
    & $\bol{c}_K$ & $r_K$\\
    \hline
        $K$ obtuse &  $(\vv{v}_2 + \vv{v}_1)/2$ & $\left|\vv{v}_2 - \vv{v}_1\right|/2$ \\
  $K$ acute / right & circumcenter of $K$ & circumradius of $K$
\end{tabular}
\vspace{.5cm}
\end{minipage}%
\hfill
\begin{minipage}[b]{0.35\linewidth}
    \centering
    \includegraphics[width=0.86\linewidth]{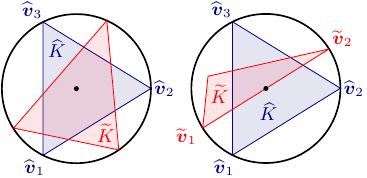}
    \end{minipage}
    \caption{Left: Selection of translation and scale parameters for triangles. Here $\vv{v}_1,\vv{v}_2$, and $\vv{v}_3$ denote the vertices of the triangle $K$ which are numbered in such a way that the side lengths $|\vv{v}_2-\vv{v}_1|$, $|\vv{v}_3-\vv{v}_2|$ and $|\vv{v}_3-\vv{v}_1|$ are in decreasing order. Right: The resulting positions of triangles in the unit ball after translation-and-scaling.}\label{fig:shiftscale}
\end{figure}

As a matter of practicality, one can observe on the one hand that directly interpolating using the physical nodes leads to highly ill-conditioned matrices (in view of the vast differences in scale in the columns), while on the other hand interpolation on the reference simplex $\widehat K$, where conditioning can be considered ideal, is simultaneously not trivially compatible with the use of normalized monomials $\{p_\alpha\}_{|\alpha|\leq n}$ introduced in~\cref{eq:mon_interp} throughout the domain (especially for curved triangles on the boundary) and can also involve a potentially inefficient change of basis calculations. We strike a middle ground by interpolating on a translated and scaled triangle
\[
  \widetilde K\coloneqq {\widetilde{\vv T}}(\widehat K)\quad\mbox{with}\quad {\widetilde{ \vv T}}(\widehat\nex) \coloneqq r_K^{-1}(\vv T(\widehat\nex)-\vv{c}_K),
\]
where $\vv{c}_K$ and $r_K$ are $K$-dependent translation and scale factors, respectively, whose prescription is given in \Cref{fig:shiftscale}. We utilize translation formulae for the purposes of using a global basis. An important consequence for the theory presented in~\cref{sec:conditioning} is that this translation and scaling operation results in the triangle $\widetilde{K}$ that possesses a minimum bounding circle of unit radius. A different translation and scaling procedure, and not with the use of translation formulae to retain a global basis, has been employed in the recent work~\cite{ShenSerkh:22}.

The interpolating polynomial is thus obtained by solving the multivariate Vandermonde system
\begin{equation}\label{eq:vandermonde_shifted}
    \widetilde{\bf{V}}_K\widetilde{\mathbf{c}}_K = \mathbf{d}_K,
\end{equation}
with $\mathbf{d}_K = [f(\nex_1), f(\nex_2), \ldots, f(\nex_{q_n})]^\top\in\C^{q_n}$ the vector of data samples at $\mathcal{I}_n$.
Here $\widetilde{\bf{V}}$ is the multivariate Vandermonde matrix corresponding to a normalized monomial basis set $\{p_\alpha\}_{|\alpha|\leq n}$
and the set $\widetilde{\mathcal{I}}_n\coloneqq \widetilde{\vv T}(\widehat{\mathcal{I}}_n)= \left\{\widetilde \nex_j\coloneqq (\widetilde x_{j,1},\widetilde x_{j,2})\right\}_{j=1}^{q_n}$ of translated-and-scaled points:
\begin{equation}\label{Vandermonde_shiftedscaled_matrix}
    \widetilde{\bf{V}}_K = \begin{pmatrix}
        1 & \widetilde x_{1,1} & \widetilde x_{1,2} & \frac{1}{2} \widetilde x_{1,1}^2 & \widetilde x_{1,1}\widetilde x_{1,2} & \ldots & \frac{1}{(n-1)!} \widetilde x_{1,1}^{n-1} \widetilde x_{1,2} & \frac{1}{n!} \widetilde x_{1,2}^n\\
        1 & \widetilde x_{2,1} & \widetilde x_{2,2} & \frac{1}{2} \widetilde x_{2,1}^2 & \widetilde x_{2,1}\widetilde x_{2,2} & \ldots & \frac{1}{(n-1)!} \widetilde x_{2,1}^{n-1} \widetilde x_{2,2} & \frac{1}{n!} \widetilde x_{2,2}^n\\
        \multicolumn{8}{c}{$\vdots$}\\
        1 & \widetilde x_{q_n,1} & \widetilde x_{q_n,2} & \frac{1}{2} \widetilde x_{q_n,1}^2 & \widetilde x_{q_n,1}\widetilde x_{q_n,2} & \ldots & \frac{1}{(n-1)!} \widetilde x_{q_n,1}^{n-1} \widetilde x_{q_n,2} & \frac{1}{n!} \widetilde x_{q_n,2}^n
    \end{pmatrix}.
\end{equation}
\Cref{sec:conditioning} develops bounds on the condition number of $\widetilde{\bf{V}}_K$ (denoted there simply as $\widetilde{\bf V}$), with the dependence tightly linked to the mesh element quality.

The resulting translated-and-scaled interpolation polynomial by construction satisfies $\widetilde f_n(\widetilde\nex_j;K)=f(\nex_j)$, $j=1,\ldots,q_n$, and can be expressed as
\begin{equation}
    \widetilde{f}_n(\widetilde\nex; K) =\sum_{|\beta| \le n} \widetilde{c}_\beta[f](K) p_{\beta}(\widetilde\nex),\qquad\lf(p_\beta(\widetilde \nex) = \frac{\widetilde \nex^\beta}{\beta!}\rg),
\end{equation}
where the coefficients $\{\widetilde c_\beta[f](K)\}_{|\beta|\le n}$, sorted in (say) lexicographical order, are contained in the vector $\widetilde {\mathbf{c}}_K\in\C^{q_n}$ solution of~\cref{eq:vandermonde_shifted}. Therefore, by setting $\widetilde\nex=r_K^{-1}(\nex-\vv{c}_K)$ above, we get that the sought Lagrange interpolation polynomial is given by
\begin{equation}\label{eq:lagrange_interp_pol}
    f_n(\nex; K) = \sum_{|\beta| \le n} \frac{\widetilde{c}_\beta[f](K)}{r_K^{|\beta|}} \frac{(\nex - \vv c_K)^\beta}{\beta!}.
\end{equation}

The needed separable expansions for $f_n(\cdot;K)$ in~\cref{eq:lagrange_interp_pol} as well as for its associated polynomial particular solution $\Phi_n(\cdot;K)$, can be found by using the binomial formula
\begin{equation}\label{eq:binom}
  (\nex-\vv{c}_K)^{\beta} = \sum_{\alpha:\alpha\leq \beta} {\binom{\beta}{\alpha}}(-\vv{c}_K)^{\beta-\alpha}\nex^\alpha= \sum_{\alpha:\alpha\leq \beta}\beta!\, p_{\beta-\alpha}(-\vv{c}_K)p_\alpha(\nex)
\end{equation}
that results in
 \begin{equation}\label{eq:sep_interp_tri}
  f_n(\nex;K) = \sum_{|\alpha|\leq n} c_\alpha[f](K) p_\alpha(\nex)
\quad\text{and}\quad \Phi_n(\nex;K) = \sum_{|\alpha|\leq n} c_\alpha[f](K)P_\alpha(\nex),
\end{equation}
with the coefficients given by
\begin{equation}\label{eq:trans_pol_tri}
  c_\alpha[f](K) \coloneqq \sum_{\substack{\beta:\alpha \leq \beta,\\|\beta|\le n}}\frac{\widetilde{c}_\beta[f](K)}{r_K^{|\beta|}}p_{\beta-\alpha}(-\bol c_K) .
\end{equation}

The potential $\mathcal{W}_k$ in~\cref{eq:deriv_volpot}, which features a vector-valued function $\mathbf{f}$ in the integrand, requires a vector-valued interpolant $\mathbf{f}_n(\cdot;K)$ that can be constructed component-wise; explicit expressions for $\mathbf{f}_n(\cdot;K)$ and $\boldsymbol{\Phi}_n(\cdot;K)$ are omitted for brevity.

\section{Numerical quadrature accuracy for regularized volume integrals}\label{sec:tri_regularization_analysis}
This section examines the regularizing effect of the density interpolant on the volume potential evaluated at some point $\nex \in K$, $K$ a (curvilinear) triangular element. Its relevance is that part of the method's efficiency owes to the avoidance of the use of (near-)singular volumetric quadrature in the entirety of the domain, as explained above.
Indeed, we develop the arguments that allow us to restrict the use of nearly-singular quadrature (for layer potentials) to a bounded region close to the boundary and still control the error introduced by using generic quadratures meant for smooth functions on regularized volume integrals. The outcome of this analysis is error estimates, given thereafter in~\Cref{thm:tri_error_analysis}. These error estimates are optimal in the sense that the order of convergence established is the same as that for the regularization error, for example $\left|\mathcal{V}_k[f - f_n](\nex) - \widebreve{\mathcal{V}}_k[f - f_n](\nex)\right|$, that is, the order of accuracy obtained by the use of exact evaluation of the nearly-singular volume integrals.
Moreover, they are closely matched by the results from numerical experiments given in \Cref{sec:numer}.
\looseness=-1

\subsection{Main convergence result}

We start by requiring our family of meshes $\mathcal T$ to satisfy the usual quality constraints~\cite{Sauter2010}.
\begin{definition}\label{def:regularuniform}
    A family of triangulations of $\Omega$ denoted by $\mathcal{T} = \{ \mathcal{T}_h \}_{h>0}$ with $h>0$ denoting the maximum element size in $\mathcal{T}_h$ is called \emph{shape-regular} if the shape-regularity constant
    \begin{equation}\label{eq:shaperegulardef}
        \kappa_{\mathcal T} = \sup_{h>0}\max_{K\in\mathcal T_h}\frac{h_K}{\rho_K}
    \end{equation}
    is finite and \emph{quasi-uniform} if the quasi-uniformity constant
    \begin{equation}\label{eq:quasiuniformdef}
        q_\mathcal{T} = \sup_{h>0}\frac{h}{\min_{K\in\mathcal T_h}h_K}
    \end{equation}
    is finite.
\end{definition}

For such mesh families, \Cref{thm:tri_error_analysis} applies for evaluation on subsets of $\overline\Omega$ that are interior to an element with respect to other elements of the mesh, an idea made precise in the next definition.
\begin{definition}\label{def:wellseparatedeval}
  Consider a family of triangulations $\mathcal{T}=\{\mathcal T_h\}_{h>0}$ of $\Omega$, and let $\mathcal S_h=\bigcup_{K\in\mathcal T_h}\p K$. A  family of evaluation point sets $\mathcal{E}=\{\mathcal E_h\}_{h>0}$, $\mathcal E_h\subset\overline{\Omega}$ for all $h>0$, possibly intersecting the boundary $\Gamma=\p\Omega$, is called \emph{well-separated} with respect to $\mathcal{T}$ if
\begin{equation}\label{eq:wellseparated_dist}
    d_{\mathcal T,\mathcal E} \coloneqq \inf_{h>0}\,\,\inf_{\substack{\ney \in \mathcal S_h\setminus\Gamma\\ \nex \in \mathcal{E}_h}} \frac{|\nex - \ney|}{h} > 0.
\end{equation}
\end{definition}

A typical example is evaluation at the point set consisting of (interior) Vioreanu-Rokhlin interpolation/quadrature nodes: the theorem thus estimates the error in the volume potential evaluated at all the interpolation/quadrature points $\left\{{\vv \xi_j}\right\}_{j=1}^N$, $N = Lq_n$. Indeed, any set of interior points on $\widehat{K}$ leads to a well-separated family of evaluation points over the triangulation, as the following proposition shows. Another use-case of the definition would be the evaluation at point-sets laying on the boundary $\Gamma=\partial \Omega$, e.g.\ for the solution of boundary value problems.\enlargethispage*{1ex}

\begin{proposition}\label{rem:interior_VR_nodes}
    Let $\mathcal T=\{\mathcal{T}_h\}_{h>0}$ be a quasi-uniform and shape-regular family of straight triangulations of $\Omega$. The family $\mathcal E=\{\mathcal E_h\}_{h>0}$ with $\mathcal{E}_h$ being a set of interior quadrature nodes over a triangulation $\mathcal T_h\subset\mathcal T$ of $\Omega$, is a well-separated family of evaluation point sets with respect to $\mathcal T$.
\end{proposition}
\begin{proof} Let $\widehat{\mathcal{I}}$ denote the quadrature nodes on $\widehat{K}$ and set $\widehat\delta =\inf_{\widehat\ney\in\p\widehat K,\widehat \nex\in \widehat{\mathcal I}}|\widehat\nex-\widehat \ney|>0$, which is clearly independent of the families $\mathcal T$ and $\mathcal E$. Consider a triangle $K\in\mathcal T_h$ and let $\bol T:\widehat K\to K$ be the affine transformation $\bol T(\widehat\nex)=\bol A\widehat\nex + \bol b$. Then $\mathcal I = \bol T(\widehat{\mathcal I})$ equals the set of quadrature nodes contained in $K$. From the bijectivity of~$\bol T$ and~\cite[Lem.\ 2]{Ciarlet:72} we have
\[
\widehat\delta  \leq |\bol T^{-1}(\nex)-\bol T^{-1}( \ney)|\leq \|\bol A^{-1}\||\nex-\ney|\leq \frac{\widehat h}{\rho_K}|\nex-\ney| \quad\text{for all}\quad \nex\in \mathcal I\subset K,\,\ney\in\p K.
\]
Now, using the fact that
$\rho_{K}^{-1}\leq \kappa_{\mathcal T}/h_K$ and that $h_K^{-1}\leq q_{\mathcal T}/h$, we obtain $\rho_K^{-1}\leq \kappa_{\mathcal T}q_{\mathcal T}/h$,
which by the shape-regularity and quasi-uniform  assumptions on $\mathcal T$ is a finite constant.
It hence follows from the inequalities above that
\[
\widehat\delta  \leq \widehat h k_{\mathcal T}q_{\mathcal T}\frac{|\nex-\ney|}{h} \quad\text{for all}\quad \nex\in \mathcal I\subset K,\,\ney\in\p K.
\]
Therefore, taking infimum over all the points $\nex\in\mathcal E_h$ and $\ney\in \mathcal S_h$, and then  infimum over all the triangulations $\mathcal T_h$, we arrive at
\[
d_{\mathcal T,\mathcal E} \geq \frac{\widehat\delta}{\widehat h}\frac{1}{ k_{\mathcal T}q_{\mathcal T}}>0,
\]
where $d_{\mathcal T,\mathcal E}$ is defined in~\eqref{eq:wellseparated_dist}.
The proof is now complete.
\end{proof}

\begin{remark}
  Similar statements can be easily made for triangulations including curved triangles as well as for families of evaluation point sets $\mathcal{E}$ that lay on $\Gamma$ (and which are yet separated from $\mathcal{S}_h\setminus \Gamma$).
\end{remark}

  A final ingredient for the statement of the theorem is an assumption on the existence and boundedness of Lagrange polynomials on mesh elements in the family.

\begin{assump}\label{assump:lagrange}
  For a given family $\mathcal{T}=\{\mathcal T_h\}_{h>0}$ of triangulations of $\Omega$, the Lagrange interpolation polynomials $\lambda_i(\ney) = \lambda_i(\ney; K_\ell)$ ($1 \le i \le q_n$) exist on each element $K_\ell\in\mathcal T_h\in\mathcal T$ and, further, are uniformly bounded; that is,
  \begin{equation}\label{eq:lagrange_assumption_bound}
    \sup_{\mathcal{T}_h \in \mathcal{T}} \max_{K_\ell \in \mathcal{T}_h} \max_{\ney \in H(K_\ell)} \max_{1 \le i \le q_n} |\lambda_i(\ney; K_\ell)| \le \Lambda(n,\mathcal{T}),
  \end{equation}
  with $H(K_\ell)$ being  the convex hull of the vertices and the set of  interpolation nodes $ \mathcal{I}_n$.
\end{assump}
This is essentially the same concept as the $\Lambda$-poisedness introduced in \Cref{sec:conditioning} for \Cref{eq:conditioning_corr}, for the purposes of analyzing the conditioning of the polynomial interpolation problem. Indeed, Assumption~\ref{assump:lagrange} is satisfied for shape-regular families
$\mathcal{T}$ of straight triangles, a fact that can be seen using elements
of the proof of \Cref{eq:conditioning_corr}: the
estimate~\cref{eq:tildeLambda_poised_bound} together with shape-regularity
providing an upper bound on the quantity $R=h_K / \rho_K$ there. We do not pursue here guarantees of the assumption over
curved triangles.

We are now in a position to state the theorem.

\begin{theorem}\label{thm:tri_error_analysis}
  Let the respective (positive weight) quadrature and interpolation degrees $m, n \in \mathbb{N}_0$, $m > n$, and the wavenumber $k\in \C$ be given and let $\mathcal{T} = \left\{\mathcal{T}_h\right\}_{h>0}$ denote a family of triangulations of the connected domain~$\Omega$  that is shape-regular and quasi-uniform and satisfies \Cref{assump:lagrange}. Let $\mathcal E=\{\mathcal E_h\}_{h>0}$ be a family of well-separated evaluation point sets for $\mathcal T$, and denote by $f_n(\cdot; K)$ (resp.\ $\mathbf{f}_n(\cdot; K)$) the Lagrange polynomial interpolant of $f\in C^{s}(\overline\Omega)$ (resp. $\bold{f}\in [C^{s}(\overline\Omega)]^2$), $s = \max\{n+3,m+1\}$, with interpolation in some $K\in\mathcal T_h$ enforced on the associated interpolation node-set $\mathcal{I}_n\subset K$. For all $\nex \in  K\cap\mathcal{E}_h$ it holds that
\begin{equation}\label{eq:triangle_error_estimate_vr}
    \begin{split}
        \left|\mathcal{V}_k\left[f - f_n(\cdot; K)\right](\nex) - \widebreve{\mathcal{V}}_k^{h,m}\left[f - f_n(\cdot; K)\right](\nex)\right| \le& C^{(1)}_{\mathcal V} h^{n+3}|\log{h}|\\
        & + C^{(2)}_{\mathcal V}
        \begin{cases}
            h^{n+3}, &\quad m > n + 2,\\
            h^{n+3}|\log h|, &\quad m = n + 2,\\
            h^{n+2}, &\quad m = n + 1,\\
        \end{cases}
    \end{split}
  \end{equation}
where $C^{(1)}_{\mathcal V}$ and $C^{(2)}_{\mathcal V}$ are positive constants independent of $h$ but dependent on $\Omega$, $\kappa_\mathcal{T}$, $q_\mathcal{T}$, $f$, $k$, $m$, $n$, $d_{\mathcal{T}, \mathcal{E}}$ and $\Lambda(n, \mathcal{T})$; the constants depend linearly on $\Lambda(n, \mathcal{T})$. Likewise, for the $\mathcal{W}_k$ operator, for all $\nex \in K\cap\mathcal{E}_h$ it holds that
  \begin{equation}\label{eq:triangle_error_estimate_Wk_vr}
      \begin{split}
          \left|\mathcal{W}_k\left[\mathbf{f} - \mathbf{f}_n(\cdot; K)\right](\nex) - \widebreve{\mathcal{W}}_k^{h,m}\left[\mathbf{f} - \mathbf{f}_n(\cdot; K)\right]\vphantom{\widebreve{\mathcal{W}}_k^{h,m}\left[\mathbf{f} - \mathbf{f}_n(\cdot; K)\right]}(\nex)\right| \le& C^{(1)}_{\mathcal W} h^{n+2}\\
          &+ C^{(2)}_{\mathcal W}
          \begin{cases}
              h^{n+2}, &\quad m \ge n + 2,\\
              h^{n+2}|\log h|, &\quad m = n + 1,\\
          \end{cases}
      \end{split}
  \end{equation}
where $C^{(1)}_{\mathcal W}$ and $C^{(2)}_{\mathcal W} $ are positive constants independent of  $h$ but dependent on $\Omega$, $\kappa_\mathcal{T}$, $q_\mathcal{T}$, $f$, $k$, $m$, $n$, $d_{\mathcal{T}, \mathcal{E}}$ and $\Lambda(n, \mathcal{T})$; the constants depend linearly on $\Lambda(n, \mathcal{T})$.
\end{theorem}
\begin{remark}\label{rem:tri_simplify_error_estimates}
The Vioreanu-Rokhlin quadrature rules serve also as interpolation schemes, allowing more concrete error expressions where $m$ and $n$ are related. For interpolation degrees $n\ge 2$, the tabulated Vioreanu-Rokhlin rules satisfy $m \ge n + 2$~\cite[Tab.\ 5.1]{Vioreanu:14}, so the error estimate~\cref{eq:triangle_error_estimate_vr} in such cases simplifies to
\begin{equation}\label{eq:Vk_optimal_error}
    \left|\mathcal{V}_k\left[f - f_n(\cdot; K)\right](\nex) - \widebreve{\mathcal{V}}_k^{h,m}\left[f - f_n(\cdot, K)\right](\nex)\right| \le C_{\mathcal V}h^{n+3}|\log h|, \quad \nex \in K\cap\mathcal{E}_h,
\end{equation}
for some constant $C_{\mathcal V}>0$. Similarly, for $n\geq 2$ the the error estimate~\cref{eq:triangle_error_estimate_Wk_vr} simplifies to
\begin{equation}\label{eq:Wk_optimal_error}
    \left|\mathcal{W}_k\left[\mathbf{f} - \mathbf{f}_n(\cdot; K)\right](\nex) - \widebreve{\mathcal{W}}_k^{h,m}\left[\mathbf{f} - \mathbf{f}_n(\cdot; K)\right](\nex)\right| \le C_{\mathcal W}h^{n+2}, \quad \nex \in K\cap\mathcal{E}_h,
\end{equation}
for some constant $C_{\mathcal W}>0$. Estimates~\cref{eq:triangle_error_estimate_vr} and~\cref{eq:triangle_error_estimate_Wk_vr} yield concrete error estimates when $n < 2$ via the relations~\cite[Tab.\ 5.1]{Vioreanu:14} $m=1$ when $n=0$ and $m=2$ when $n=1$; note that the case $n = 1$ happens to be, comparing to Figures~\ref{fig:tris_splines} and~\ref{fig:tris_polygons}, the case where super-convergence is observed in numerical experiments for the $\mathcal{V}_k$ operator. As a consequence, the analysis that follows shows, for $\mathcal{V}_k$ and $n \ge 2$ or for $\mathcal{W}_k$ and $n \ge 1$, that these estimates for the numerical Vioreanu-Rokhlin quadrature evaluation of the regularized volume integral yield precisely (modulo a $|\log h|$ loss for $\mathcal{W}_k$ with $n = 1$, $m = n + 1$)  the order of accuracy in~\cref{Vk_analytical_regularization_error} and~\cref{Wk_analytical_regularization_error}, respectively, of the regularization error alone---that is, the approximation error introduced by neglecting the regularized weakly-singular integral over $K\ni\nex$.\enlargethispage*{1ex}
\end{remark}

\begin{figure}
    \centering
    \includegraphics[scale=1]{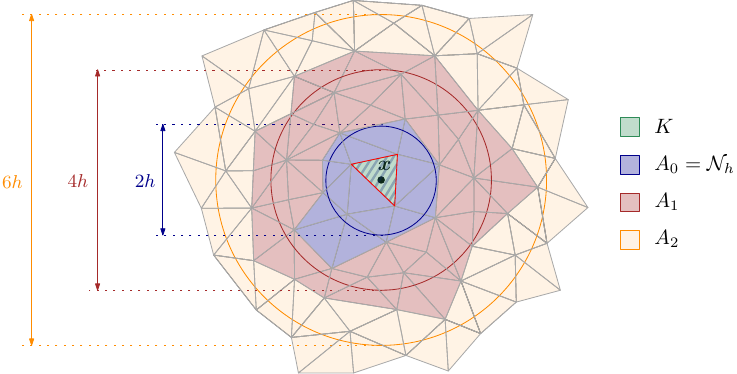}
    \caption{Depiction of the meshed domains employed in the proof of~\Cref{thm:tri_error_analysis}.}
    \label{fig:rings}
\end{figure}

In preparation for proving \Cref{thm:tri_error_analysis}, we note that we must account for contributions from not only singular and near-singular integration but also more distant regions, and for these purposes it is helpful to handle integration regions
separately. We thus define
    \begin{equation}
 \mathcal N_h \coloneqq \mathcal N_h(\nex) = \bigcup\{K_\ell \in \mathcal{T}_h:K_\ell \cap B_{h}(\nex) \ne \emptyset\},
    \end{equation}
    i.e.\ $\mathcal{N}_h$ is an
$\mathcal{O}(h)$ diameter-size neighborhood of $K \ni \nex$ (see \Cref{fig:rings}). Here and in the sequel, $B_{\delta}:=B_\delta(\nex)=\{\ney\in\R^2:|\nex-\ney|<\delta\}$ denotes the Euclidean ball of radius $\delta>0$.

Writing the volume potential $\mathcal{V}_k[f - f_n]$ as
\begin{equation}\label{eq:Vk_decomp_four}
  \mathcal{V}_k\left[ f - f_n(\cdot; K)\right](\nex) = \left( \int_{K} + \int_{\mathcal{N}_h\setminus K} + \int_{\Omega \setminus \mathcal{N}_h} \right) G_k(\nex, \ney) \left[f(\ney) - f_n(\cdot; K)\right]\,\mathrm{d} \ney,
\end{equation}
the analysis that follows (1)~establishes that the contribution of the first integral for the potential is small (the regularization error), (2)~shows that employing standard quadratures for the second integral leads to errors that depend on the interpolation degree and that match the order of accuracy of the regularization error, and (3)~develops estimates for the accuracy of standard quadrature rules when applied to the last integral. Based on the estimates that we develop here we could, for evaluation at a given $\nex$, neglect all elements in the region $\mathcal{N}_h(\nex)$ and retain an optimal order of convergence; in practice, for both simplicity and accuracy, only the contribution from $K$ will be discarded, consistent with our definitions~\cref{eq:punctured_green_function} of the punctured Green's function and~\cref{eq:breveop} of the punctured volume potentials.

Before proving \Cref{thm:tri_error_analysis}, we recall the latter definitions,
and for convenience write them in the form
\begin{equation*}
  \widebreve{\mathcal{V}}_k[f](\nex) = \left(\int_{\Omega \setminus \mathcal N_h} + \int_{\mathcal{N}_h\setminus K} \right) G_k(\ney, \nex)f(\ney)\de\ney
\end{equation*}
and
\begin{equation*}
\widebreve{\mathcal{W}}_k[f](\nex) = \left(\int_{\Omega \setminus \mathcal N_h} + \int_{ \mathcal{N}_h\setminus K} \right) \nabla_{\ney} G_k(\nex, \ney) \cdot \mathbf{f}(\ney)\de\ney.
\end{equation*}
(Also, recall from \Cref{sec:method_overview} that $\widebreve{\mathcal{V}}_k^{h,m}$ and $\widebreve{\mathcal{W}}_k^{h,m}$ denote numerical quadrature approximations to these operators over $\mathcal{T}_h$ with a quadrature rule capable of integrating polynomials of total degree at most $m$.)

\begin{proof}[Proof of \Cref{thm:tri_error_analysis}]
The proof, with exposition primarily restricted to the operator $\mathcal{V}_k$ and relevant differences for the $\mathcal{W}_k$ case mentioned as needed, is divided into in two main parts. The first part proceeds by showing that the error incurred by neglecting the weakly-singular integral over $K$ in $\widebreve{\mathcal{V}}_k$ is an $\mathcal{O}(h^{n+3}\left|\log h\right|)$ quantity as $h\to 0$, and then next that the near-singular volume integrals over each of the elements within the neighborhood $\mathcal{N}_h$ of $\nex$ are each $\mathcal{O}(h^{n+3}\left|\log h\right|)$ quantities so that the integral over $\mathcal{N}_h \setminus K$  approximates the integral over $\mathcal{N}_h$ with errors that behave as $\mathcal{O}(h^{n+3}\left|\log h\right|)$ as $h \to 0$. Then, the second part considers the error from numerical quadrature that arises both from integration on the elements in the near- and far-fields $\Omega\setminus\mathcal{N}_h $. Throughout the proof, ${C}_j$, $j=0,\ldots,6$, will denote positive constants independent of $h$.\enlargethispage*{1ex}

\proofstep{Part 1: contributions from $\mathcal{N}_h$.}
Let $\nex \in K\cap\mathcal{E}_h$ denote a fixed evaluation point. In what follows the Lagrange interpolation error estimate~\cite[Thm.\ 1]{Strang:72} (that holds since the interpolant here is polynomial and thus easily satisfies the required uniformity condition needed there) will prove useful; it provides (cf.\ also~\cite[Thm.\ 2]{Ciarlet:72})
\begin{equation}\label{f_fn_interp_near_estimate}
  \lf|f(\ney) - f_n(\ney; K)\rg| \le {C}_0 h^{n+1}\quad\mbox{for all}\quad \ney \in \mathcal{N}_h(\nex),
\end{equation}
  and, for vector-valued functions the corresponding estimate
\begin{equation}\label{f_fn_interp_near_estimate_vec}
  |\mathbf{f}(\ney) - \mathbf{f}_n(\ney; K)| \le {C}_1 h^{n+1}\quad\mbox{for all}\quad \ney \in \mathcal{N}_h(\nex),
\end{equation}
where ${C}_0$ and ${C}_1$ denote constants independent of $h$ but possibly dependent on the quotient of the diameter of $\mathcal{N}_h$ and $h$.  Note that these estimates hold not merely on $K$ but as well over an $\mathcal{O}(h)$ neighborhood; indeed, the result of~\cite[Thm.\ 2]{Ciarlet:72} can be applied over the entirety of $\mathcal{N}_h$ with interpolation conditions enforced at $\mathcal{I}_n \subset K$.

  The first step of the proof provides bounds on components of the regularized volume potential. Considering each of the integrals over the elements of $\mathcal{N}_h$, we have from a change to polar variables centered at  $\nex\in K\cap\mathcal E_h\subset\mathcal{N}_h$ and the bound~\cref{f_fn_interp_near_estimate}, the estimate
\begin{equation}\label{eq:int_near_sing}
  \bigg|\, \int_{\mathcal{N}_h} G_k(\nex, \ney)\left[ f(\ney) - f_n(\ney; K) \right]\de \ney \,\bigg|\leq \int_{\mathcal{N}_h} \left|G_k(\nex, \ney)\left[ f(\ney) - f_n(\ney; K) \right] \right| \de \ney \le {C}_2 h^{n+3}|\log{h}|,
\end{equation}
which in view of the fact that $K\subset\mathcal{N}_h(\nex)$ implies that
\begin{equation}\label{Vk_analytical_regularization_error}
  \left|\mathcal{V}_k\left[f - f_n(\cdot; K)\right](\nex) - \widebreve{\mathcal{V}}_k\left[f - f_n(\cdot; K)\right](\nex)\right| \le {C}_2 h^{n+3}|\log{h}|.
\end{equation}

In a similar vein, it is easy to see that
\begin{equation}\label{Wk_analytical_regularization_error}
  \left|\mathcal{W}_k\left[\mathbf{f} - \mathbf{f}_n(\cdot; K)\right](\nex) - \widebreve{\mathcal{W}}_k\left[\mathbf{f} - \mathbf{f}_n(\cdot; K)\right](\nex)\right| \le {C}_3 h^{n+2},
\end{equation}
where we used the fact $|\nabla_{\ney} G_k(\nex, \ney)|= \mathcal O(|\nex - \ney|^{-1})$ for $\nex,\ney\in K$.

Similarly, since $\mathcal{N}_h\setminus K\subset\mathcal{N}_h$ we get from the latter inequality of~\cref{eq:int_near_sing} the bound
  \begin{equation}\label{Gk_near_estimate_analytical}
    \begin{split}
      \bigg|\, \int_{\mathcal{N}_h\setminus K} G_k(\nex, \ney)\left[ f(\ney) - f_n(\ney; K) \right]\de \ney \,\bigg| &\le {C}_2 h^{n+3}|\log{h}|
    \end{split}
\end{equation}
and a similar argument for $\mathcal{W}_k$ yields
\begin{equation}\label{nablaGk_near_estimate_analytical}
  \bigg|\, \int_{\mathcal{N}_h \setminus K} \nabla_{\ney} G_k(\nex, \ney)\cdot \left[\mathbf{f}(\ney) - \mathbf{f}_n(\ney; K)\right]\de \ney \,\bigg| \le {C}_3 h^{n+2}.
\end{equation}

  Turning to numerical quadratures and using the rule~\eqref{eq:standard_QR}, for the elements comprising $\mathcal{N}_h\setminus K$, using~\cref{Gk_near_estimate_analytical} and~\cref{f_fn_interp_near_estimate} in conjunction with the triangle inequality, we find
\begin{equation}\label{near_quadrature_estimate}
  \begin{split}
    \bigg|\,\int_{\mathcal{N}_h \setminus K} G_k(\nex, \ney) \left[f(\ney)\;-\right. & \left. f_n(\ney; K)\right]\de\ney - \sum_{\substack{(\vv{\xi}_j,\omega_j) \in \scalebox{1.0}{$\chi$}:\\ \vv\xi_j\in  \mathcal{N}_h \setminus K}} \omega_{j} G_k(\nex,\vv\xi_j) \left[f(\vv\xi_j) - f_n(\vv\xi_j; K)\right] \,\bigg| \le\\
    &\le {C}_2h^{n+3}|\log{h}| + {C}_0 h^{n+1} \max_{\substack{(\vv{\xi}_j,\omega_j) \in \scalebox{1.0}{$\chi$}:\\ \vv\xi_j\in  \mathcal{N}_h \setminus K}} \left|G_k(\nex, \vv\xi_j)\right|\sum_{\substack{(\vv{\xi}_j,\omega_j) \in \scalebox{1.0}{$\chi$}:\\ \vv\xi_j\in  \mathcal{N}_h \setminus K}} \omega_{j}\\
    &\le C_5 h^{n+3}|\log{h}|,
  \end{split}
\end{equation}
  where the inequalities above follow because, firstly, the quadrature rule has positive weights that satisfy $\sum_{\substack{(\vv{\xi}_j,\omega_j) \in \scalebox{1.0}{$\chi$}:\vv\xi_j\in  \mathcal{N}_h \setminus K}}\omega_{j} \lesssim h^2$ and secondly, since $\nex \in \mathcal{E}_h$ lies at a distance from $\mathcal{N}_h\setminus K$ that scales linearly with $h$, $\sup_{\substack{(\vv{\xi}_j,\omega_j) \in \scalebox{1.0}{$\chi$}: \vv\xi_j\in  \mathcal{N}_h \setminus K}} \left|G_k(\nex, \vv\xi_j)\right| \lesssim |\log{h}|$, with an implied constant dependent on $d_{\mathcal{T},{\mathcal{E}}}$ (see \Cref{def:wellseparatedeval}).\enlargethispage*{1ex}

  In a similar vein for the $\widebreve{\mathcal{W}}_k$ kernel, using the triangle inequality and~\cref{nablaGk_near_estimate_analytical}
  together with the estimate~\cref{f_fn_interp_near_estimate_vec} yields the bound
\begin{equation}\label{near_quadrature_estimate_Wk}
  \begin{split}
    \left|\int_{\mathcal{N}_h \setminus K} \nabla_{\ney} G_k(\nex, \ney) \cdot \left[\mathbf{f}(\ney)\;-\right.\right.&\left.\left. \mathbf{f}_n(\ney; K)\right]\de\ney - \sum_{\substack{(\vv{\xi}_j,\omega_j) \in \scalebox{1.0}{$\chi$}:\\ \vv\xi_j\in  \mathcal{N}_h \setminus K}} \omega_{j} \nabla_{\ney} G_k(\nex,\vv\xi_j) \cdot \left[\mathbf{f}(\vv\xi_j) - \mathbf{f}_n(\vv\xi_j; K)\right]\right| \le\\
    &\le C_3 h^{n+2} + C_1 h^{n+1} \max_{\substack{(\vv{\xi}_j,\omega_j) \in \scalebox{1.0}{$\chi$}:\\ \vv\xi_j\in  \mathcal{N}_h \setminus K}} \left|\nabla_{\ney}G_k(\nex, \vv\xi_j)\right|\sum_{\substack{(\vv{\xi}_j,\omega_j) \in \scalebox{1.0}{$\chi$}:\\ \vv\xi_j\in  \mathcal{N}_h \setminus K}}\omega_{j}\\
    &\le {C}_6 h^{n+2},
  \end{split}
\end{equation}
  where the last inequality follows because, recalling $\nex \in \mathcal{E}_h$,  $\sup_{\substack{(\vv{\xi}_j,\omega_j) \in \scalebox{1.0}{$\chi$}:\vv\xi_j\in  \mathcal{N}_h \setminus K}} \left|\nabla_{\ney} G_k(\nex, \vv\xi_j)\right| \lesssim h^{-1}$ with an implied constant dependent on $d_{\mathcal{T},{\mathcal{E}}}$.

\proofstep{Part 2: contributions from $\Omega\setminus\mathcal{N}_h$.}
This part relies on a result concerning the accuracy of ordinary quadratures over  $\Omega \setminus \mathcal N_h$, given thereafter in \Cref{lem:ordquad_convergence_farfield_optimal}. That result in turn requires
sharp estimates on the derivatives of the regularized volume integrands
  \begin{equation}\label{eq:phi1_def}
    \phi_1(\nex, \ney) \coloneqq
    G_k(\nex, \ney) [f(\ney) - f_n(\ney; K)], \quad\quad \nex, \ney \in \widebar{\Omega},
  \end{equation}
  and
  \begin{equation}\label{eq:phi2_def}
    \phi_2(\nex, \ney) \coloneqq
    \nabla_{\ney} G_k(\nex, \ney) \cdot [\mathbf{f}(\ney) - \mathbf{f}_n(\ney; K)], \quad\quad \nex, \ney \in \widebar{\Omega}.
  \end{equation}
Such estimates for the derivatives of $\phi_l(\cdot,\ney)$, $l=1,2$, are given in the following  lemma whose proof (deferred to Section~\ref{proof:taylor_lagrange_interpolation_lemma}) is based on error estimates for multivariate interpolation.

\begin{lemma}\label{taylor_lagrange_interpolation_lemma}
  Take as given the assumptions and setting of \Cref{thm:tri_error_analysis}. In particular, $n\in\N_0$ is the interpolation degree and $m\in\N$, $m > n$, is the degree of exactness of the quadrature rule. Letting  $\alpha\in\N_0^2$ be a multi-index satisfying $|\alpha| = m + 1$ and for $\phi_1$ given in~\eqref{eq:phi1_def}, the following estimate holds for all $\delta\in(0,\min\{1,(\operatorname{diam}\Omega)^{1/(n-m)}\}]$   and any $\nex \in K$:
      \begin{equation}\label{eq:interpolation_estimate}
         \left|D_{\ney}^{\alpha} \phi_1(\nex, \ney)\right| \le
 C_{\mathscr{T}}\delta^{n-m}+  C_{\mathscr{L}} \Lambda \sum_{\gamma:|\gamma| \le n} h^{n + 1 - |\gamma|} \delta^{|\gamma| - m - 1},\quad \ney \in\Omega\setminus B_\delta(\nex),
    \end{equation}
where $C_{\mathscr{T}} = C_{\mathscr{T}}(m,n, k, f,\Omega) $ and $C_{\mathscr{L}}
    = C_{\mathscr{L}}(m,n, k, \kappa_\mathcal{T}, q_\mathcal{T}, f, \Omega)$ denote positive constants both  independent of
    $\nex$, $\delta$, $K$, $\mathcal{T}$, $h$ and $\Lambda$ ($\Lambda$ defined in \Cref{assump:lagrange}).

    Additionally, for $\phi_2$ given in~\eqref{eq:phi2_def}, for the same $\alpha\in\N_0^2$ and  for all $\delta\in(0,\min\{1,(\operatorname{diam}\Omega)^{1/(n-m)}\}]$  and any $\nex \in K$:
\begin{equation}\label{eq:interpolation_estimate_vec}
        \left|D_{\ney}^{\alpha} \phi_2(\nex, \ney)\right| \le  C_{\mathscr{T}}' \delta^{n-m-1} + C_{\mathscr{L}}' \Lambda \sum_{\gamma:|\gamma| \le n} h^{n+1 - |\gamma|} \delta^{|\gamma| - m - 2}, \quad \ney \in \Omega\setminus B_\delta(\nex),
    \end{equation}
    where $C_{\mathscr{T}}' = C_{\mathscr{T}}'(m,n, k, f,\Omega)$ and
    $C_{\mathscr{L}}' = C_{\mathscr{L}}'(m, n, k, \kappa_\mathcal{T}, q_\mathcal{T}, f, \Omega)$  again denote positive
    constants both independent of $\nex$, $\delta$, $K$, $\mathcal{T}$, $h$ and
    $\Lambda$.
\end{lemma}

The estimates given in \Cref{taylor_lagrange_interpolation_lemma} for the derivatives of the integrands $\phi_1(\cdot,\nex)$ for $\mathcal{V}_k[f - f_n(\cdot; K)]$ (resp.\ $\phi_2(\cdot, \nex)$ for $\mathcal{W}_k[\mathbf{f} - \mathbf{f}_n(\cdot; K)]$) are on themselves insufficient to obtain
convergence rates of higher order than $\mathcal{O}(h^{n+1})$ (resp.\ $\mathcal{O}(h^n)$) in connection with classical error
estimates~\cite[Sec.\ 7.4]{IsaacsonKeller}: while the integrand is globally smooth in
$\Omega \setminus \mathcal{N}_h$, the growth of derivatives $\left|D_{\ney}^{\alpha} \phi_j(\nex,
\cdot) \right|$ ($j = 1, 2$) in the \Cref{taylor_lagrange_interpolation_lemma} estimates, for derivative orders $|\alpha|=m+1> n + 1$, exactly
offsets any possible gain in order of accuracy from use of a higher-order quadrature rule.
The following lemma, whose proof is given in Section~\ref{proof:ordquad_convergence_farfield_optimal}, addresses this issue, recovering near-optimal convergence rates (optimal up to a possible $|\log h|$ multiplicative factor in some limited cases). It does so by developing sharp bounds for the  quadrature error over  $\Omega\setminus\mathcal{N}_h$.

\begin{lemma}\label{lem:ordquad_convergence_farfield_optimal}
  Take as given the assumptions and setting of \Cref{thm:tri_error_analysis}. In particular, $n\in\N_0$ is the interpolation degree and $m\in\N$, $m > n$, is the degree of exactness of the reference quadrature rule leading to the global composite quadrature rule $\mathcal{Q}_\Omega$  in~\eqref{eq:globnodesweights}-\eqref{eq:standard_QR}.
For all $\nex \in \mathcal{E}_h$ and for sufficiently small $h>0$, it holds that
\begin{equation}\label{far_quadrature_optimal_estimate_V}
  \bigg|\, \int_{ \Omega \setminus \mathcal N_h} \phi_1(\nex, \ney)\de \ney - \mathcal{Q}_{\Omega \setminus \mathcal{N}_h}[\phi_1(\nex, \cdot)] \,\bigg| \le \mathcal{C}_1 \Lambda
    \begin{cases}
        h^{n + 3}, &\quad m > n + 2,\\
        h^{n+3}|\log h|, &\quad m = n + 2,\\
        h^{n+2}, &\quad m = n + 1,
    \end{cases}
\end{equation}
and
\begin{equation}\label{far_quadrature_optimal_estimate_W}
  \bigg|\, \int_{ \Omega \setminus \mathcal{N}_h} \phi_2(\nex, \ney)\de \ney - \mathcal{Q}_{\Omega \setminus \mathcal{N}_h}[\phi_2(\nex, \cdot)] \,\bigg| \le \mathcal{C}_2 \Lambda
        \begin{cases}
            h^{n+2}, &\quad m \ge n + 2,\\
            h^{n+2}\left|\log h\right|, &\quad m = n + 1,
        \end{cases}
\end{equation}
    with positive constants $\mathcal{C}_j = \mathcal{C}_j(m,n, k, \kappa_\mathcal{T}, q_\mathcal{T}, f, \Omega)$ ($j = 1,2$) independent of $h$ and $\Lambda$ ($\Lambda$ defined in \Cref{assump:lagrange}).
\end{lemma}

\proofstep{Proof completion.}
Collecting~\cref{Vk_analytical_regularization_error}, \cref{near_quadrature_estimate}, and~\cref{far_quadrature_optimal_estimate_V} from \Cref{lem:ordquad_convergence_farfield_optimal}, we conclude that for all $\nex \in  K\cap\mathcal{E}_h$ it holds that
\begin{equation*}
    \begin{split}
      \left|\mathcal{V}_k\left[f - f_n(\cdot; K)\right](\nex) - \widebreve{\mathcal{V}}_k^{h,m}\left[f - f_n(\cdot; K)\right](\nex)\right| &\le \left|\mathcal{V}_k\left[f - f_n(\cdot; K)\right](\nex) - \widebreve{\mathcal{V}}_k\left[f - f_n(\cdot; K)\right](\nex)\right|\\
      &\quad +\left|\widebreve{\mathcal{V}}_k\left[f - f_n(\cdot; K)\right](\nex) - \widebreve{\mathcal{V}}_k^{h,m}\left[f - f_n(\cdot; K)\right](\nex)\right|\\
      &\le {C}^{(1)}_{\mathcal V} h^{n+3}|\log{h}| + C^{(2)}_{\mathcal{V}}
        \begin{cases}
            h^{n+3}, &\quad m > n + 2,\\
            h^{n+3}|\log h|, &\quad m = n + 2,\\
            h^{m+1}, &\quad m = n + 1.
        \end{cases} \quad
    \end{split}
\end{equation*}
Similarly, for $\widebreve{\mathcal{W}}_k$ by collecting~\cref{Wk_analytical_regularization_error}, \cref{near_quadrature_estimate_Wk}, and~\cref{far_quadrature_optimal_estimate_W} from \Cref{lem:ordquad_convergence_farfield_optimal} we find that for all $\nex \in  K\cap\mathcal{E}_h$ it holds that
  \begin{equation*}
    \begin{split}
      \left|\mathcal{W}_k\left[\mathbf{f} - \mathbf{f}_n(\cdot; K)\right](\nex) - \widebreve{\mathcal{W}}_k^{h,m}\left[\mathbf{f} - \mathbf{f}_n(\cdot; K)\right](\nex)\right| &\le \left|\mathcal{W}_k\left[\mathbf{f} - \mathbf{f}_n(\cdot; K)\right](\nex) - \widebreve{\mathcal{W}}_k\left[\mathbf{f} - \mathbf{f}_n(\cdot; K)\right](\nex)\right|\\
      &\quad +\left|\widebreve{\mathcal{W}}_k\left[\mathbf{f} - \mathbf{f}_n(\cdot; K)\right](\nex) - \widebreve{\mathcal{W}}_k^{h,m}\left[\mathbf{f} - \mathbf{f}_n(\cdot; K)\right](\nex)\right|\\
      &\le {C}^{(1)}_{\mathcal W} h^{n+2}  + C^{(2)}_{\mathcal{W}}
      \begin{cases}
              h^{n+2}, &\quad m \ge n + 2,\\
              h^{n+2}|\log h|, &\quad m = n + 1.
      \end{cases}
    \end{split}
\end{equation*}
The proof is complete.
\end{proof}

\begin{remark}
    Weaker versions of \Cref{thm:tri_error_analysis} can be easily established using a corresponding version of \Cref{taylor_lagrange_interpolation_lemma} with quadrature order $m \le n$ (and, indeed, with no need for \Cref{lem:ordquad_convergence_farfield_optimal}). However, since the resulting quadrature error of $\mathcal{O}(h^{m+1})$ would not saturate the order of accuracy of the regularization error $\left|\mathcal{V}_k[f - f_n] - \widebreve{\mathcal{V}}_k[f - f_n]\right|$ that is established in \Cref{thm:tri_error_analysis}, the use of such low-powered quadratures is disfavored and so not studied in this paper.
\end{remark}

\begin{remark}\label{rem:composite_quads}
   It may be useful for intuition to grasp the basic implication of Lemma~\ref{lem:ordquad_convergence_farfield_optimal} in a simplified context, ignoring all regularization techniques employed in this paper. The proof of Lemma~\ref{lem:ordquad_convergence_farfield_optimal} (see \Cref{proof:ordquad_convergence_farfield_optimal}) in essence establishes that composite quadrature of $\log|\ney|$ over $B_1(0)$, when ignoring the integral over the mesh element that contains the origin and when using a rule of exactness degree $m\ge2$, results in errors of $\mathcal{O}(h^2|\log h|)$; the (pessimistic) classical error estimate~\cite[Sec.\ 7.4]{IsaacsonKeller} does not predict convergence. The techniques involve use of bounds on concentric rings (see \Cref{fig:rings}) surrounding the singularity over which classical estimates are repeatedly used, implementing the intuition that the integrand is `well-behaved' on `most' mesh elements. Such \emph{composite} quadrature error estimates over unstructured meshes in the presence of an isolated singularity appear to be novel. We find also low-order accurate estimates on quadrature over uniform structured grids in the presence of weak singularities~\cite{vainikko2006multidimensional} (some of the proof techniques appear similar to those that we independently developed and employed); all higher-order estimates provided there (which are, additionally, restricted to parallelepipeds) make essential use of composite quadratures on graded meshes. We also find interesting connection to a different body of work in the broader context of \emph{non-composite} quadrature in the presence of isolated singularities, much of it involving Rabinowitz~\cite{lubinsky1984rates, rabinowitz1986rates,rabinowitz1967gaussian,xiang2012convergence,davis1965ignoring}, for Gauss or Clenshaw-Curtis type quadratures.
\end{remark}

\subsection{Proof of \Cref{taylor_lagrange_interpolation_lemma}}
\label{proof:taylor_lagrange_interpolation_lemma}

  Our goal is to find suitable bounds for $D^{\alpha}_{\ney}\phi_1(\nex,\ney)$, which, applying the product rule, can be expressed as
  \begin{equation}\label{leibniz_rule_phi1}
    D^{\alpha}_{\ney}\phi_1(\nex,\ney) = \sum_{\beta: \beta \leq \alpha} \binom{\alpha}{\beta} \left(D^\beta_{\ney} G_k(\nex,\ney)\right)\left(D^{\alpha-\beta} (f(\ney)-f_n(\ney;K))\right).
  \end{equation}
  We start by noticing that the partial derivatives of $G_k$ in~\eqref{eq:Green_Function} satisfy the bounds
    \begin{equation}\label{eq:Gk_bounds}
    \left|D_{\ney}^\eta G_k(\nex, \ney)\right| \le
      \begin{cases}
          C_G |\nex - \ney|^{-|\eta|}, & |\eta| > 0,\\
          C_G(1 +  \left|\log |\nex - \ney|\right|), & |\eta| = 0,
      \end{cases} \quad\quad C_{G} = C_G(k,\eta) > 0, \quad \nex\neq\ney,
  \end{equation}
  that  follow from the properties of the Hankel/logarithm function and induction~\cite[\S 10.6, 10.7]{NIST:DLMF}.

  Our proof strategy to bound the other factors in~\eqref{leibniz_rule_phi1} is based on expressing
  \begin{equation}\label{eq:decomp_taylor_lagrange}
      f(\ney) - f_n(\ney; K) = R_{\nex, n}^{(1)}[f](\ney) + R_{\nex, n}^{(2)}[f](\ney),
  \end{equation}
  where, denoting by $\mathscr{T}_{\nex, n}[f](\ney)$ the $n$th total degree Taylor
    polynomial of $f$ centered at $\nex$, we define the two remainder functions
  \begin{equation}\label{eq:decomp_taylor}
      R_{\nex, n}^{(1)}[f](\ney) = f(\ney) - \mathscr{T}_{\nex,n}[f](\ney),\quad\mbox{and}\quad  R_{\nex, n}^{(2)}[f](\ney) = \mathscr{T}_{\nex,n}[f](\ney) - f_n(\ney; K).
  \end{equation}
  The term
    $R_{\nex, n}^{(1)}[f]$ can be simply expressed in terms of the Taylor
    remainder, while $R_{\nex, n}^{(2)}[f]$ will be rewritten
    using a certain `multipoint' Taylor formula for Lagrange polynomials~\cite{Ciarlet:72}.

    As expected, for $\eta\in\N_0^2$, the Taylor remainder $R_{\nex, n}^{(1)}[f](\ney)$ satisfies
     \begin{equation}\label{eq:taylor_error_bounds}
      \left|  D^{\eta}_{\ney} R^{(1)}_{\nex, n}[f](\ney)\right| \le
      \begin{cases}
          C_{R^{(1)}} |\ney - \nex|^{n + 1 - |\eta|},  & n + 1 > |\eta|, \\
          C_{R^{(1)}}, &n + 1 \le |\eta| \leq m + 1,
      \end{cases} \quad\quad C_{R^{(1)}} = C_{R^{(1)}}(\eta, n, f) > 0.
  \end{equation}
This estimate can be proven for $|\eta|<n+1$   by first
   leveraging the identity
   $ D^{\eta} R^{(1)}_{\nex, n}[f] =R^{(1)}_{\nex, n-|\eta|}[D^\eta f]
 $    which follows directly from the fact that $D^\eta \mathscr{T}_{\nex,n}[f]$ is the  $(n-|\eta|)$th-degree Taylor polynomial  of $D^\eta f$. From the integral form of the Taylor remainder $R^{(1)}_{\nex, n-|\eta|}[D^\eta f]$~\cite[Sec.1.1]{taylor2013partial}, we obtain
   \begin{equation}\label{eq:taylor_error}
      \begin{split}
          D^{\eta}_{\ney} R^{(1)}_{\nex, n}[f](\ney)&= (n+1-|\eta|)\sum_{\gamma:|\gamma| = n+1-|\eta|}  \frac{(\ney - \nex)^\gamma}{\gamma!}\int_0^1 (1 - t)^n D^{\gamma+\eta} f(\nex + t(\ney - \nex))\de t.
      \end{split}
  \end{equation}
(For the sake of simplicity of presentation, we have assumed here that $\nex$ and $\ney$ can be joined by a straight line, an assumption that is fulfilled for convex domains $\Omega$. However, a simple generalization of the remainder formula \cite[Thm.A1]{driveranalysis}, which is valid for path-connected domains, can be used instead.)
Then, from the bounds
$$\left|\int_0^1 (1 - t)^n  D^{\gamma+\eta} f(\nex + t(\ney - \nex))\de t \right|\leq \frac{\|f\|_{C^{|\gamma|+|\eta|}(\overline\Omega)}}{n+1}$$
and
\begin{equation}\label{eq:bound_mon}
|(\ney-\nex)^{\gamma}|=\prod_{j=1}^2|y_j-x_j|^{\gamma_j} \leq \prod_{j=1}^2|\ney-\nex|^{\gamma_j} = |\ney-\nex|^{\sum_{j=1}^2\gamma_j} = |\ney-\nex|^{|\gamma|} ,
\end{equation}
we readily arrive at
\begin{equation}
      \left| D^{\eta}_{\ney} R^{(1)}_{\nex, n}[f](\ney)\right|\leq \left\{\|f\|_{C^{n+1}(\overline\Omega)}\left(\frac{n+1-|\eta|}{n+1}\right)
          \sum_{\gamma:|\gamma| = n + 1-|\eta|} \frac{1}{\gamma!} \right\} |\ney-\nex|^{n+1-|\eta|},
  \end{equation}
 which gives an expression for the constant $C_{R^{(1)}}$ in the case $|\eta|\leq n+1$. For $n+1<|\eta|\leq m+1$, on the other hand, the bound follows directly from the fact that $D^\eta \mathscr{T}_{\nex,n}[f]=0$, which allows us to select $C_{R^{(1)}}=\|f\|_{C^{|\eta|}(\overline\Omega)}$.

  Remembering that $\alpha\in\N_0^2$ is a multi-index satisfying $|\alpha| = m + 1>n+1$ and utilizing the bounds~\eqref{eq:taylor_error_bounds} with $\eta=\beta$ in conjunction with~\eqref{eq:Gk_bounds} with $\eta = \alpha -
  \beta$ , we have
  \begin{equation}\label{eq:prod_Taylor}
      \begin{split}
      \left| D^{\alpha}_{\ney} \left(G_k(\nex, \ney)R_{\nex, n}^{(1)}[f](\ney) \right)\right| =& \left| \sum_{\beta: \beta \le \alpha} \binom{\alpha}{\beta} \left(D_{\ney}^{\alpha-\beta} G_k(\nex, \ney)\right) \left(D^{\beta}_{\ney} R_{\nex, n}^{(1)}[f](\ney)\right)\right| \\
          \le& \sum_{\beta: |\beta| \le n+1} \binom{\alpha}{\beta} C_G |\,|\nex - \ney|^{-|\alpha-\beta|} C_{R^{(1)}} |\nex - \ney|^{n+1-|\beta|}\\
&\phantom{\le} +\sum_{\substack{\beta: \beta < \alpha\\ n+1<|\beta|}} \binom{\alpha}{\beta} C_G |\nex - \ney|^{-|\alpha-\beta|} C_{R^{(1)}}\\
&\phantom{\le} +C_{R^{(1)}}C_G(1+|\log|\nex-\ney||)\qquad (\nex\neq\ney).
      \end{split}
  \end{equation}

To bound the term with the logarithm, we first note that $\log\delta\leq\log|\nex-\ney|\leq |\nex-\ney|\leq \operatorname{diam}\Omega$, for all $\nex,\ney\in\Omega$, $\delta\leq |\nex-\ney|$. If $|\nex-\ney|\geq 1$, the fact that  $\delta\leq(\operatorname{diam}\Omega)^{1/(n-m)}$ implies $\delta^{n-m}\geq \operatorname{diam}\Omega\geq |\nex-\ney|\geq \log|\nex-\ney|=|\log|\nex-\ney||$ since $m > n$. On the other hand, if $|\nex-\ney|<1$,  it follows from $\delta\leq 1$ and $m>n$ that $\delta^{n-m}\geq \delta^{-1}\geq\log(\delta^{-1})\geq -\log|\nex-\ney|=|\log|\nex-\ney||$. From these inequalities we readily conclude that
$$
C_{R^{(1)}}C_G(1+|\log|\nex-\ney||)\lesssim (1+\delta^{n-m})\lesssim\delta^{n-m},
$$
for a suitable implied constant.

To bound the other two terms in~\eqref{eq:prod_Taylor}, we note that since $|\nex-\ney|\geq\delta$, for  $m>n$ it holds that $|\nex - \ney|^{n-m} \leq  \delta^{n-m}$. Then, for $|\alpha|=m+1$ and for suitable implied constants we  have, using the relation $n + 1 - |\alpha| = n - m$,
$$\sum_{\beta: |\beta| \le n+1} \binom{\alpha}{\beta} C_G |\nex - \ney|^{-|\alpha-\beta|} C_{R^{(1)}} |\nex - \ney|^{n+1-|\beta|}\lesssim |\nex-\ney|^{n-m}\lesssim \delta^{n-m},$$ and
$$\sum_{\substack{\beta: \beta < \alpha\\ n+1<|\beta|}} \binom{\alpha}{\beta} C_G |\nex - \ney|^{-|\alpha-\beta|} C_{R^{(1)}}\lesssim \sum_{j=n+2}^{m+1} \delta^{-m-1+j}\lesssim \delta^{n+1-m}\lesssim \delta^{n-m},$$ where in the latter the fact that $\delta\leq 1$ was also used.

Therefore, using the above relations we conclude that~\eqref{eq:prod_Taylor} simplifies   to
    \begin{equation}\label{eq:taylor_error_final}
      \left| D^{\alpha}_{\ney} \left(G_k(\nex, \ney)R_{\nex, n}^{(1)}[f](\ney) \right)\right|  \le
      C_{\mathscr{T}}\delta^{n-m} ,\quad \ney\in \Omega\setminus B_\delta(\nex), \ \delta\in(0,\min\{1,(\operatorname{diam}\Omega)^{1/(n-m)}\}],
  \end{equation}
  where $ C_{\mathscr{T}} =  C_{\mathscr{T}}(m,n, k, f,\Omega)$.  This completes the estimate concerning $R^{(1)}_{\nex, n}[f]$.

  Proceeding now with the case of $R_{\nex, n}^{(2)}[f]$ we will rewrite it
  using~\cite[Thm.\ 1]{Ciarlet:72}; we note that the hypotheses of
  that theorem are (i)~That the interpolation set $\mathcal{I}_n$ is poised
  (termed `$k$-unisolvent' there) and (ii)~That $H(K)$ (defined in
  \Cref{assump:lagrange}) is star-shaped with respect to every element in
  $\mathcal{I}_n$ (termed `$K$ is a $\Sigma$-admissible set' there). Both are assured by the present hypotheses of \Cref{assump:lagrange}; note that by construction $H(K)$ is always star-shaped with respect to each interpolation node $\nex_i \in \mathcal{I}_n$. Thus, proceeding first with estimates for $R_{\nex, n}^{(2)}[f]$, the theorem provides  for $\gamma\in\N_0^2$ that the derivatives of the Lagrange interpolation error can be expressed as
    \begin{equation}\label{eq:multipoint_expr}
      {D}^\gamma f(\ney) - {D}^\gamma f_n(\ney; K) = -\frac{1}{(n+1)!} \sum_{i=1}^{q_n} \left\{{D}^{(n+1)} f(\eta_i(\ney);\nex_i-\ney)\right\} D^\gamma \lambda_i(\ney), \quad \ney \in H(K),
  \end{equation}  with $\eta_i(\ney) = \theta_i \ney + (1 - \theta_i) \nex_i$ for some $0 <
  \theta_i < 1$ and $\nex_i\in\mathcal I_n$. Here ${D}^{(n+1)} f(\eta_i(\ney);\nex_i-\ney)$ denotes the $(n+1)$th-order total differential of $f$ at $\eta_i(\ney)$~\cite[Ch.7]{marsden1993elementary}; i.e.,
  \begin{equation}\label{eq:total_differential}
  {D}^{(n+1)} f(\eta_i(\ney);\nex_i-\ney) = \sum_{j_1=1}^2\cdots\sum_{j_{n+1}=1}^2 \frac{\p^{n+1} f}{\p x_{j_1}\cdots\p x_{j_{n+1}}}(\eta_i(\ney))h_{i_1}\cdots h_{j_{n+1}},\quad (h_1,h_2)=\nex_i-\ney.
  \end{equation}
  (In the event that $H(K) \not\subset \Omega$ --- that is, if $K$ is a curved triangle with concave boundary and so $f$ is not defined for all possible $\eta_i(\ney)$ --- we identify $f$ here and in what follows with its Sobolev extension to $\R^2$ assured by the continuous linear extension operator $\mathscr{E}_s: H^s(\Omega) \to H^s(\R^2)$, $s \in \N_0$, provided by~\cite[\S 6, Thm.\ 5]{stein1970singular} for Lipschitz $\Omega$. This, together with the Sobolev lemma~\cite{folland} accounts for the assumed additional technical need of $f \in C^{n+3}(\widebar{\Omega})$ so that the extension is also an element of $C^{n+1}(H(K))$ with norm bounded by $\left\|\mathscr{E}_{n+3}\right\|$ times the norm of $f \in C^{n+3}(\overline{\Omega})$.)

  Now, since both $\mathscr{T}_{\nex,n}[f](\ney)$ and $f_n(\ney; K)$ are polynomials of degree $n$ in $\ney$, we can write
    \begin{equation}\label{eq:remainder_2_multipoint_power_series}
          R_{\nex, n}^{(2)}[f](\ney) = \mathscr{T}_{\nex,n}[f](\ney) - f_n(\ney; K)       = \sum_{\gamma:|\gamma| \le n} c_\gamma (\ney - \nex)^\gamma,
  \end{equation}
  where the coefficients $\{c_\gamma\}_{|\gamma|\leq n}$ that measure error in approximating the Taylor coefficients are given by
  \[
  c_\gamma = \frac{1}{\gamma!} D^\gamma R_{\nex, n}^{(2)}[f](\nex) =\frac{1}{\gamma!} \left(D^\gamma \mathscr{T}_{\nex,n}[f](\nex)  -D^\gamma f_n(\nex;K)\right)=\frac{1}{\gamma!} \left(D^\gamma f(\nex)  -D^\gamma f_n(\nex;K)\right).
  \]
  Note that we have used here the fact that $\mathscr{T}_{\nex,n}[f]$ is the $n$th-degree Taylor polynomial of $f$ at $\nex$ so that $D^\gamma \mathscr{T}_{\nex,n}[f](\nex)=D^\gamma f(\nex)$. We then obtain from~\eqref{eq:multipoint_expr} that
\begin{equation}\label{eq:cgamma_def}
    c_\gamma = -\frac{1}{\gamma!(n+1)!} \sum_{i=1}^{q_n} \left\{{D}^{(n+1)} f(\eta_i(\nex);\nex_i-\nex)\right\} D^{\gamma} \lambda_i(\nex).
\end{equation}
We next provide bounds on the  coefficients $\{c_\gamma\}_{|\gamma|\leq n}$. In view of~\eqref{eq:total_differential} we first observe that
  \begin{equation}\label{eq:deriv_bounds}
      \left|{D}^{(n+1)} f(\eta_i(\ney);\nex_i-\nex)\right| \le M_{n+1} h^{n+1},\quad i=1,\ldots,q_n,
  \end{equation}
  with $M_{n+1} = 2^{n+1} \left\|\mathscr{E}_{n+3}\right\|\|f\|_{C^{n+3}(\overline\Omega)}$. In order to bound the derivatives of the Lagrange polynomials $\lambda_i$, we use the following
  Markov-like inequality: for the convex hull $H(K)$, $K \in \mathcal{T}_h$, and a polynomial $p$ of total degree at
  most $n$ satisfying $\sup_{\nex \in H(K)} |p(\nex)| \le \Lambda$, it holds that
  \begin{equation}\label{eq:markov_iterated}
    \sup_{\nex \in H(K)} |D^\gamma p(\nex)| \ \le\ (4n^2\kappa_\mathcal{T}q_{\mathcal T})^{|\gamma|}\,h^{-|\gamma|} \sup_{\nex \in H(K)} |p(\nex)| \ \le\ C_1 \kappa_\mathcal{T}^{|\gamma|}q_{\mathcal T}^{|\gamma|} h^{-|\gamma|} \Lambda.
  \end{equation}
(where $\kappa_\mathcal{T}$ and $q_\mathcal{T}$ are the shape-regularity and quasi-uniformity constants of the mesh family). This inequality stems from iterating $|\gamma|$ times the inequality chain~\cite[eq.~3.2]{wilhelmsen:74}
\[
  \sup_{\nex \in H(K)}|\partial_{j} p(\nex)| \leq \sup_{\nex \in H(K)} \| \nabla p(\nex)\|_{\infty}
  \leq \sup_{\nex \in H(K)} \| \nabla p(\nex)\|_2 \leq \frac{4n^2}{\omega_{H(K)}}\sup_{\nex \in H(K)}|p(\nex)| \qquad (j=1,2),
\]
in which the \emph{width} $\omega_{H(K)}$ of a compact $H(K) \subset \R^2$, defined~\cite{wilhelmsen:74} as the smallest distance between two parallel lines touching $\partial H(K)$ without intersecting its interior, here verifies $h/(\kappa_\mathcal{T}q_{\mathcal T})\leq\omega_{H(K)}\leq h$ on account of the fact that two parallel lines touching $\partial H(K)$ are separated by at least the diameter of the largest inscribed circle in $H(K)$ and are separated by at most the diameter $h$ of the triangle (see also \Cref{def:regularuniform}, from which it follows that $h / q_\mathcal{T}$ provides a lower bound on the mesh element size in the preceding argument).
The inequalities~\cref{eq:deriv_bounds} and~\cref{eq:markov_iterated} together imply that $c_\gamma$ in~\cref{eq:cgamma_def} can be bounded as
\begin{equation}\label{taylor_lagrange_error_coeff_bound}
    |c_\gamma| \le C_1 \frac{q_n}{\gamma! (n+1)!} M_{n+1} \kappa_\mathcal{T}^{|\gamma|}q_{\mathcal T}^{|\gamma|} h^{n+1-|\gamma|} \Lambda \le C_2 h^{n+1-|\gamma|} \Lambda,
  \end{equation}
  by using the assumed uniform bound~\cref{eq:lagrange_assumption_bound} on the Lagrange polynomials $\{\lambda_i\}_{i=1}^{q_n}$.
  Thus, using~\eqref{eq:remainder_2_multipoint_power_series} and~\eqref{eq:bound_mon} it follows immediately that for a multi-index $\eta\in\N_0^2$ and for all $\ney \ne \nex$,
  \begin{equation}\label{eq:lagrange_taylor_error}
    \begin{split}
        \left| D^{\eta}_{\ney} R^{(2)}_{\nex, n}[f](\ney)\right| &= \left|D^{\eta}_{\ney} \left(\mathscr{T}[f]_{\nex, n}(\ney) - f_n(\ney; K)\right)\right|\\
        &\le \sum_{\substack{\gamma:|\gamma| \le n\\ \eta\leq\gamma}}  \frac{\gamma!}{(\gamma - \eta)!} \left|c_\gamma(\ney - \nex)^{\gamma-\eta}\right|\\
        &\le C_{R^{(2)}} \Lambda\sum_{\substack{\gamma:|\gamma| \le n\\ \eta\leq\gamma}} h^{n + 1 - |\gamma|} \left|\ney - \nex\right|^{|\gamma| - |\eta|}, \quad C_{R^{(2)}} = C_{R^{(2)}}(n, \eta, \kappa_\mathcal{T}, q_{\mathcal{T}}, f, \Omega) > 0.
    \end{split}
  \end{equation}
  for $|\eta|\leq n$ and $\left| D^{\eta}_{\ney} R^{(2)}_{\nex, n}[f](\ney)\right|=0$ otherwise.
  Again utilizing the bounds~\eqref{eq:Gk_bounds}
  together with~\eqref{eq:lagrange_taylor_error}, we have for $|\alpha| = m + 1>n+1 $ and again for $\ney \ne \nex$,
\begin{subequations}  \begin{equation}\label{eq:lagrange_multipoint_error_final}
      \begin{split}
      \left| D^{\alpha}_{\ney} \left(G_k(\nex, \ney)R_{\nex, n}^{(2)}[f](\ney) \right)\right| &= \left| \sum_{\beta: \beta \le \alpha} \binom{\alpha}{\beta} \left(D_{\ney}^{\alpha-\beta} G_k(\nex, \ney)\right) \left(D^{\beta}_{\ney} R_{\nex, n}^{(2)}[f](\ney)\right)\right| \\
          &\le \sum_{\substack{\beta: \beta \le \alpha,\\ |\beta|\leq n}} \binom{\alpha}{\beta} C_G |\nex - \ney|^{-|\alpha|+|\beta|} C_{R^{(2)}} \Lambda\sum_{\substack{\gamma:|\gamma| \le n\\\beta\leq\gamma}} h^{n + 1 - |\gamma|} \left|\ney - \nex\right|^{|\gamma| -  |\beta|}\\
          &\le  C_{\mathscr{L}} \Lambda \sum_{\gamma:|\gamma| \le n} h^{n + 1 - |\gamma|} |\ney - \nex|^{|\gamma| - m - 1},
      \end{split}
  \end{equation}
  where $C_{\mathscr{L}} = C_{\mathscr{L}}( m,n, k, \kappa_\mathcal{T}, q_\mathcal{T}, f, \Omega) > 0$. For $|\nex-\ney|\geq\delta$, $0<\delta<1$ the above estimate, which holds for all $\nex,\ney\in\Omega$, $\nex\neq\ney$, simplifies to
    \begin{equation}\label{eq:lagrange_multipoint_error_final_2}
      \begin{split}
      \left| D^{\alpha}_{\ney} \left(G_k(\nex, \ney)R_{\nex, n}^{(2)}[f](\ney) \right)\right|
          &\le  C_{\mathscr{L}} \Lambda \sum_{\gamma:|\gamma| \le n} h^{n + 1 - |\gamma|} \delta^{|\gamma| - m - 1}.
      \end{split}
  \end{equation}\label{eq:both_esti}\end{subequations}

  Estimates~\eqref{eq:taylor_error_final}
  and~\eqref{eq:both_esti}
  imply~\eqref{eq:interpolation_estimate}. The proof
  of~\eqref{eq:interpolation_estimate_vec} proceeds component-wise in identical fashion in view of~\eqref{eq:Gk_bounds}, and a
  detailed presentation is omitted for brevity.

\subsection{Proof of \Cref{lem:ordquad_convergence_farfield_optimal}}
\label{proof:ordquad_convergence_farfield_optimal}

    Let $\nex \in \mathcal E_h\subset\overline{\Omega}$ be given and $K\ni \nex$. Let $\delta_0 =\min\{1, (\operatorname{diam}\Omega)^{1/(n-m)}\}$ and assume that $h<\delta_0$.     Consider a sequence of radii $r_j = jh$, ($j=0,\ldots,J+1$), where $J$ is the largest integer such that $hJ<\delta_0$ (i.e., $(J+1)h\geq \delta_0$), and let
    $$
    \widetilde{A}_{j}=\widetilde{A}_j(\nex) :=     \{ \ney\in\Omega: r_{j}^2 \le |\nex -\ney|^2 \le r_{j+1}^2\}, \quad j=0,\ldots,J,
    $$
    in order to define the `meshed' rings (see \Cref{fig:rings}):
    \begin{equation}
A_j=A_j(\nex) := \begin{cases}
 \mathcal N_h=\bigcup\{K_\ell\in\mathcal T_h:K_\ell\cap \widetilde{A}_j(\nex)\neq\emptyset\}, \quad &j = 0,\\
            \bigcup\{K_\ell\in\mathcal T_h:K_\ell \cap \widetilde{A}_j(\nex)\neq \emptyset, K_\ell \not\subset A_{j-1}(\nex)\}, \quad &j=1,\ldots,J,\\
            \overline{\Omega}\setminus\bigcup_{j=0}^J A_j,&j=J+1.
        \end{cases}\label{eq:meshedring_def}
    \end{equation}

Note that with the definitions above we have  $\bigcup_{j=1}^{J+1}A_j = \overline\Omega\setminus\mathcal N_h$. By the triangle inequality it then follows that
\begin{equation}\label{eq:inte_split}\begin{split}
\left|\int_{\Omega\setminus\mathcal N_h} \phi_l(\nex, \ney)\de \ney - \mathcal{Q}_{\Omega \setminus \mathcal{N}_h}[\phi_l(\nex, \cdot)] \right|\leq &
\sum_{j=1}^{J}\left|\int_{A_j} \phi_l(\nex, \ney)\de \ney - \mathcal{Q}_{A_j}[\phi_l(\nex, \cdot)]\right| \\
  &+ \left|\int_{A_{J+1}} \phi_l(\nex, \ney)\de \ney - \mathcal{Q}_{A_{J+1}}[\phi_l(\nex, \cdot)]\right|
\end{split}
\end{equation}
for $l=1,2$. In what follows we derive error estimates for each summand.

    Now, on each of $A_j$ ($j = 1, \ldots, J+1)$, the integrand
    $\phi_l(\nex, \cdot)$ ($l=1,2$) is an infinitely-smooth function, for which a classical
    composite quadrature error estimate~\cite[Sec.\ 7.4]{IsaacsonKeller} provides that,
    given a quadrature formula with non-negative coefficients and a degree of precision $m \ge 0$, the error estimate
    \begin{equation}\label{eq:Ai_estimate_1}
        \begin{split}
        \left|\int_{A_j} \phi_l(\nex, \ney)\de \ney - \mathcal{Q}_{A_j}[\phi_l(\nex, \cdot)]\,\right|  &\le \frac{(2h)^{m+1}}{(m+1)!} 2M_{m+1}^{(l)}(j) \int_{A_j} \de \ney,
        \end{split}
        \quad\quad l = 1, 2,
    \end{equation}
    holds for $j = 1, \ldots, J+1$, with
    \begin{equation}\label{eq:deriv_sup}
        M_{m+1}^{(l)}(j) = \max_{\alpha: |\alpha| = m + 1} \sup_{\ney \in  A_j} |D^\alpha \phi_l(\nex, \ney)|, \quad l = 1, 2.
    \end{equation}
    To obtain a bound for $M_{m+1}^{(1)}(j)$ for $j = 1, \ldots, J$, we use the fact that $jh\leq|\ney - \nex|$ for $\ney\in A_j$, noting that $jh$ satisfies $jh\in (0,\delta_0)$. Therefore, using $\delta=jh$ in~\eqref{eq:interpolation_estimate} in
    \Cref{taylor_lagrange_interpolation_lemma},  we get
    \begin{equation}\label{eq:deriv_estimate_ring}
        \begin{split}
            M_{m+1}^{(1)}(j) &\le C_{\mathscr{T}} \frac{1}{(jh)^{m-n}} + C_{\mathscr{L}} \Lambda \sum_{|\gamma| \le n} h^{n+1 - |\gamma|} \frac{1}{(jh)^{m+1-|\gamma|}}\\
            &\le \frac{1}{h^{m-n}}\left(\frac{C_{\mathscr{T}}}{j^{m-n}}  +C_{\mathscr{L}} \Lambda\sum_{|\gamma| \le n} \frac{1}{j^{m+1-|\gamma|}}\right).
        \end{split}
    \end{equation}

To obtain a bound for $\int_{A_j} \de\ney$ in~\eqref{eq:Ai_estimate_1}, on the other hand, in view of~\eqref{eq:meshedring_def} and since every element $K_\ell \subset \mathcal{T}_h$ has maximum diameter $h$, for all $\ney \in A_j$ the inequality $jh  \leq  |\nex - \ney| < (j+2)h$ holds from which it follows that
    \begin{equation}\label{eq:area_estimate}
        \int_{A_j} \de \ney \le \pi \left(j+2 \right)^2 h^2 - \pi j^2 h^2 = 4\pi (j + 1) h^2, \quad j=1,\ldots,J.
    \end{equation}
    Thus, using~\eqref{eq:deriv_estimate_ring} together
    with~\eqref{eq:area_estimate}, the estimate~\eqref{eq:Ai_estimate_1} becomes
    \begin{equation}\label{eq:Ai_estimate_2}
        \begin{split}
            \left|\int_{A_j} \phi_1(\nex, \ney)\de \ney- \mathcal{Q}_{A_j}[\phi_1(\nex, \cdot)]\,\right|  &\le  4\pi \frac{(2h)^{n+3}}{(m+1)!} (j + 1)\left(\frac{C_{\mathscr{T}}}{j^{m-n}}  + C_{\mathscr{L}}\Lambda \sum_{|\gamma| \le n} \frac{1}{j^{m+1-|\gamma|}}\right)\\
            &\le C_1 \Lambda \frac{(2h)^{n+3}}{(m+1)!} \frac{1}{j^{m-n-1}},\quad j=1,\ldots,J.
        \end{split}
    \end{equation}
    To obtain a quadrature error estimate over $ \bigcup_{j=1}^J A_j$ one simply sums over the rings $A_j$, $j = 1, \ldots, J$:
\begin{equation}\label{showing_far_quadrature_optimal_estimate_V}
        \begin{split}
            \sum_{j=1}^{J} \left|\int_{A_j} \phi_1(\nex, \ney)\de \ney - \mathcal{Q}_{A_j}[\phi_1(\nex, \cdot)]\,\right|  &\le C_2 \Lambda  h^{n+3}\sum_{j=1}^{J} \frac{1}{j^{m-n-1}}\\
            &\le C_3 \Lambda
        \begin{cases}
            h^{n+3}, &\quad m > n + 2\\
            h^{n+3}\left|\log h\right|, &\quad m = n + 2\\
            h^{n+2}, &\quad m = n + 1,
        \end{cases}
        \end{split}
    \end{equation}
    the last sum being estimated above, for $m > n + 2$ by a constant (as $\sum_{j=1}^\infty j^{-\zeta}$ converges for $\zeta > 1$), for $m = n + 2$ by the harmonic numbers $H_J$ that grow as $\log J\lesssim |\log h|$, and for $m = n + 1$  by $J<\delta_0 h^{-1}$.

   To estimate the quadrature error associated with the integral over $A_{J+1}$, we observe that $\delta_0<(J+1)h\leq|\ney - \nex|$ for $\ney\in A_{J+1}$. Therefore, using $\delta=\delta_0$ in~\eqref{eq:interpolation_estimate} and from the assumption $h<\delta_0$ (see beginning of proof) we obtain
      \begin{equation}\label{eq:deriv_estimate_ring_J+1}
        \begin{split}
            M_{m+1}^{(1)}(J+1) &\le C_{\mathscr{T}} \frac{1}{\delta_0^{m-n}} + C_{\mathscr{L}} \Lambda \sum_{|\gamma| \le n} h^{n+1 - |\gamma|} \frac{1}{\delta_0^{m+1-|\gamma|}}\\
            &\le C_{\mathscr{T}} \frac{1}{\delta_0^{m-n}} + C_{\mathscr{L}} \Lambda \sum_{|\gamma| \le n}  \frac{1}{\delta_0^{m-n}}\\
         &= \frac{1}{\delta_0^{m-n}}\left(C_{\mathscr{T}} + C_{\mathscr{L}} \Lambda q_n\right).
        \end{split}
    \end{equation}
It hence follows from~\eqref{eq:Ai_estimate_1} that
    \begin{equation}\label{eq:Ai_estimate_3}
        \begin{split}
          \left|\int_{A_{J+1}} \phi_1(\nex, \ney)\de \ney- \mathcal{Q}_{A_{J+1}}[\phi_1(\nex, \cdot)]\,\right|
         &\le \frac{(2h)^{m+1}}{(m+1)!} \frac{2}{\delta_0^{m-n}}\left(C_{\mathscr{T}} + C_{\mathscr{L}} \Lambda q_n\right) \int_{A_{J+1}} \de \ney\\
           &       \le \frac{(2h)^{m+1}}{(m+1)!} \frac{2}{\delta_0^{m-n}}\left(C_{\mathscr{T}} + C_{\mathscr{L}} \Lambda q_n\right) \int_{\Omega} \de \ney\\
           &    \le  C_4 \Lambda h^{m+1}.
        \end{split}
    \end{equation}
 In view of~\eqref{eq:inte_split}, from~\eqref{showing_far_quadrature_optimal_estimate_V} and~\eqref{eq:Ai_estimate_3} we have that the desired bound in~\eqref{far_quadrature_optimal_estimate_V} is established with $\mathcal{C}_1 = C_3 + C_4$.

    The proof for $\phi_2$ proceeds in identical fashion. To bound $M_{m+1}^{(2)}(j)$ for each $j = 1, \ldots, J$, we use the fact that $jh\leq |\ney - \nex|$ for all $\ney\in A_{j}$, so from~\eqref{eq:interpolation_estimate_vec} with $\delta=jh<\delta_0$ we obtain
    \begin{equation}\label{eq:deriv_estimate_ring_vec}
        \begin{split}
            M_{m+1}^{(2)}(j) &\le C_{\mathscr{T}}' \frac{1}{(jh)^{m+1-n}} + C_{\mathscr{L}}' \Lambda \sum_{|\gamma| \le n} h^{n+1 - |\gamma|} \frac{1}{(jh)^{m+2-|\gamma|}}\\
            &\le h^{n-m-1}\left( \frac{{C}_{\mathscr{T}}'}{j^{m+1-n}} + {C}_{\mathscr{L}}'\Lambda \sum_{|\gamma| \le n} \frac{1}{j^{m+2-|\gamma|}}\right).
        \end{split}
    \end{equation}
    Using~\eqref{eq:deriv_estimate_ring_vec} together
    with~\eqref{eq:area_estimate}, the estimate~\eqref{eq:Ai_estimate_1} in the case $l = 2$
    becomes
    \begin{equation}\label{eq:Ai_estimate_2_vec}
        \begin{split}
            \left|\int_{A_j} \phi_2(\nex, \ney)\de \ney- \mathcal{Q}_{A_{j}}[\phi_2(\nex, \cdot)]\,\right|  &\le 4\pi \frac{(2h)^{n+2}}{(m+1)!} (j + 1)\left( \frac{C_{\mathscr{T}}'}{j^{m+1-n}} +{C}_{\mathscr{L}}'\Lambda \sum_{|\gamma| \le n} \frac{1}{j^{m+2-|\gamma|}}\right)\\
            &\le C_5 \Lambda \frac{(2h)^{n+2}}{(m+1)!} \frac{1}{j^{m-n}}.
        \end{split}
    \end{equation}
The quadrature error estimate over  $\bigcup_{j=1}^J A_j$ results  from summing over the rings $A_j$, $j = 1, \ldots, J$:
    \begin{equation}
        \begin{split}
            \sum_{j=1}^{J} \left|\int_{A_j} \phi_2(\nex, \ney)\de \ney - \mathcal{Q}_{A_{j}}[\phi_2(\nex, \cdot)]\,\right|  &\le C_6 \Lambda h^{n+2}\sum_{j=1}^{J} \frac{1}{j^{m-n}}\\
            &\le C_7\Lambda
        \begin{cases}
            h^{n+2}, &\quad m \ge n + 2\\
            h^{n+2}\left|\log h\right|, &\quad m = n + 1,
        \end{cases}
        \end{split}
    \end{equation}
where  we have used once again the boundedness of the series $\sum_{j=1}^\infty j^{-\zeta}$, $\zeta>1$, as well as  the relation $J < \delta_0h^{-1}$ together with the logarithmic growth of the harmonic series for $m = n + 1$. Finally, to estimate the quadrature error over $A_{J+1}$ we use the fact that $|\nex-\ney|>(J+1)h>\delta_0$ for all $\ney\in A_{J+1}$ which allows us to use the bound~\eqref{eq:interpolation_estimate_vec} with $\delta=\delta_0$. Therefore, from that bound and the assumption $h<\delta_0$ we obtain
      \begin{equation}\label{eq:deriv_estimate_ring_J+1_2}
        \begin{split}
            M_{m+1}^{(2)}(J+1) &\le C'_{\mathscr{T}} \frac{1}{\delta_0^{m+1-n}} + C'_{\mathscr{L}} \Lambda \sum_{|\gamma| \le n} h^{n+1- |\gamma|} \frac{1}{\delta_0^{m+2-|\gamma|}}\\
            &\le C'_{\mathscr{T}} \frac{1}{\delta_0^{m+1-n}} + C'_{\mathscr{L}} \Lambda \sum_{|\gamma| \le n}  \frac{1}{\delta_0^{m+1-n}}\\
         &= \frac{1}{\delta_0^{m+1-n}}\left(C'_{\mathscr{T}} + C'_{\mathscr{L}} \Lambda q_n\right),
        \end{split}
    \end{equation}
so from~\eqref{eq:Ai_estimate_1} in the case $l = 2$ it follows that
    \begin{equation}\label{eq:Ai_estimate_3_4}
        \begin{split}
          \left|\int_{A_{J+1}} \phi_2(\nex, \ney)\de \ney- \mathcal{Q}_{A_{J+1}}[\phi_2(\nex, \cdot)] \,\right|
         &\le \frac{(2h)^{m+1}}{(m+1)!} \frac{2}{\delta_0^{m+1-n}}\left(C'_{\mathscr{T}} + C'_{\mathscr{L}} \Lambda q_n\right) \int_{A_{J+1}} \de \ney\\
           &       \le \frac{(2h)^{m+1}}{(m+1)!} \frac{2}{\delta_0^{m+1-n}}\left(C'_{\mathscr{T}} + C'_{\mathscr{L}} \Lambda q_n\right) \int_{\Omega} \de \ney\\
           &    \le  C_{8}\Lambda h^{m+1}.
        \end{split}
    \end{equation}

    This establishes~\eqref{far_quadrature_optimal_estimate_W} with $\mathcal{C}_2 = C_7 + C_8$, and the proof is thus complete.

\section{Numerical examples}\label{sec:numer}
\subsection{Validation examples}
In order to validate the proposed methodology, we manufacture potential evaluations with scalar and vectorial volume densities, $f:\Omega\to\C$ and $\mathbf{f} \coloneqq f\vv d:\Omega\to\C^2$ ($\vv d\in\R^2$ being a constant vector), respectively,  that do not involve the evaluation of volume integrals and so are more easily computable as a reference solution for assessing error levels. To do this, we start off by selecting a known smooth function $u:\Omega\to\C$ and define the scalar source density as
$f \coloneqq (\Delta+k^2)u:\Omega\to\C.$
By Green's third identity it readily follows that $\mathcal V_k [f]=v_{\rm ref}$ in $\Omega$, where the reference potential is given by \begin{equation}\label{eq:GF_validation}
  v_{\rm ref}(\nex) \coloneqq -u(\nex)-\int_{\Gamma} \left\{\frac{\p G_k(\nex,\ney)}{\p \nu(\ney)}   u(\ney)-G_k(\nex,\ney) \frac{\p u}{\p \nu}(\ney)  \right\}\de s(\ney)
\end{equation}
in terms of the known function, $u$, and the single- and double-layer Laplace/Helmholtz potentials applied to its Neumann and Dirichlet traces, respectively. Similarly, we have $\mathcal W_{k}[\mathbf{f}]=w_{\rm ref}$ in $\Omega$, where
\begin{equation}\label{eq:GF_validation_W}
w_{\rm ref}(\nex) \coloneqq \vv d \cdot \nabla\left\{ u(\nex)+\int_{\Gamma} \left\{\frac{\p G_k(\nex,\ney)}{\p \nu(\ney)}   u(\ney)-G_k(\nex,\ney) \frac{\p u}{\p \nu}(\ney) \right\}\de s(\ney)\right\}.
\end{equation}
Since both $v_{\rm ref}$ and $w_{\rm ref}$ only entail the evaluation of layer potentials and their gradients, a highly accurate numerical approximation of the boundary integrals yields a reliable expression for the volume potentials that can be used as a reference to measure the numerical errors.

In all the numerical examples considered in this section, we employ
$$u(\nex) =\frac{1}{k^2-k_0^2}\e^{ik_0\nex\cdot{\vv \theta}}+\frac{1}{5(k^2-s^2)}\sum_{j=1}^4\e^{-s|\nex-\nez_j|^2},$$
where $|\vv \theta|=1$ and $\nez_j\in\Gamma$, $j=1,\ldots,4$, which gives rise to the density function
\begin{equation}\label{eq:sources_ex}
f(\nex) =\e^{ik_0\nex\cdot{\vv \theta}}++\frac{1}{5(k^2-s^2)}\sum_{j=1}^4(4s^2|\nex-\nez_j|^2-4s+k^2)\e^{-s|\nex-\nez_j|^2},
\end{equation}
and make the selections $k=2\pi$, $k_0=3\pi$, $\vv d =(1,1)$, $s=10$, and $\vv\theta=(\cos\frac{\pi}3,\sin\frac{\pi}3)$. The numerical results are demonstrated on six specific domains. We first consider three domains with smooth boundary parametrizations, namely, the unit disk, the kite~\cite{COLTON:2012}, and the jellyfish~\cite{gomez2021regularization}, whose boundaries are parametrized, respectively, by the smooth curves
\begin{subequations}\begin{align}
\bol\gamma_0(t)&=(\cos t, \sin t),\\
\bol\gamma_1(t)&=(\cos t+0.65 \cos 2 t-0.65,1.5 \sin t),\quad\text{and}\label{eq:kite}\\
\bol\gamma_2(t)&=\{1+0.3 \cos (4 t+2 \sin t)\}(\sin t,-\cos t),\qquad 0\leq t \leq 2 \pi.
\end{align}
\label{eq:params}\end{subequations}
These  domains encompass curved triangles as described in~\Cref{sec:curved_triangles}. Finally, our last three domains, namely, the windmill, the Nazca bird, and the snowflake, are polygonal. The latter encompasses a large number of small-scale features, its meshes are quite heterogeneous ($\max_{K\in\mathcal T_h} h_K/\min_{K\in\mathcal T_h}h_K>10$), and it includes several poor-quality triangles in some cases.

\begin{figure}
  \centering
  \begin{subfigure}[b]{0.31\textwidth}
      \centering
      \includegraphics[scale=0.2]{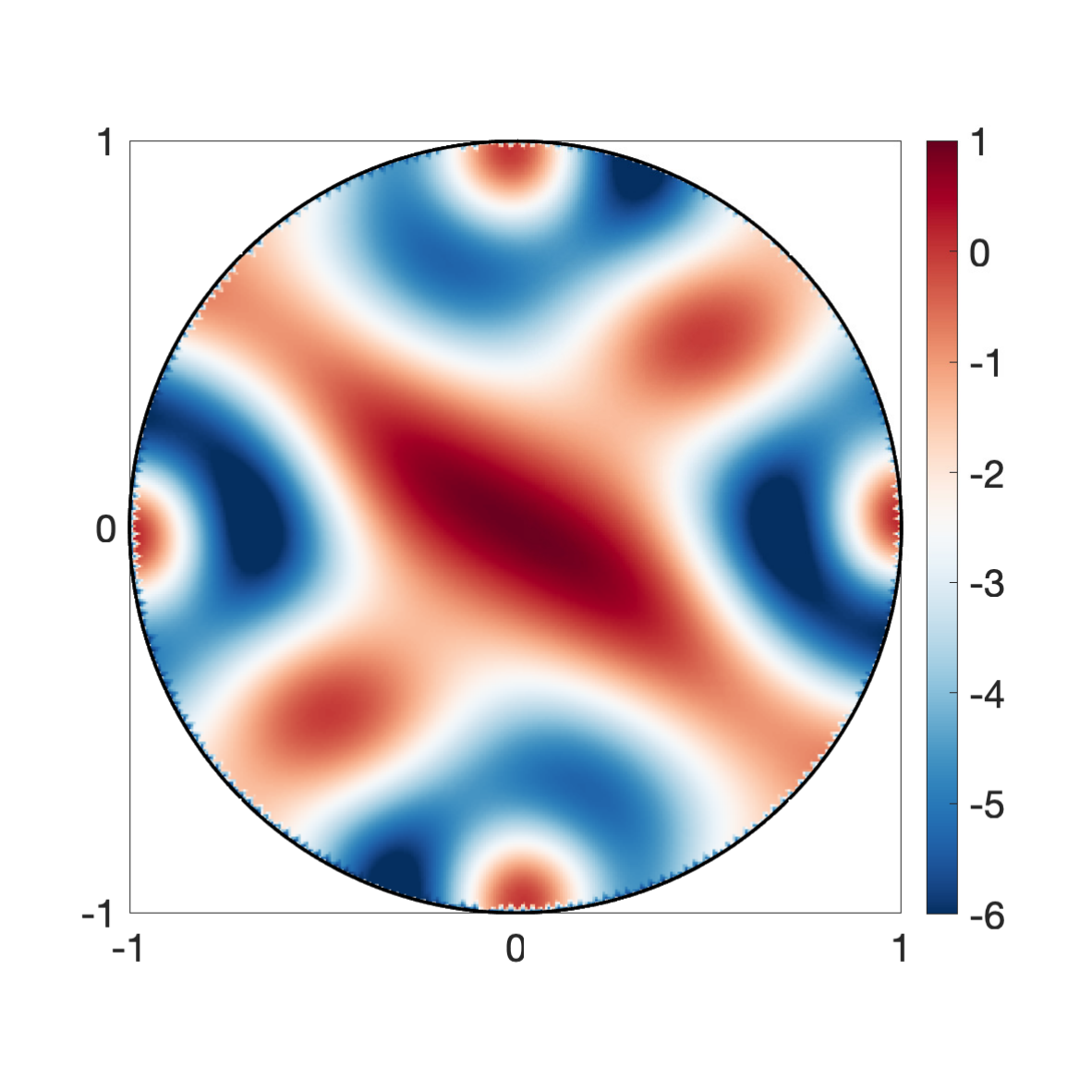}\\
      \includegraphics[width=\linewidth]{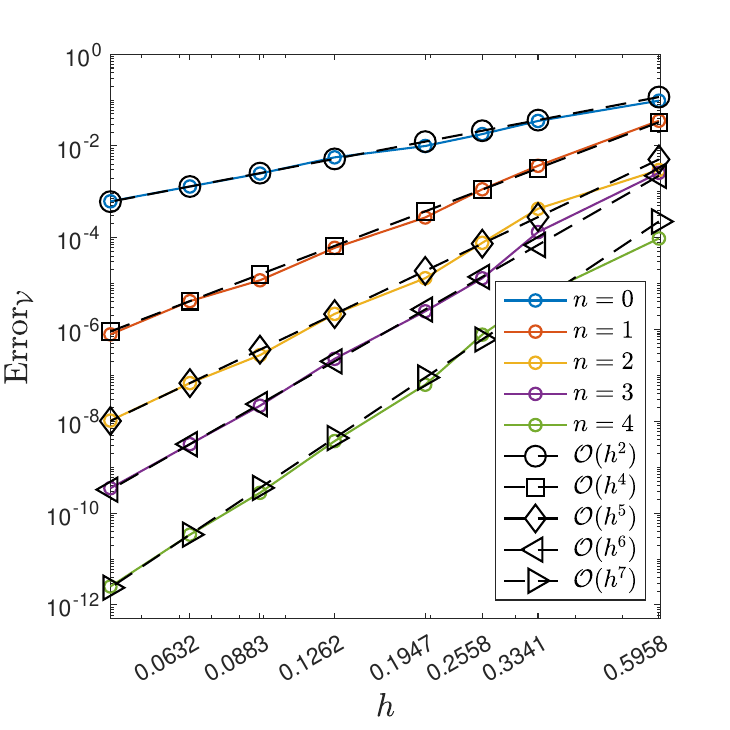}\\
      \includegraphics[width=\linewidth]{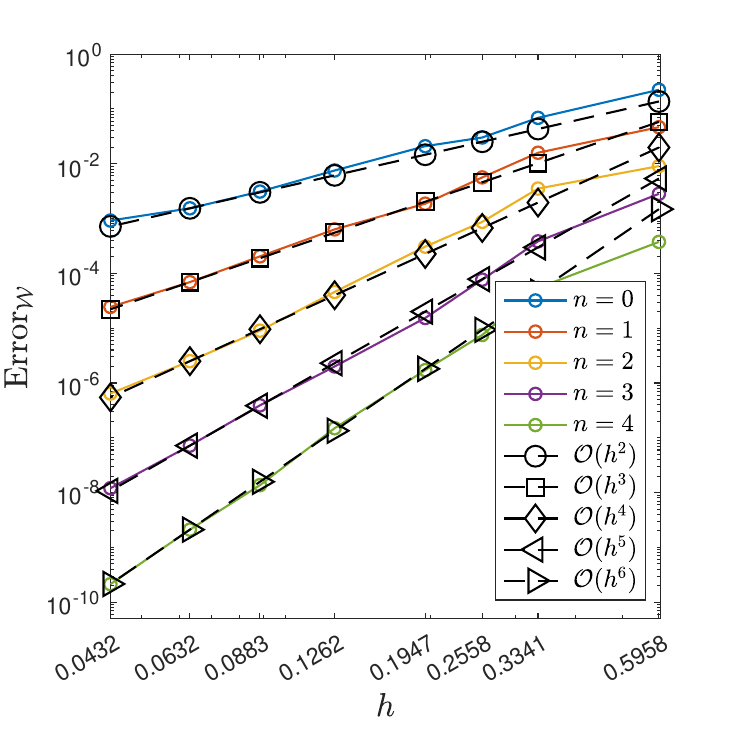}
      \caption{Circle}
      \label{fig:circle_tris}
  \end{subfigure}\quad
  \begin{subfigure}[b]{0.31\textwidth}
      \centering
      \includegraphics[scale=0.2]{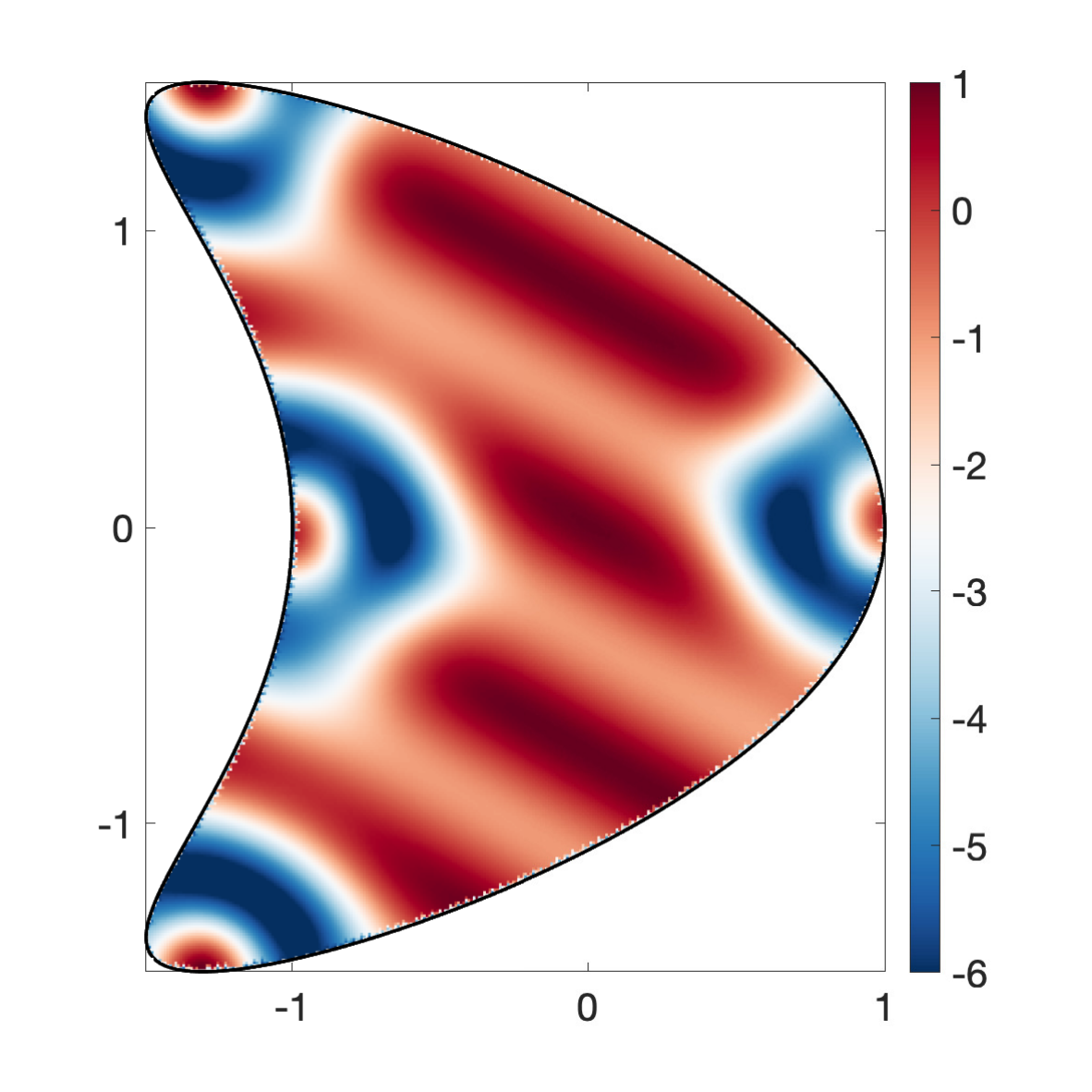}\\
      \includegraphics[width=\linewidth]{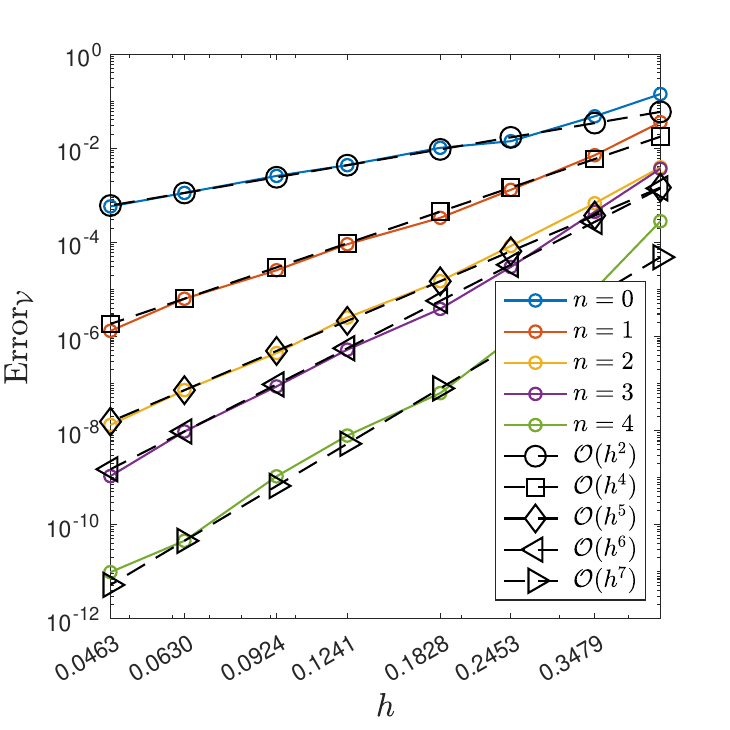}\\
      \includegraphics[width=\linewidth]{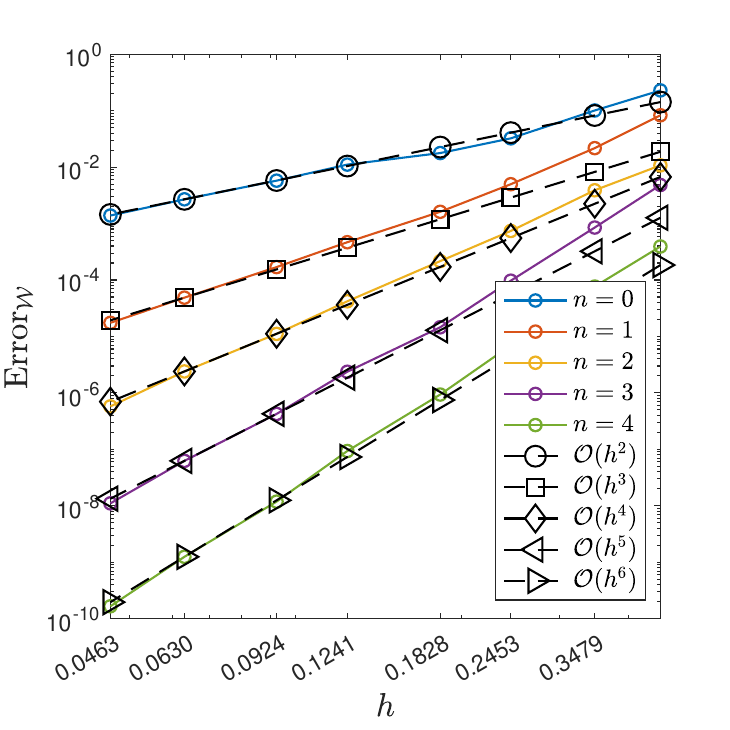}
      \caption{Kite}
      \label{fig:kite_tris}
  \end{subfigure}
  \quad
  \begin{subfigure}[b]{0.31\textwidth}
      \centering
      \includegraphics[scale=0.2]{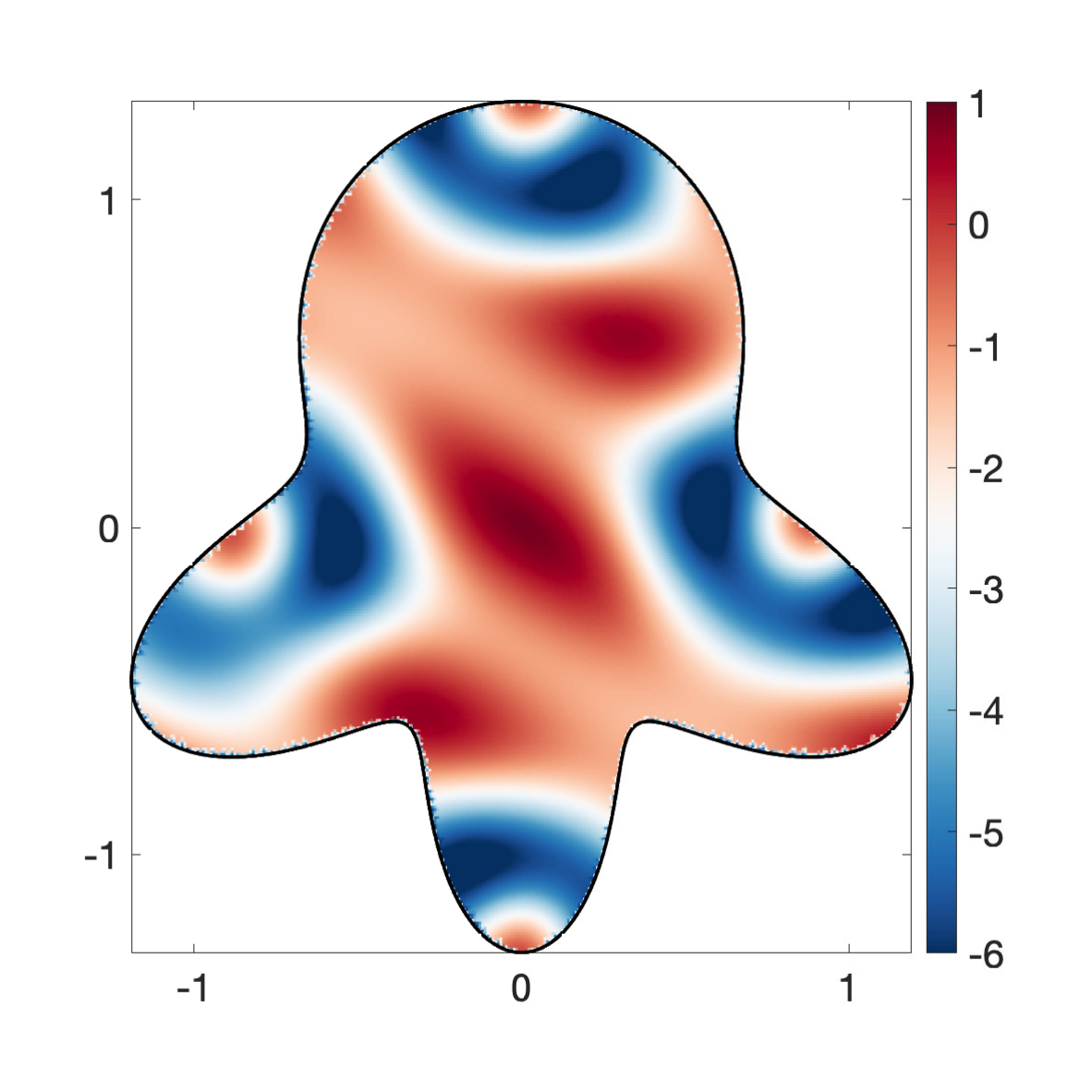}\\
      \includegraphics[width=\linewidth]{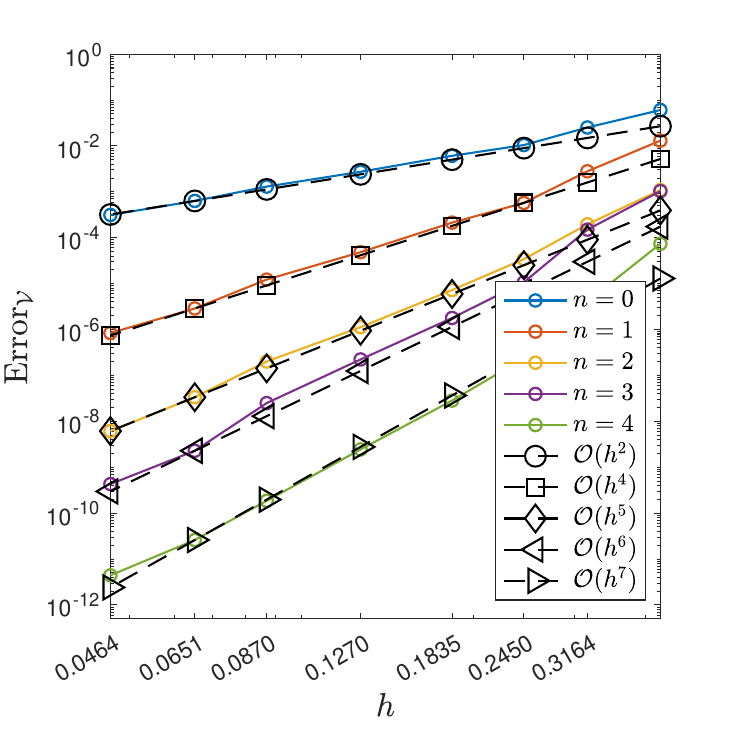}\\
      \includegraphics[width=\linewidth]{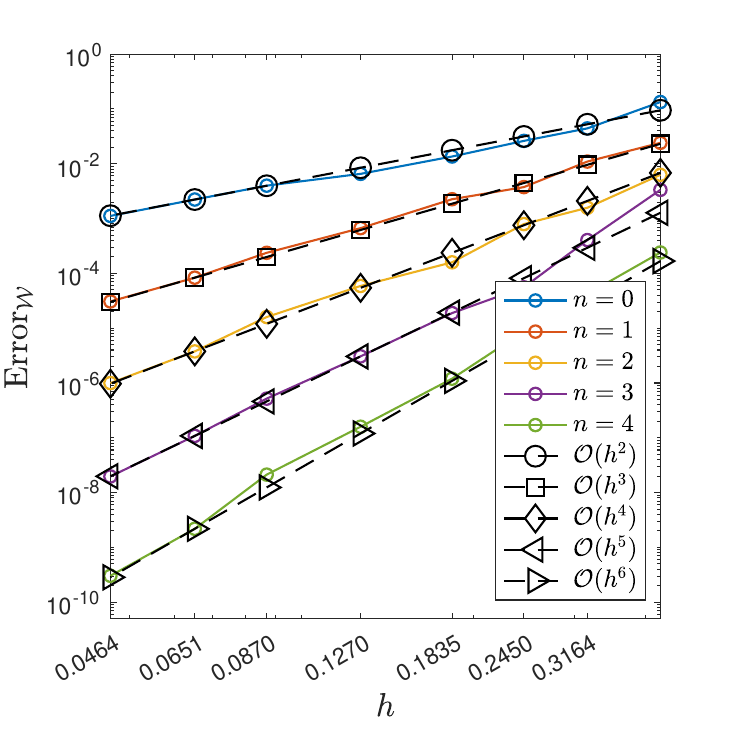}
      \caption{Jellyfish}
      \label{fig:jellyfish_tris}
  \end{subfigure}
     \caption{Numerical accuracy in the evaluation of volume potentials $\mathcal V_k$ and $\mathcal W_k$. Top Row: Plot of the real part of the density functions $f$, defined in~\eqref{eq:sources_ex}, that are used in the numerical experiments for each of the three smooth domains $\Omega$. Relative errors in the numerical evaluation of $\mathcal V_k$ (middle row) and $\mathcal W_k$ (bottom row) for various discretization sizes.}
     \label{fig:tris_splines}
\end{figure}

\begin{figure}
  \centering
  \begin{subfigure}[b]{0.31\textwidth}
      \centering
      \includegraphics[scale=0.2]{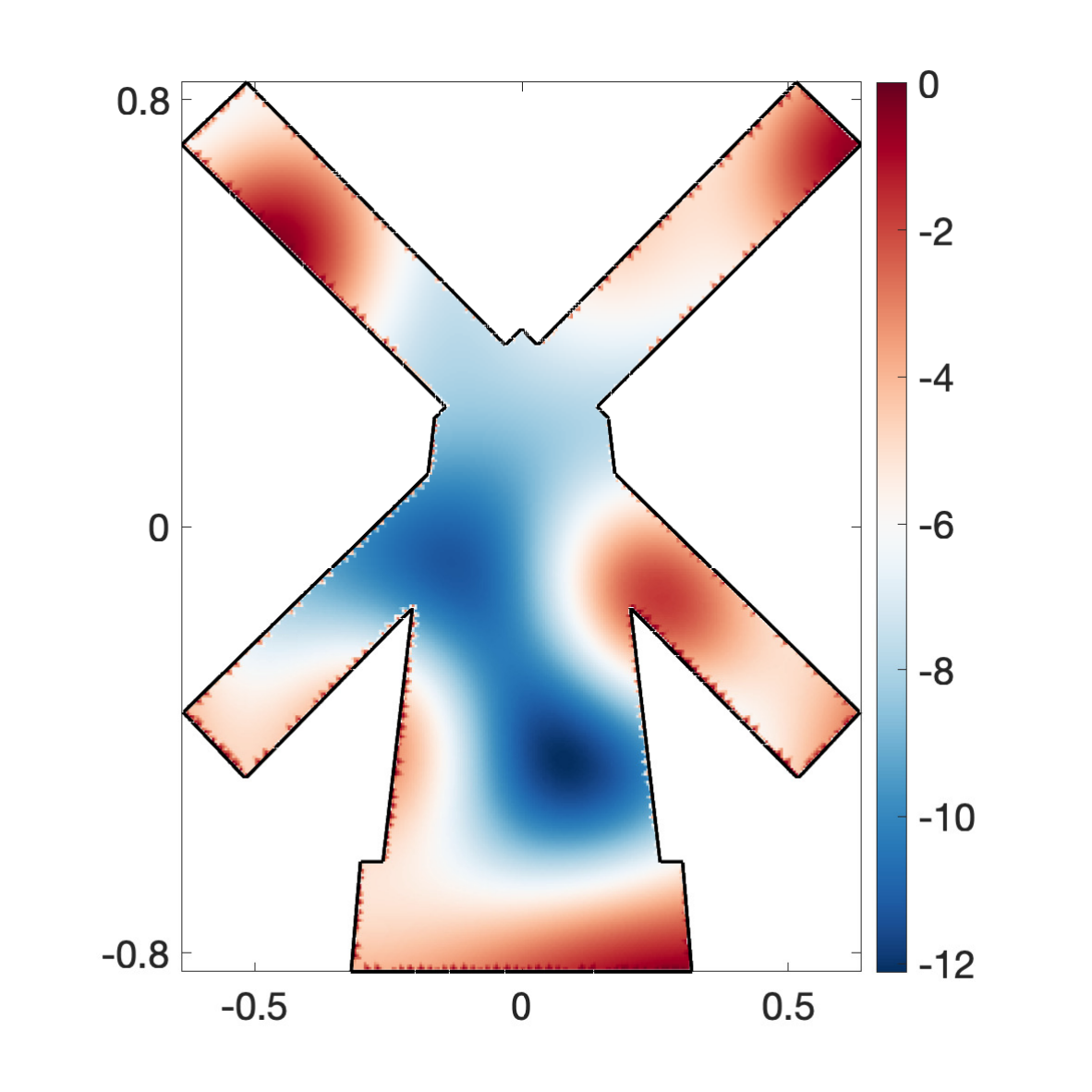}\\
      \includegraphics[width=\linewidth]{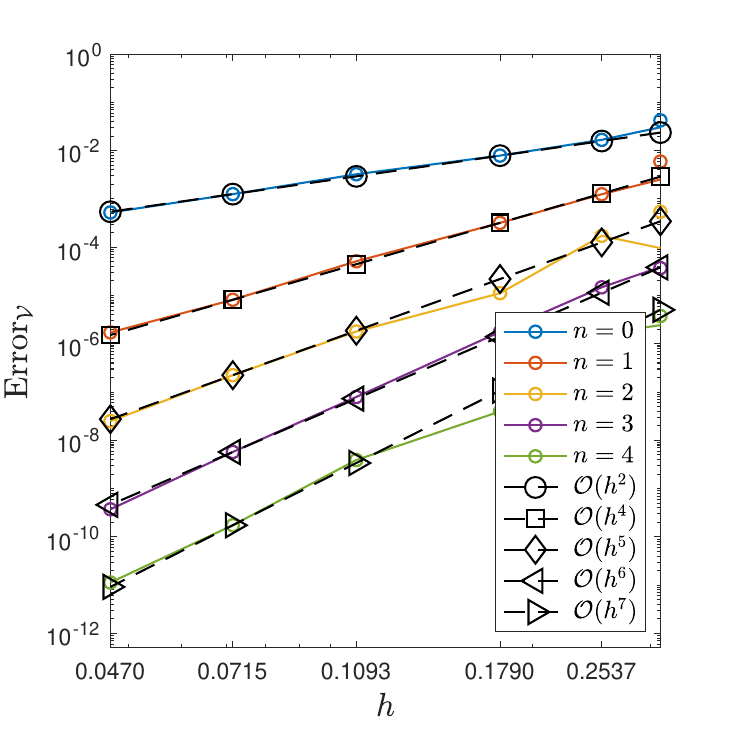}\\
      \includegraphics[width=\linewidth]{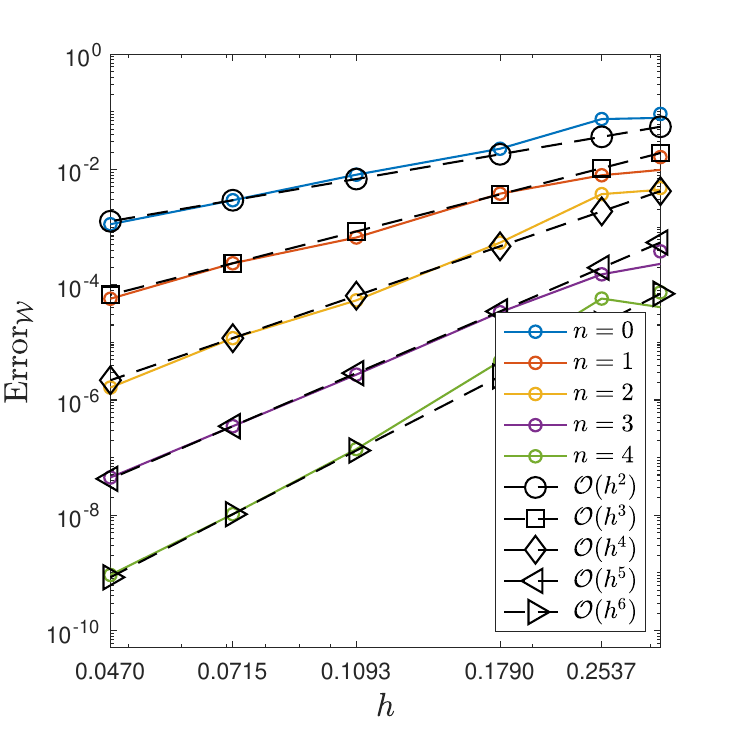}
      \caption{Windmill}
      \label{fig:windmill_tris}
  \end{subfigure}
  \quad
  \begin{subfigure}[b]{0.31\textwidth}
      \centering
      \includegraphics[scale=0.2]{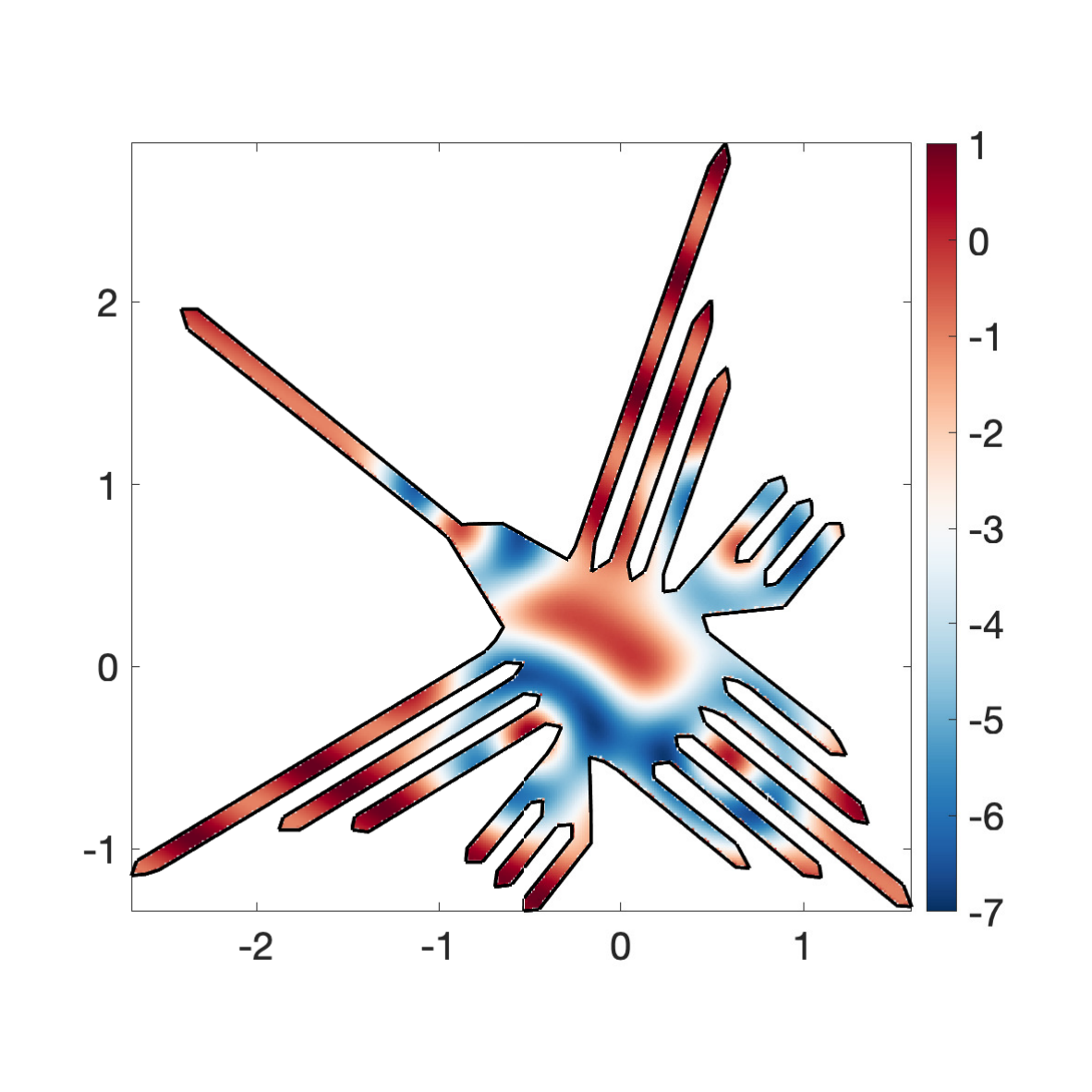}\\
      \includegraphics[width=\linewidth]{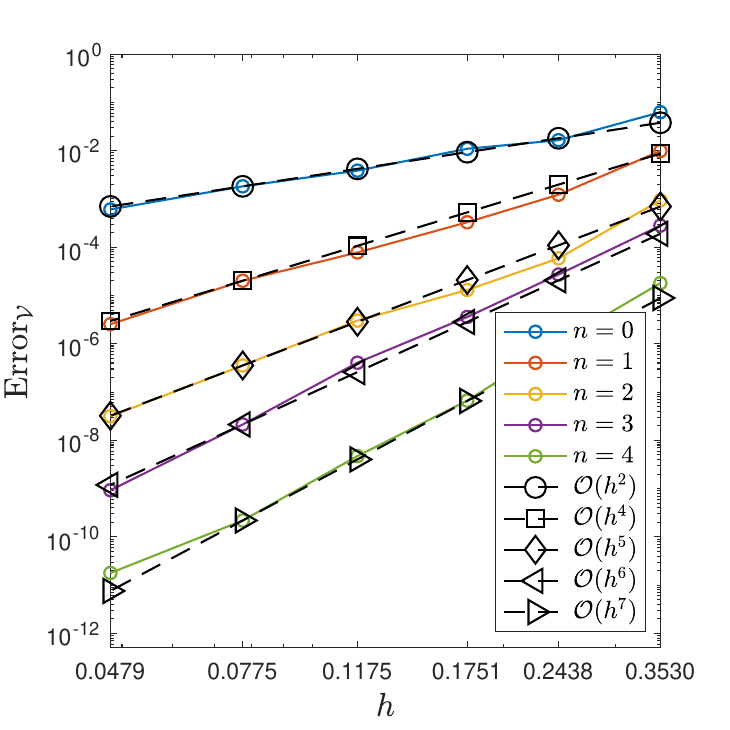}\\
      \includegraphics[width=\linewidth]{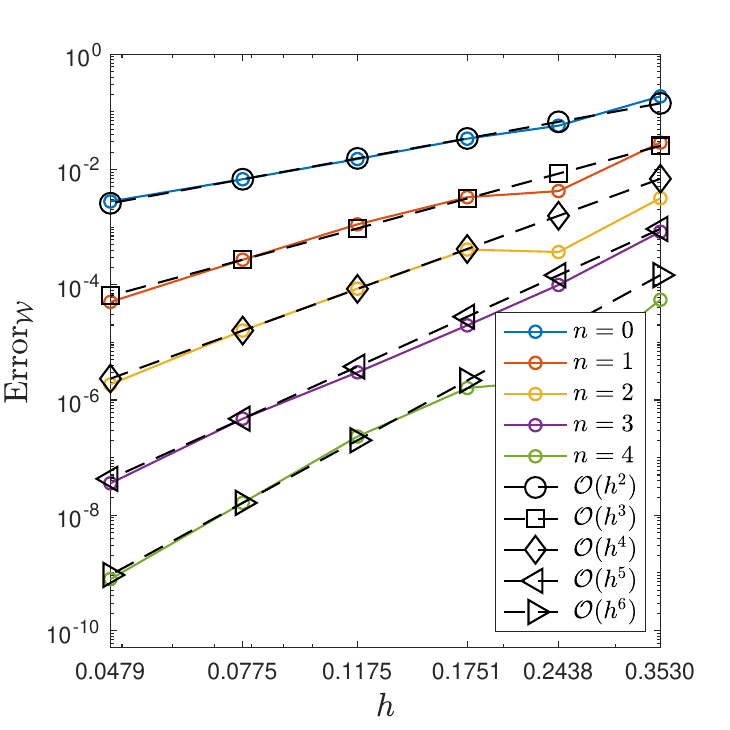}
      \caption{Nazca Bird}
      \label{fig:nazcabird_tris}
  \end{subfigure}
  \quad
  \begin{subfigure}[b]{0.31\textwidth}
      \centering
      \includegraphics[scale=0.2]{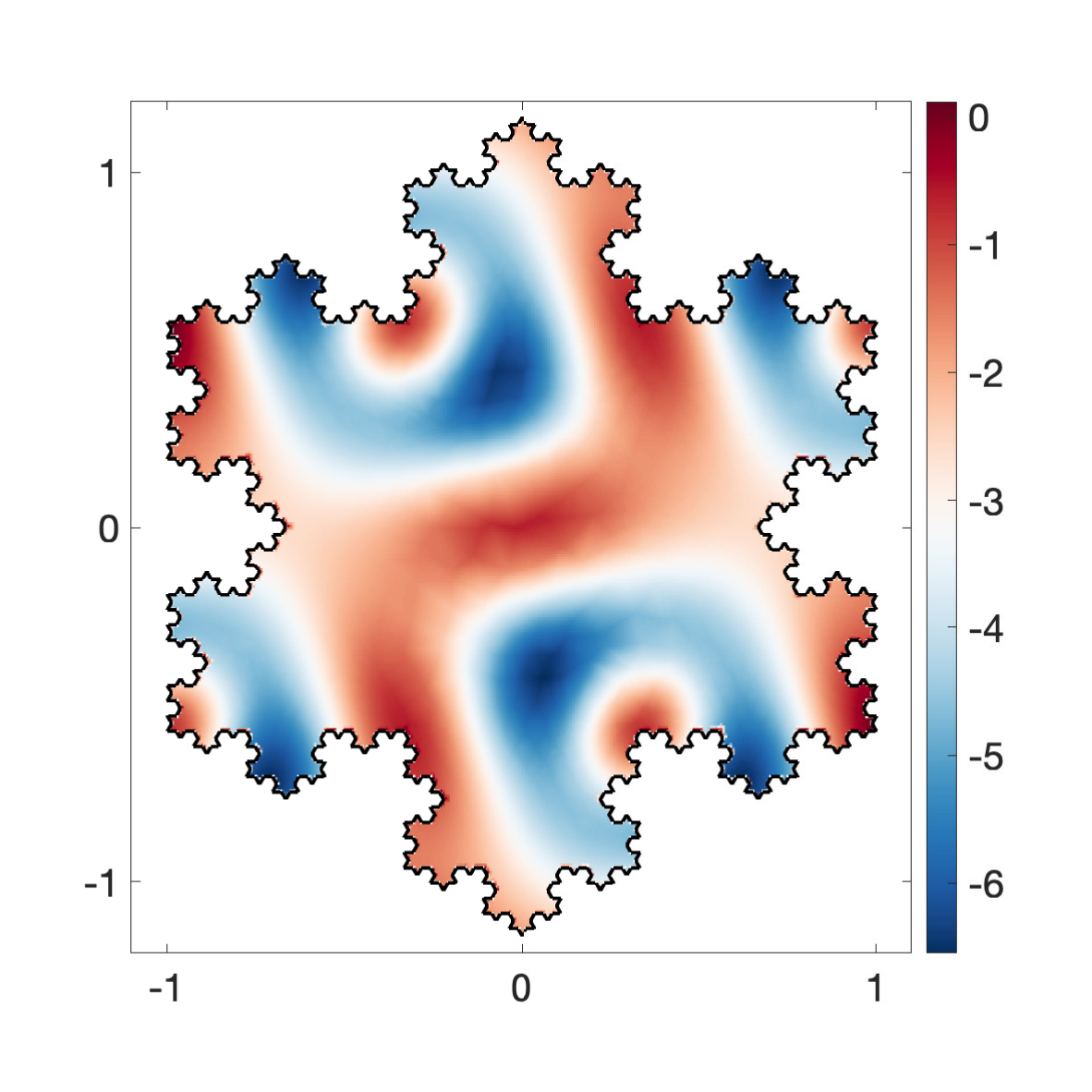}\\
      \includegraphics[width=\linewidth]{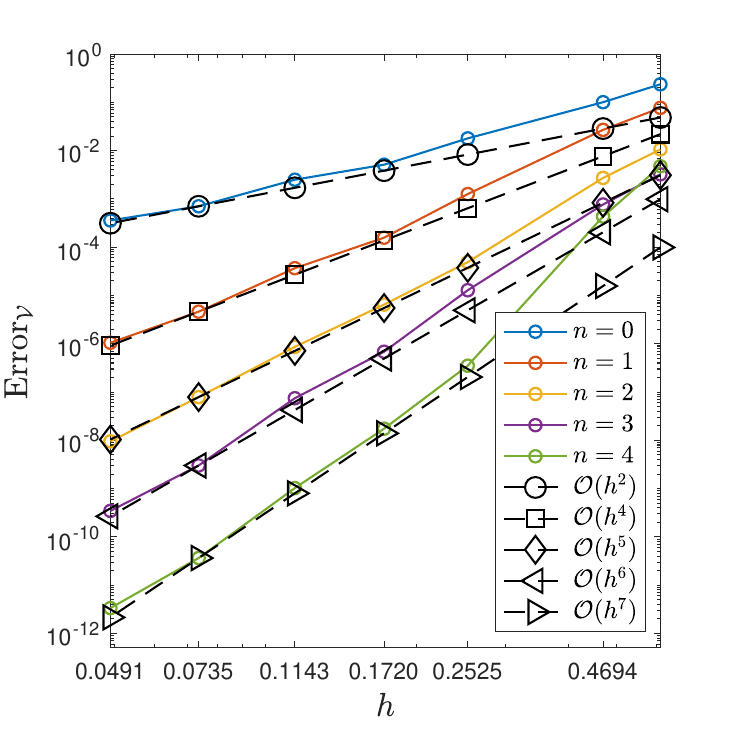}\\
      \includegraphics[width=\linewidth]{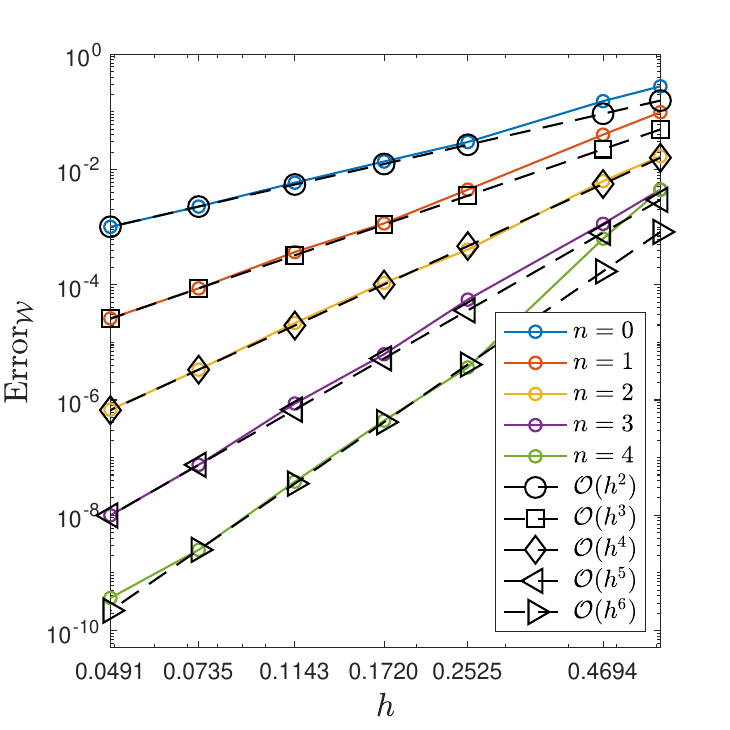}
      \caption{Snowflake}
      \label{fig:snowflake_tris}
  \end{subfigure}
     \caption{Numerical accuracy in the evaluation of volume potentials $\mathcal V_k$ and $\mathcal W_k$ using the proposed polynomial density interpolation method. Top Row: Plot of the real part of the density functions $f$, defined in~\eqref{eq:sources_ex}, that are used in the numerical experiments for each of the three polygonal domains $\Omega$. Middle and Bottom Rows: Relative errors in the numerical evaluation of $\mathcal V_k$ (middle row) and $\mathcal W_k$ (bottom row) for various discretization sizes.}
     \label{fig:tris_polygons}
\end{figure}

In the numerical examples that follow, the relative numerical errors are measured by means of the formulae
\begin{equation}\label{eq:rel_error}
\text{Error}_{\mathcal V} \coloneqq \frac{\displaystyle\max_{j\in\{1,\ldots,N\}}\left|v_{\rm approx}(\bol\xi_j)-v_{\rm ref}(\bol\xi_j)\right|}{\displaystyle\max_{j\in\{1,\ldots,N\}}\left|v_{\rm ref}(\bol\xi_j)\right|}\andtext \text{Error}_{\mathcal W} \coloneqq \frac{\displaystyle\max_{j\in\{1,\ldots,N\}}\left|w_{\rm approx}(\bol\xi_j)-w_{\rm ref}(\bol\xi_j)\right|}{\displaystyle\max_{j\in\{1,\ldots,N\}}\left|w_{\rm ref}(\bol\xi_j)\right|},
\end{equation}
which provide approximate relative errors in  the natural $C^0(\overline{\Omega})$-norm, where $\{\xi_j\}_{j=1}^N$ are the volumetric quadrature nodes and where $v_{\rm approx}$ and $w_{\rm approx}$ denote respectively the approximate potentials $\mathcal V_{k}[f]$ and $\mathcal W_{k}[\mathbf{f}]$. The evaluation of each operator is performed in the experiments that follow using Vioreanu-Rokhlin quadrature nodes for integration and interpolation. The results are presented in \Cref{fig:tris_splines} and~\Cref{fig:tris_polygons} for the smooth and polygonal domains, respectively, showcasing the obtained errors, ${\rm Error}_{\mathcal V}$ and ${\rm Error}_{\mathcal W}$, for meshes of varying sizes $h$ and interpolation degrees $n\in\{0,1,2,3,4\}$. We limit our focus to these interpolation degrees for two primary reasons. On one hand, the method's efficiency decreases as the interpolation degree increases, and on the other, to manage the error-amplifying effect of the translation formula~\eqref{eq:trans_pol_tri}, discussed in detail below in \Cref{sec:limitations}.  Notably, we observe clear convergence orders for interpolation degrees $n\in\{0,2,3,4\}$, that match the error estimates  established in \Cref{thm:tri_error_analysis}. However, interestingly, for the $\mathcal V_{k}$ operator, we observe evidence of super-convergence in the case $n=1$ even as the expected rate is observed for the $\mathcal W_k$ operator (see \Cref{rem:tri_simplify_error_estimates}).

\subsection{Timing results}
  The method has been implemented in Julia, in the open-source \texttt{Inti.jl}\cite{Inti} package also developed by the authors.  Some single-core double-precision performance statistics of the method for computing $\mathcal{V}_k$ are given in Tables~\ref{tab:fmm_laplace} and~\ref{tab:fmm_helmholtz} for a volume potential computed over a unit circle; all times reported are in seconds. The code was executed in Julia v1.10 on a Linux machine with an AMD Ryzen 7950X3D CPU. The column marked $T_{\mathcal{B}_{\Gamma, \mathrm{close}}}$ lists the precomputation cost of setting up near-singular operator evaluation $\mathcal{B}_\Gamma[\varphi, \psi]$. The column $T_{\mathcal{Q}_\Omega,\mathcal{B}_\Gamma, P_\alpha}$ shows the cost of volume forward maps, boundary forward maps and particular solution at all points in the domain; $T_C$ denotes the cost of (i) Factorization of the Vandermonde-like matrix $\bold{V}$, of (ii) Computing translation and scaling quantities $\vv{c}_K$, $r_K$, and $\bold{P}\bold{V}^{-1}$, and finally of (iii) Construction of the sparse correction matrix. The column $T_{\mathrm{pre}}$ displays the cost of all precomputations. The reported times do not make use of vectorized (or batched) use of the FMM, wherein the same FMM call works on multiple densities in a single pass.
  Turning to online costs, $T_{\mathrm{corr}}$ denotes the time required for application of the sparse correction matrix while $T_{\mathrm{app}}$ denotes the full cost of applying $\mathcal{V}_k$ during the online phase. The online throughput is $N/T_{\mathrm{app}}$; as evidenced by the $\%\mathrm{fmm}$ column the method captures the full speed of the FMM at all orders.
\begin{table}[h]
    \centering
    \begingroup
    \setlength{\tabcolsep}{4pt}
    \begin{tabular}{ c|c|c|c|c|c||c||c|c|c|c }
    $h$ & $N$ & $T_{\mathcal{B}_\Gamma, \mathrm{close}}$ & $T_{\mathcal{Q}_\Omega, \mathcal{B}_\Gamma, P_{\alpha}}$ & $T_C$ & $T_{\mathrm{pre}}$ & Mem. & $T_\mathrm{corr}$ & $T_{\mathrm{app}}$ & $\%$fmm & $N / T_{\mathrm{app}}$ \\
        \bottomrule
        \multicolumn{11}{c}{$n = 2$}\\
        \toprule
        $5.00$e$-2$ & $17832$ & $7.2$e$-1$ & $7.4$e$-1$ & $1.3$e$-2$ & $1.5$e$+0$ & $1.8$e$+0$ & $5.4$e$-5$ & $8.2$e$-2$ & $99.9$ & $1.7$e$+5$\\
        $2.50$e$-2$ & $70656$ & $1.3$e$+0$ & $2.7$e$+0$ & $6.9$e$-2$ & $4.1$e$+0$ & $7.0$e$+0$ & $1.9$e$-4$ & $3.8$e$-1$ & $99.9$ & $1.9$e$+5$\\
        $1.25$e$-2$ & $280218$ & $2.7$e$+0$ & $1.1$e$+1$ & $2.6$e$-1$ & $1.4$e$+1$ & $2.8$e$+1$ & $8.3$e$-4$ & $1.5$e$+0$ & $99.9$ &  $1.9$e$+5$\\
        $6.25$e$-3$ & $1118412$ & $6.2$e$+0$ & $4.3$e$+1$ & $1.2$e$+0$ & $5.1$e$+1$ & $1.1$e$+2$ & $4.2$e$-3$ & $6.0$e$+0$ & $99.9$ & $1.9$e$+5$\\
        $3.13$e$-3$ & $4464126$ & $1.5$e$+1$ & $1.7$e$+2$ & $4.5$e$+0$ & $1.9$e$+2$ & $4.4$e$+2$ & $1.9$e$-3$ & $2.4$e$+1$ & $99.9$ & $1.8$e$+5$\\
        $1.60$e$-3$ & $17846688$ & $5.5$e$+1$ & $6.9$e$+2$ & $1.8$e$+1$ & $7.7$e$+2$ & $1.8$e$+3$ & $1.1$e$-1$ & $9.1$e$+1$ & $99.9$ & $2.0$e$+5$ \\
        \bottomrule
        \multicolumn{11}{c}{$n = 4$}\\
        \toprule
        $1.00$e$-1$ & $11355$ & $1.0$e$+0$ & $1.2$e$+0$ & $2.2$e$-2$ & $2.3$e$+0$ & $2.7$e$+0$ & $6.2$e$-5$ & $5.4$e$-2$ & $99.9$ & $2.1$e$+5$\\
        $5.00$e$-2$ & $44580$ & $1.7$e$+0$ & $4.1$e$+0$ & $6.1$e$-2$ & $5.9$e$+0$ & $1.1$e$+1$ & $2.5$e$-4$ & $2.2$e$-1$ & $99.9$ & $2.0$e$+5$\\
        $2.50$e$-2$ & $176640$ & $3.6$e$+0$ & $1.7$e$+1$ & $4.6$e$-1$ & $2.1$e$+1$ & $4.2$e$+1$ & $1.1$e$-3$ & $9.1$e$-1$ & $99.9$ & $1.9$e$+5$\\
        $1.25$e$-2$ & $700545$ & $7.7$e$+0$ & $6.8$e$+1$ & $1.2$e$+0$ & $7.7$e$+1$ & $1.7$e$+2$ & $4.4$e$-3$ & $3.8$e$+0$ & $99.9$ & $1.9$e$+5$ \\
        $6.25$e$-3$& $2796030$ & $1.9$e$+1$ & $2.8$e$+2$ & $5.3$e$+0$ & $3.1$e$+2$ & $6.6$e$+2$ & $1.9$e$-2$ & $1.6$e$+1$ & $99.9$ & $1.7$e$+5$\\
        $3.13$e$-3$& $11160315$ & $4.9$e$+1$ & $1.1$e$+3$ & $2.1$e$+1$ & $1.2$e$+3$ & $2.6$e$+3$ & $9.6$e$-2$ & $6.6$e$+1$ & $99.9$ & $1.7$e$+5$\\
        \bottomrule
    \end{tabular}
    \endgroup
  \caption{Performance table for Laplace using FMM ($\epsilon_{\mathrm{fmm}} = 10^{-13}$); all experiments are performed on a single core, and $h$ values are chosen to result in a range of $N$ roughly from $10$ thousand to $10$ million. The offline throughput ($N/T_{\mathrm{pre}}$) at large $N$ is steady around $23,000$ for $n = 2$ and around $9200$ targets/core/second for $n = 4$. Memory usage of the volume correction matrix is displayed in MB.}\label{tab:fmm_laplace}
\end{table}

\begin{table}[h]
    \centering
    \begingroup
    \setlength{\tabcolsep}{4pt}
    \begin{tabular}{ c|c|c|c|c|c||c||c|c|c|c }
            $h$ & $N$ & $T_{\mathcal{B}_\Gamma, \mathrm{close}}$ & $T_{\mathcal{Q}_\Omega, \mathcal{B}_\Gamma, P_{\alpha}}$ & $T_C$ & $T_{\mathrm{pre}}$ & Mem. & $T_\mathrm{corr}$ & $T_{\mathrm{app}}$ & $\%$fmm & $N / T_{\mathrm{app}}$ \\
            \bottomrule
            \multicolumn{11}{c}{$n = 2$}\\
            \toprule
            $5.00$e$-1$ & $17382$ & $1.9$e$+0$ & $2.3$e$+0$ & $1.4$e$-2$ & $4.1$e$+0$ & $2.6$e$+0$ & $6.4$e$-5$ & $2.8$e$-1$ & $99.9$ & $6.1$e$+4$\\
            $2.50$e$-2$ & $70656$ & $3.6$e$+0$ & $7.6$e$+0$ & $6.3$e$-2$ & $1.1$e$+1$ & $1.0$e$+1$ & $2.6$e$-4$ & $1.0$e$+0$ & $99.9$ & $7.0$e$+4$\\
            $1.25$e$-2$ & $280218$ & $6.4$e$+0$ & $3.0$e$+1$ & $3.0$e$-1$ & $3.6$e$+1$ & $4.1$e$+1$ & $1.2$e$-3$ & $4.3$e$+0$ & $99.9$ & $6.5$e$+4$\\
            $6.25$e$-3$ & $1118412$ & $1.5$e$+1$ & $1.2$e$+2$ & $1.3$e$+0$ & $1.4$e$+2$ & $1.6$e$+2$ & $6.7$e$-3$ & $1.6$e$+1$ & $99.9$ & $7.1$e$+4$\\
            $3.13$e$-3$ & $4464126$ & $3.1$e$+1$ & $4.4$e$+2$ & $5.8$e$+0$ & $4.8$e$+2$ & $6.5$e$+2$ & $3.6$e$-2$ & $6.6$e$+1$ & $99.9$ & $6.8$e$+4$\\
            $1.60$e$-3$ & $17846688$ & $8.8$e$+1$ & $1.8$e$+3$ & $2.1$e$+1$ & $1.9$e$+3$ & $2.6$e$+3$ & $1.5$e$-1$ & $2.5$e$+2$ & $99.9$ & $7.1$e$+4$\\
            \bottomrule
            \multicolumn{11}{c}{$n = 4$}\\
            \toprule
            $1.00$e$-1$ & $11355$ & $3.5$e$+0$ & $5.0$e$+0$ & $2.8$e$-2$ & $8.6$e$+0$ & $4.0$e$+0$ & $9.4$e$-5$ & $2.5$e$-1$ & $99.9$ & $4.5$e$+4$ \\
            $5.00$e$-2$ & $44580$ & $6.6$e$+0$ & $1.6$e$+1$ & $1.0$e$-1$ & $2.3$e$+1$ & $1.6$e$+1$ & $3.4$e$-4$ & $8.3$e$-1$ & $99.9$ & $5.4$e$+4$ \\
            $2.50$e$-2$ & $176640$ & $1.3$e$+1$ & $6.3$e$+1$ & $5.5$e$-1$ & $7.7$e$+1$ & $6.2$e$+1$ & $1.5$e$-3$ & $3.6$e$+0$ & $99.9$ & $5.4$e$+4$\\
            $1.25$e$-2$ & $700545$ & $2.7$e$+1$ & $2.6$e$+2$ & $1.9$e$+0$ & $2.9$e$+2$ & $2.5$e$+2$ & $7.0$e$-3$ & $1.3$e$+1$ & $99.9$ & $5.3$e$+4$\\
            $6.25$e$-3$ & $2796030$ & $5.9$e$+1$ & $9.4$e$+2$ & $6.5$e$+0$ & $1.0$e$+3$ & $9.8$e$+2$ & $3.5$e$-2$ & $5.2$e$+1$ & $99.9$ & $5.4$e$+4$\\
            $3.13$e$-3$ & $11160315$ & $1.2$e$+2$ & $3.7$e$+3$ & $2.5$e$+1$ & $3.9$e$+3$ & $3.9$e$+3$ & $1.4$e$-1$ & $2.0$e$+2$  & $99.9$ & $5.5$e$+4$ \\
            \bottomrule
    \end{tabular}
    \endgroup
  \caption{Performance table for Helmholtz ($k = 10\pi$) using FMM ($\epsilon_{\mathrm{fmm}} =
  10^{-7}$); all experiments are performed on a single core, and $h$ values are chosen to result in a range of $N$ roughly from $10$ thousand to $10$ million. The offline throughput ($N/T_{\mathrm{pre}}$) at large $N$ is approximately $9300$ targets/core/second for $n = 2$ and approximately $3000$ targets/core/second for $n = 4$. Memory usage of the volume correction matrix is displayed in MB.}\label{tab:fmm_helmholtz}
\end{table}

\subsection{Lippmann-Schwinger equation}
To demonstrate the practical application of our proposed methodology, we employ it to solve the scattering problem of a planewave $u^{\rm inc}(\nex) = \exp(ik\nex\cdot\bol\theta)$, $|\bol\theta|=1$, interacting with a scatterer $\Omega\subset\R^2$ with the examples thus showing how the high-order accurate operator evaluation previously demonstrated can lead to a high-order volume integral equation solver. We consider a scatterer that has a smooth refractive index function $\eta>0$ defined within its boundaries. Outside of~$\Omega$, we assume a constant refractive index of $\eta=1$. The scattering problem seeks a volumetric total field $u:\R^2\to\C$, $u\in C^1(\R^2)\cap H^2_{\rm loc}(\R^2)$, that satisfies the Helmholtz equation
\begin{align}
    \Delta u+\eta(\nex)k^2 u=0&\quad\text{in}\quad\R^2\setminus\Gamma,
\end{align}
with the scattered field  $u^{\rm sc} \coloneqq u-u^{\rm inc}$ satisfying Sommerfeld's radiation condition. It can also be reformulated as the Lippmann-Schwinger volume integral equation for $u$ inside the scatterer:
\begin{equation}\label{eq:LS_eqn}
    u(\nex)+k^2 \mathcal{V}_{k}\left[(1-\eta) u\right](\nex)=u^{\mathrm{inc}}(\nex), \quad \nex \in \Omega,
\end{equation}
with the solution outside the scatterer given by the representation formula
\begin{equation}\label{eq:ext_LS_eqn}
u(\nex) = u^\inc(\nex)-k^2\mathcal V_k[(1-\eta)u](\nex),\quad\nex\in\R^2\setminus\overline\Omega.
\end{equation}

The Lippmann-Schwinger equation~\eqref{eq:LS_eqn} can be effectively solved using Nystr\"om methods based on the high-order discretization approaches presented in this work using iterative methods, wherein the resulting $N\times N$ linear system for the approximate values of the total field at the quadrature nodes $\{\boldsymbol{\xi}_j\}_{j=1}^N$ is solved via GMRES~\cite{saad1986gmres} and the volume integral operator $\mathcal{V}_k$ is used as a forward map. This methodology, which only necessitates repeated applications of the operator $\mathcal{V}_k$ over the fixed domain $\Omega$, allows us to exploit the algorithm's potential fully, as all source-independent steps can be precomputed and reused at each operator application (see \Cref{tab:complexity-estimates}).\enlargethispage*{3ex}
\begin{figure}[h!]
  \centering
    \begin{subfigure}[b]{\textwidth}
      \centering
      \includegraphics[width=0.35\linewidth]{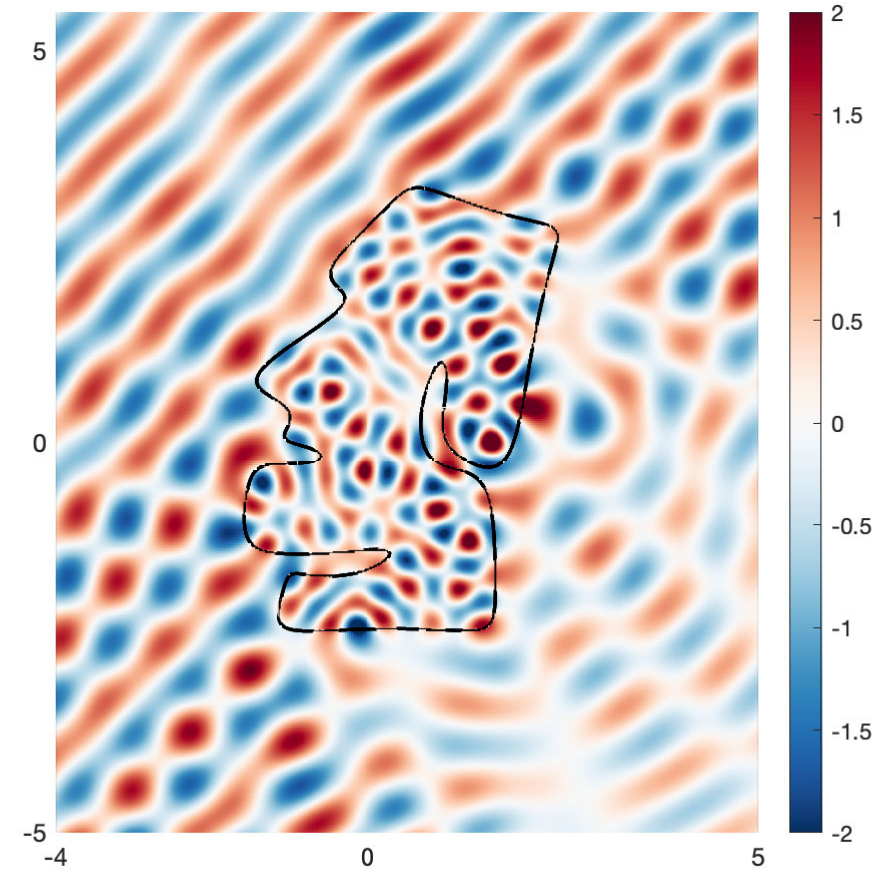}\qquad\qquad
      \includegraphics[width=0.35\linewidth]{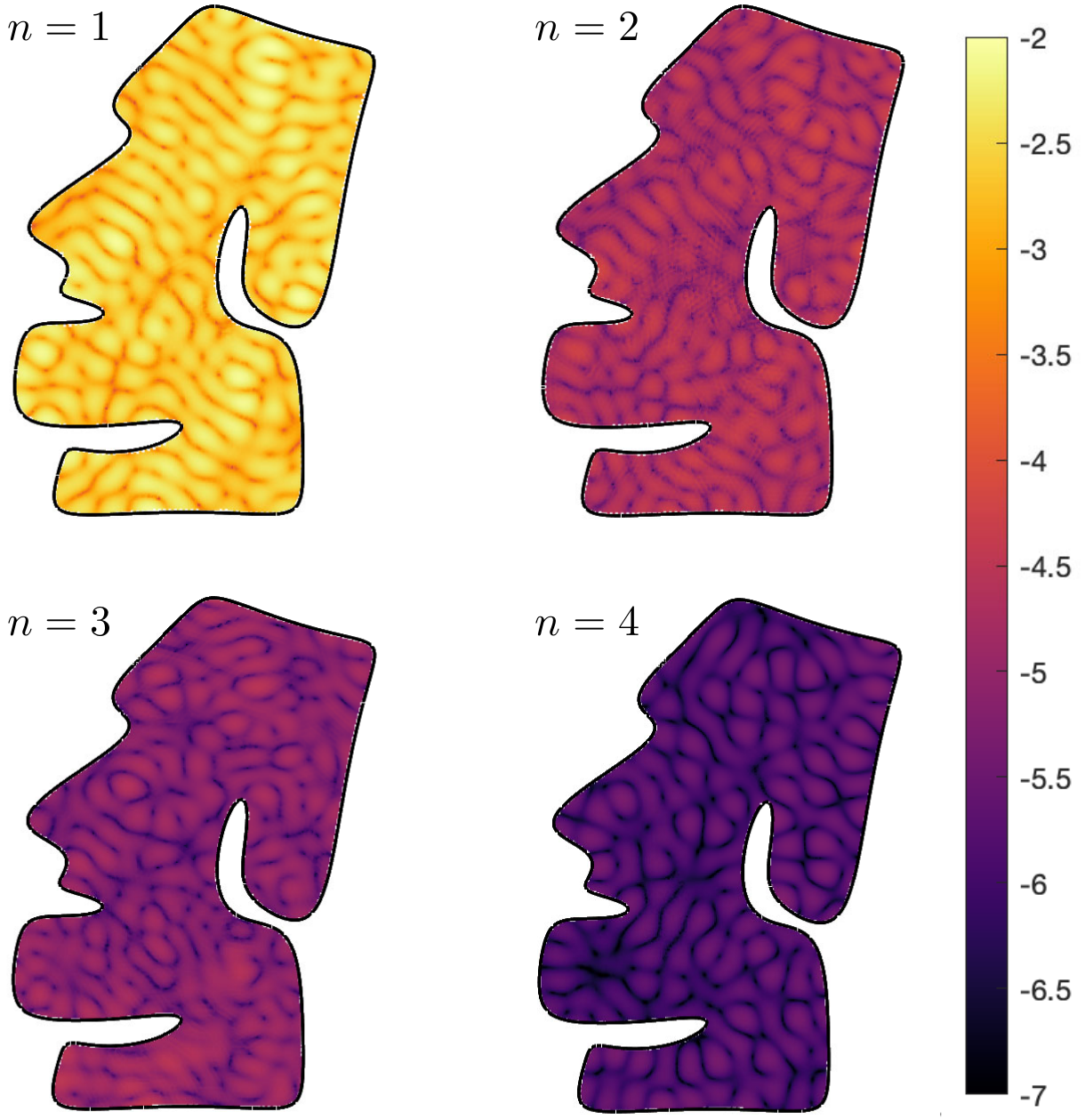}
      \caption{Planewave scattering by a homogeneous obstacle, for $k=2\pi$ and $\eta=4$ in $\Omega$. Left: Real part of the total field $u$ solution of the transmission problem~\eqref{eq:trans_problem}; nearly-singular volume integration arising in near-field evaluation of the representation formula~\cref{eq:ext_LS_eqn} is handled seamlessly by VDIM per \Cref{rem:at_boundary} Right: Logarithm in base ten of the absolute pointwise error in the Lippmann-Schwinger equation solution for interpolation degrees $n=1,2,3$ and $4$ over a fixed mesh of $\Omega$ with $h\approx 0.09$. The maximum poinwise errors in those figures are around $3.0\cdot 10^{-3}$, $2.3\cdot 10^{-5}$, $8.8\cdot 10^{-6}$, and $9.2\cdot 10^{-7}$, respectively.}
      \label{fig:LS_error_exmaple}
  \end{subfigure}
  \begin{subfigure}[b]{\textwidth}
      \centering
      \includegraphics[width=0.35\linewidth]{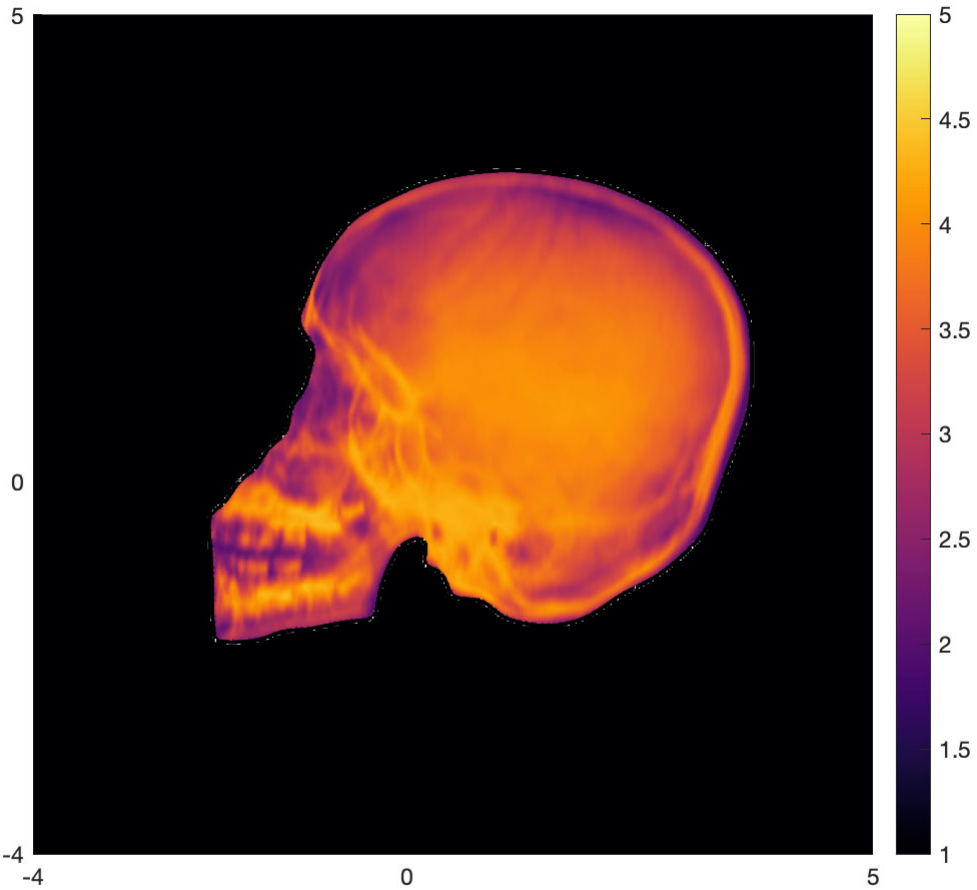}\qquad\qquad
      \includegraphics[width=0.35\linewidth]{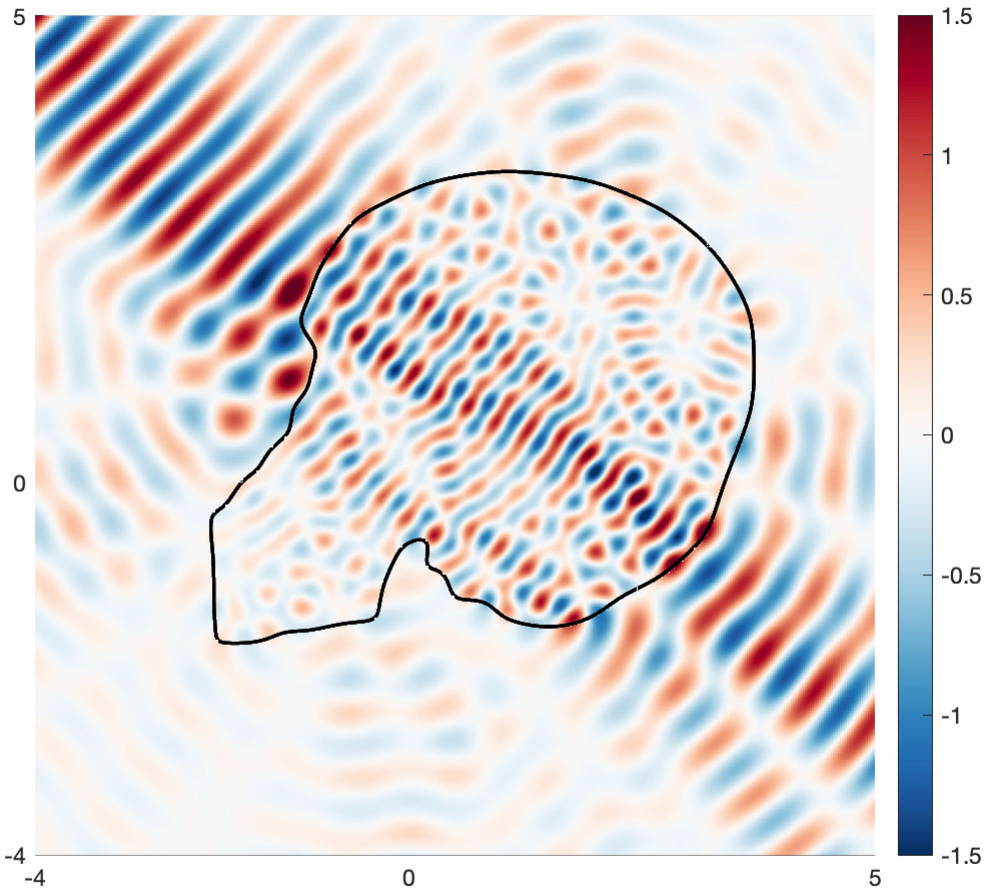}
      \caption{Gaussian beam scattering by an inhomogeneous obstacle, for $k=3\pi$ and a variable refractive index. Left: Piecewise smooth refractive index~$\eta$ generated from a skull x-ray image. Right: Real part of the total field, resulting from the scattering of a Gaussian beam by the skull-like inhomogeneous obstacle, obtained from an approximate Lippmann-Schwinger equation solution.}
      \label{fig:scattering_exmaple}
  \end{subfigure}
     \caption{Application of the proposed polynomial density interpolation method to the Lippmann-Schwinger equation. Convergence experiments demonstrate the same high convergence orders for solution of the volume integral equation that are illustrated for the operator evaluation in Figure~\ref{fig:tris_splines}, and to our knowledge represent the first reported provably high-order accurate solutions of such integral equations with piecewise-smooth material parameters.}
     \label{fig:lippmann-schwinger}
\end{figure}

To assess the accuracy of such a Nystr\"om method for the Lippmann-Schwinger equation, we first focus on a specific scenario where $\eta$ maintains a constant value of $4$ within $\Omega$. In this case, the total field can be determined by solving the transmission problem
\begin{equation}\begin{aligned}
\Delta u +\eta k^2 u=0 \ \text { in } \Omega,&\qquad \Delta u^{\rm sc}+k^2 u^{\rm sc}=0\ \text { in } \R^2\setminus\overline\Omega,\\
u=u^{\rm inc}+u^{\rm sc}\quad &\text{and}\quad \partial_{\nu} u=\partial_{\nu} \left(u^{\rm inc}+u^{\rm sc}\right) \text { on } \Gamma.
\end{aligned}\label{eq:trans_problem}\end{equation}
that can itself be recast into a second-kind boundary integral equation~\cite{kress1978transmission}. By solving the boundary integral equation (and employing an exponentially convergent Nystr\"om method, based on the Martensen-Kussmaul quadrature rule~\cite[Sec. 3.5]{COLTON:2012} for kernels with logarithmic singularities to do so), we obtain a highly accurate reference solution that serves as the benchmark to quantify the error in the solution of the volume integral equation for this piecewise constant refractive index medium. \Cref{fig:LS_error_exmaple} displays the real part of the total field $u$ in this example as well as the pointwise errors within $\Omega$ in the volume integral equation solutions obtained utilizing the proposed method for interpolation degrees $n=1,2,3$ and $4$ for a fixed mesh with $h\approx 0.09$. The resulting linear system was solved via GMRES with a tolerance of $10^{-7}$, which was achieved in approximately the same number of iterations ($\sim 340$) in all the examples presented in that figure.

We consider, finally, a challenging scattering problem involving a piecewise-smooth refractive index $\eta$ generated from an x-ray image of a human skull, seen in the left panel of \Cref{fig:scattering_exmaple}. The right panel of that figure shows the real part of the total field resulting from the scattering of a Gaussian beam that impinges on the skull obstacle from the top-left corner. The Lippmann-Schwinger integral equation solution displayed in that figure, was obtained using the proposed method with $n=2$. The scatterer diameter is about $9\lambda$, where the wavelength is $\lambda=2\pi/k=2/3$ in this example. GMRES convergence was achieved after $753$ iterations for a relative error tolerance of $10^{-7}$.  In order to generate a reference solution to roughly estimate the error, the same problem was solved over the same mesh ($h=0.11$) but with $n=4$. The relative error, measured on a rectangular curve enclosing the obstacle, is approximately $10^{-4}$.

\section{Limitations and Future Work\label{sec:limitations}}
The scheme developed in this paper is not meant to be directly applied as-is on arbitrarily-large domains or for sources with arbitrarily-fine features leading to highly anisotropic meshes; this section explores the causes of these limitations and proposes means of managing them. Unfortunately, the ill-conditioning issues that arise when directly solving the Vandermonde system on the physical element $K$, which were in principle avoided by using the translation and scaling strategy, manifest in a different form as numerical instabilities when evaluating~\eqref{eq:sep_interp_tri} and the translation formula~\eqref{eq:trans_pol_tri}. The formulas in~\eqref{eq:sep_interp_tri} may encounter classical cancellation errors typical in polynomial evaluation, and, additionally, significant numerical errors can arise in computation of the coefficients $\{c_\alpha[f](K)\}_{|\alpha|\leq n}$ arising from the translation formula in~\eqref{eq:trans_pol_tri}. The latter may indeed amplify the errors in the coefficients $\{\widetilde c_{\beta}[f](K)\}_{|\beta|\leq n}$ to a considerable extent and amounts to the main limitation of the proposed methodology.

To estimate the extent of this undesired amplification effect we can bound the norm of the linear transformation that maps $\{\widetilde c_\beta[f](K)\}_{|\beta|\le n}$ to $\{c_\alpha[f](K)\}_{|\alpha|\le n}$. In view of~\eqref{eq:trans_pol_tri},  this mapping is represented by the matrix $\bold P_K\in\R^{q_n\times q_n}$ with coefficients
\[
(\bold P_K)_{\alpha,\beta} := \begin{cases}\displaystyle\frac{p_{\beta-\alpha}(-\bol c_K)}{r_K^{|\beta|}},&\beta\geq \alpha,\\
0,&\beta<\alpha,\end{cases}
\]
where $\alpha=(\alpha_1,\alpha_2)$ and $\beta=(\beta_1,\beta_2)$ are multi-indices satisfying $|\alpha|=\alpha_1+\alpha_2\le n$ and $|\beta|=\beta_1+\beta_2\le n$.
In order to estimate the matrix norm $\|{\bold P_K}\|_{\infty}$,  we note from the definition of $p_{\beta-\alpha}$ and the inequality $|\bol c_K^{\beta-\alpha}|\leq |\bol c_K|^{|\beta-\alpha|}=|\bol c_K|^{\beta_1-\alpha_1}|\bol c_K|^{\beta_2-\alpha_2}$, that
\[
\sum_{\beta:|\beta|\leq n}|(\bold P_K)_{\alpha,\beta}| = \sum_{\substack{\beta:\alpha\leq\beta,\\|\beta| \le n}}\frac{|\bol c_K^{\beta-\alpha}|}{r_K^{|\beta|}(\beta-\alpha)!}
\leq\sum_{\beta_1=\alpha_1}^n\sum_{\beta_2=\alpha_2}^{n-\beta_1}\frac{|\bol c_K|^{\beta_1-\alpha_1}}{r_K^{\beta_1}(\beta_1-\alpha_1)!}\frac{|\bol c_K|^{\beta_2-\alpha_2}}{r_K^{\beta_2}(\beta_2-\alpha_2)!}.
\]
Translating the summation indices by $\alpha=(\alpha_1,\alpha_2)$ and restricting our attention to the relevant case $|\bol c_K|/r_K>1$, we then obtain
\begin{equation*}\begin{split}
        \sum_{\beta:|\beta|\leq n}|(\bold P_K)_{\alpha,\beta}|\leq&~\frac{1}{r_K^{|\alpha|}}\sum_{\eta_1=0}^{n-\alpha_1}\frac{|\bol c_K|^{\eta_1}}{r_K^{\eta_1}\eta_1!}\sum_{\eta_2=0}^{n-\eta_1-|\alpha|}\frac{|\bol c_K|^{\eta_2}}{r_K^{\eta_2}\eta_2!}\leq \frac{1}{r_K^{|\alpha|}}\sum_{\eta_1=0}^{n-\alpha_1}\frac{|\bol c_K|^{n-|\alpha|}}{r_K^{n-|\alpha|}\eta_1!}\sum_{\eta_2=0}^{n-\eta_1-|\alpha|}\frac{1}{\eta_2!}\\
\leq&~\e^2\frac{|\bol c_K|^{n-|\alpha|}}{r_K^{n}}.
\end{split}
\end{equation*}
It therefore  follows from the bound above that
\begin{equation}\label{eq:norm_translation}
\|\bold P_K\|_{\infty} =\max_{|\alpha|\leq n}\sum_{\beta:|\beta|\leq n}|(\bold P_K)_{\alpha,\beta}|\leq \e^2\begin{cases}\displaystyle\left(\frac{|\bol c_K|}{r_K}\right)^{n},& |\bol c_K|\geq 1,\\
r_K^{-n},&|\bol c_K|<1.\end{cases}
\end{equation}

\begin{figure}[h]
    \centering
    \includegraphics[width=0.43\textwidth]{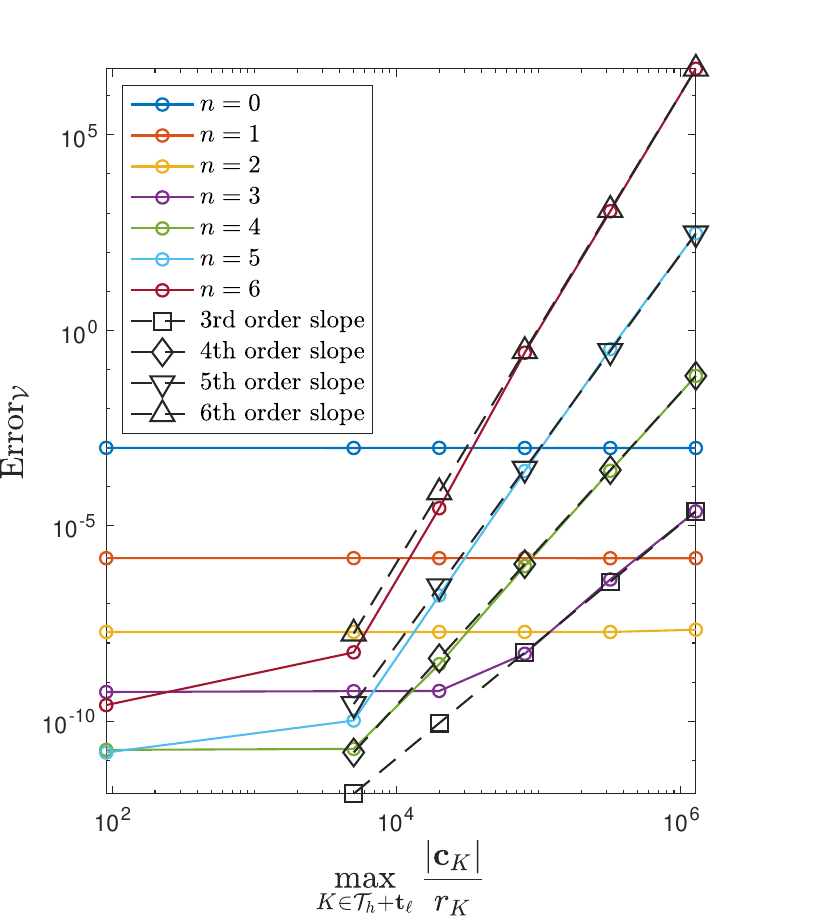}\quad
    \includegraphics[width=0.4\textwidth]{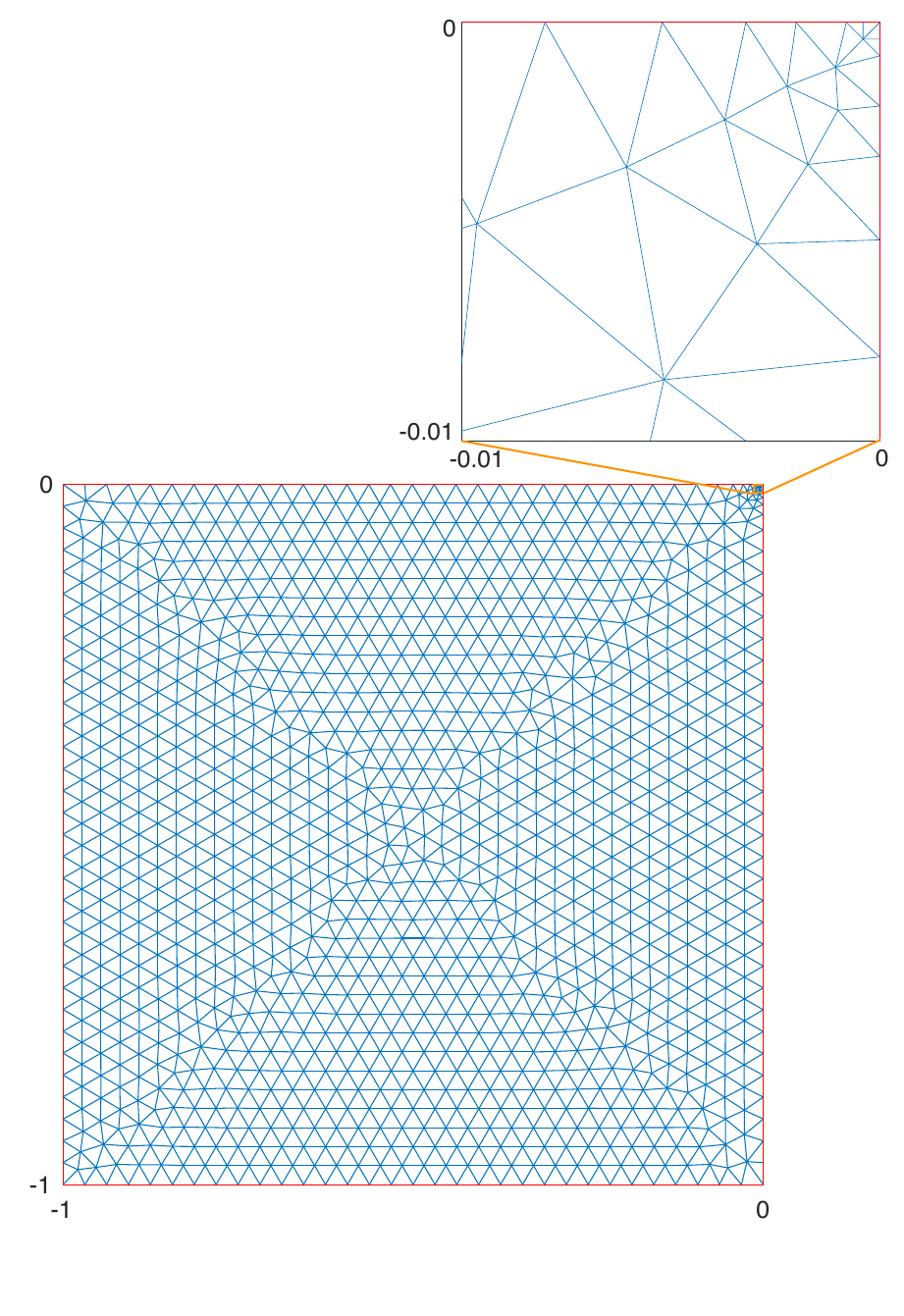}
    \caption{Effect of the translation formula~\eqref{eq:trans_pol_tri} on the errors in the evaluation of the operator $\mathcal V_k$. Left: Relative errors in the evaluation of $\mathcal V_k$ integrating over translations of the square domain $(-1,0)\times(-1,0)$. Right: Locally refined mesh of the square domain used to resolve the peaked density function~\eqref{eq:peakedSource} utilized in the corresponding numerical experiments.}\label{fig:stabilityExperiment}
\end{figure}

This analysis provides a rough indication of the impact of the error amplifying effect of the translation formula on the overall error in the proposed methodology and it reveals a critical quantity on which the error for a given domain discretization depends, namely, $\max_{K\in\mathcal T_h}|\bol c_K|/r_K$. To examine this in more detail we consider a numerical experiment where ${\rm Error}_{\mathcal V}$ is computed for integral operators $\mathcal V_k$ integrating over  a sequence of translated square domains $\Omega_\ell=(-1,0)\times(-1,0)+\bol t_\ell$ where $\bol t_0=(0,0)$ and $\bol t_\ell=2^{\ell-1}(1,1)$ for $\ell=1,\ldots,5$, which are applied to  densities  $f(\nex-\bol t_\ell)$, where $f$ is given by
\begin{equation}\label{eq:peakedSource}
f(\nex) = \mathrm{e}^{i k_0 \boldsymbol{x} \cdot \boldsymbol{\theta}}+\frac{1}{\left(k^2-4s\right)} \left(4 s^2\left|\boldsymbol{x}\right|^2-4 s+k^2\right) \mathrm{e}^{-s\left|\boldsymbol{x}\right|^2},
\end{equation}
for the parameters $ s=100$, $k=2\pi$, $k_0=3\pi$, and $\bol\theta=(\cos\frac{\pi}2,\sin\frac{\pi}3)$. To effectively handle the peaked Gaussian source component in $f$, which is centered at the square's top-right corner, the mesh $\mathcal{T}_h$ of $\Omega_0$ is locally refined around $(0,0)$ resulting in a ratio $\max_{K\in\mathcal{T}_h} h_K / \min_{K\in\mathcal{T}_h} h_K \approx (4.12 \times 10^{-2}) / (5.65 \times 10^{-4}) = 72.92$ (the meshed domain is displayed on the right panel in \Cref{fig:stabilityExperiment}). The remaining meshes of the domains $\Omega_\ell$, $\ell = 1, \ldots, 5$, are constructed by a direct translation of $\mathcal{T}_h$, i.e., as $\mathcal{T}_h + \bol{t}_\ell = \bigcup\{K + \bol{t}_\ell : K \in \mathcal{T}_h\}$. The plot in  \Cref{fig:stabilityExperiment} displays ${\rm Error}_{\mathcal{V}}$ for each one of the domains $\Omega_\ell$ versus the parameter $\max_{K\in \bol{t}_\ell + \mathcal{T}_h}  |\mathbf{c}_K| / r_K$ that, in view of~\eqref{eq:norm_translation}, appears to be relevant to assess the impact of translation formula on the error. These results show that for sufficiently large $\max_{K\in \bol{t}_\ell + \mathcal{T}_h}  |\mathbf{c}_K| / r_K$ values, the error in evaluating the interpolation polynomial coefficients is amplified to the extent that it becomes dominant, leading to ${\rm Error}_{\mathcal{V}} \propto \max_{K\in \bol{t}_\ell + \mathcal{T}_h}  (|\mathbf{c}_K| / r_K)^n$ as expected from the estimate~\eqref{eq:norm_translation}. This verifies that the translation formulae limits the achievable accuracy of the proposed methodology, especially for interpolation degrees $n>4$. Indeed, these results show that even in the most benign case, $\Omega_0$, where $\max_{K\in\mathcal T_h}|\bol c_K|\leq 1$, for interpolation degrees $n>4$ the overall error  is dominated by the amplified errors in the coefficients of the interpolation polynomial in the monomial basis, a phenomenon also revealed by the estimate~\eqref{eq:norm_translation}.

The analysis and experiments presented above suggest a strategy to control the accuracy deterioration due to the use of the translation formula, at least to some extent. Such strategy consists of decomposing the mesh $\mathcal T_h$ as $\bigcup_{j=1}^D \mathcal T^{(j)}_h$, and  find $\bol t_j$ that minimizes $\max_{K \in \mathcal{T}^{(j)}_h} | \bol{c}_K - \bol t_j | / r_K$ for $j=1,\ldots,D$. This min-max problem, however, can be difficult to solve in practice, but one may instead choose to use a heuristic for the location of the translation point $\bol t_j$ to approximate it as a weighted center of mass of $\mathcal{T}^{(j)}_h$ that accounts for the $1/r_K$ factor, i.e., $\bar{\bol t}_j = \left(\sum_{K\in\mathcal T_h^{(j)}}\bol c_K/r_K\right)/\left(\sum_{K\in\mathcal T_h^{(j)}}1/r_K\right)$. Provided $D$ is sufficiently large (i.e. there are enough subdomains), one can manage the size of $\max_{1 \leq j \leq D} \max_{K\in \mathcal{T}^{(j)}_h} \bol |\bol c_K - \bol t_j|/r_K$,
thus alleviating the error-amplifying effect of the translation formula stemming from the resulting $(|\bol c_K-\bol t_j|/r_K)^n$ factor in~\eqref{eq:norm_translation}. In the limit case when $\mathcal T_h^{(j)}=\{K_j\}$, $j=1,\ldots,D=L$, this domain decomposition procedure reduces to the methodology put forth in~\cite{ShenSerkh:22} which is known not to suffer from the aforementioned instabilities. It is worth noticing, though, that the experiments in \Cref{sec:numer} show that, in practice, for interpolation degrees $n\leq 4$, one can achieve $\mathrm{Error}_{\mathcal V} \approx 10^{-11}$ without any rescaling or domain decomposition, for domains exhibiting both $\max_{K\in\mathcal T_h}|\bol c_K|/r_K$ and $\max_{K\in\mathcal T_h}r^{-1}_K$ well above $5\cdot 10^3$.

In fact, the error analysis of regularized volume integrals developed here justifies a slightly modified alternative approach where the regularization is done more locally, with the boundary of a regularization region laying on a boundary of the union of elements at a distance sufficiently far from the evaluation point so that those boundary integrals can be treated with regular quadrature, while the volume integrals in the complementary region are sufficiently smooth so as to not require regularization. As a byproduct of the locality, no fast algorithms would be needed in the off-line phase of such a variant. We are currently pursuing this idea, and anticipate that this will allow the positive features of the method to seamlessly hold in settings with extreme mesh anisotropy or in large domains with no additional cost.

\section{Conclusions}\label{sec:conclusions}
We have presented a provably high-order scheme for the numerical evaluation of volume integral operators that is amenable to fast algorithms. Future work may consider the macro forms of domain subdivision mentioned in \Cref{sec:limitations} and their effect on computational efficiency. Applications to other PDEs will be explored in future work and is in some cases quite natural given the generality afforded by~\cite{faria2021general,anderson2023particular}. The reduction in the dimensionality of the region of singular quadrature should prove especially effective in three dimensions. We look forward to more finely tuned and parallelized schemes built on these concepts and coupling them to more adaptive methods.

\appendix \section{Invertibility and conditioning of the translated-and-scaled Vandermonde matrix\label{sec:conditioning}}
This appendix addresses the invertibility and conditioning of the multivariate Vandermonde matrix $\widetilde{\bf{V}}:=\widetilde{\bf{V}}_K$ in~\cref{eq:vandermonde_shifted} and~\cref{Vandermonde_shiftedscaled_matrix} corresponding to a set of distinct points $\widetilde{\mathcal I}_n$. It is well-known~\cite{Sauer:95,Sauer:00,Olver:06} that there is a subtle (and increasingly so at high order) interplay between the geometry of the nodal set and the solvability of the interpolation problem (nevertheless, this is often still the case~\cite{Sauer:95}). The poisedness for node sets transformed from known-poised sets is, in our case, a consequence of reference~\cite[Lem.\ 1]{Ciarlet:72}, by an argument which is partially used in the proof of \Cref{eq:conditioning_corr} below. However, until the effect of the transformation on the node-sets is understood, \emph{the conditioning could in principle be arbitrarily poor}.

\begin{figure}[h]
    \centering
    \includegraphics[width=0.6\textwidth]{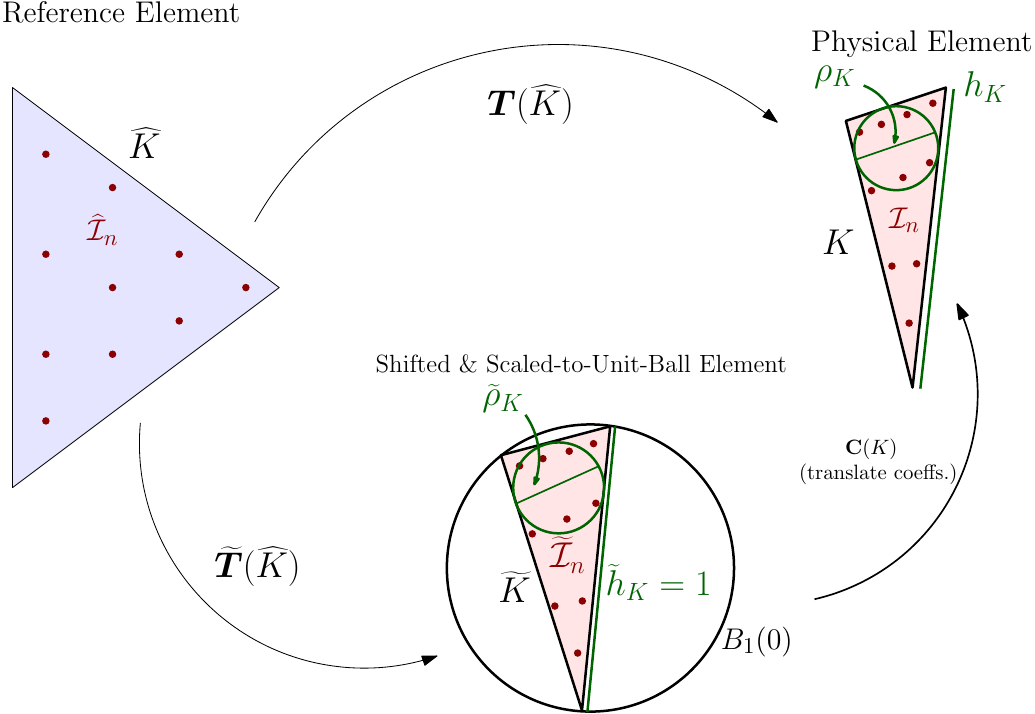}
    \caption{Coordinate systems and transformations. The condition number of the Vandermonde-like system for interpolation in the unit ball (lower) can be controlled.}\label{fig:shiftscale0}
\end{figure}

To address these matters, drawing in part on results in~\cite{Conn:08}, in \Cref{eq:conditioning_corr} below we show that the linear system in~\cref{eq:vandermonde_shifted} \emph{for an arbitrary element $\widetilde K$} generated by an affine transformation $\widetilde{\vv T}$, is not only necessarily invertible ($\widetilde{\mathcal{I}}_n$ is poised) but has a condition number that can be bounded above by a function that depends explicitly on a standard measure of triangle quality and the degree $n$ of the polynomial basis. For simplicity, we restrict the theorem statement to straight elements; a similar statement could be developed for curved elements generated by (or approximated by) polynomial transformations. It will be useful to consider the geometrical parameters (see \Cref{fig:shiftscale0})
\begin{equation}
    \widetilde h_K \coloneqq \operatorname{diam}\widetilde{K}\quad\mbox{and}\quad \widetilde\rho_K \coloneqq \sup_{a>0} \left\{2a: B_a(\nex) \subset \widetilde{K}, \nex \in \R^2 \right\},
\end{equation}
that denote the diameter (length of longest edge) of the triangle and the maximum diameter of inscribed disks $B_a(\nex) \coloneqq \{\ney\in\R^2: |\nex-\ney|<a\}$ in the triangle, respectively. They, together through their ratio, provide a common measure of triangle quality. Note that by construction we have $\sqrt 3\leq \widetilde h_K\leq 2$; the parameters $\widehat h_K=\sqrt{3}$ and $\widehat\rho_K=1$ denote analogous quantities for the reference triangle $\widehat{K}$.

In order to state and prove the result, it will help to recall a few facts from the work~\cite{Conn:08}. First, a given nodal interpolation set $\mathcal Y = \{\vv{\zeta}_1, \vv{\zeta}_2, \ldots, \vv{\zeta}_{q_n}\} \subset B_1(\bol 0) \subset \R^2$ is said to be \emph{$\Lambda$-poised} if and only if the Euclidean norm $\| \vv{\lambda}(\cdot) \|_2$ of the vector function $\vv{\lambda}:\R^2\to\R^{q_n}$ with each element a member of the Lagrange polynomial basis for that nodal set $\mathcal Y$, satisfies $\sup_{\vv{\zeta} \in B_1(\bol 0)}\| \vv{\lambda}(\vv\zeta) \|_2 \le \Lambda$. The quantity $\sup_{\vv{\zeta} \in B_1(\bol 0)}\| \vv{\lambda}(\vv\zeta) \|_2$ can be bounded above and below by constants times the Lebesgue constant of common parlance, and so in effect $\Lambda$ functions as such. Then, using the specific factorial-normalized monomial basis set $\{p_\alpha\}_{|\alpha|\leq n}$, $p_\alpha(\ney) = \ney^\alpha / \alpha!$, the work~\cite{Conn:08} develops a \emph{bidirectional} connection between the condition number of the multivariate Vandermonde matrix $\bf{V}$ (for that basis and the nodal set $\mathcal Y$) and the size of Lagrange polynomials associated with $\mathcal Y$, i.e. the quantity $\Lambda$ for which $\mathcal Y$ is $\Lambda$-poised.

We repeatedly use the following lemma, which is a direct adaptation of~\cite[Thm.\ 1]{Conn:08} to our notation; essentially, it represents one direction of the just mentioned \emph{bidirectional} connection.
\begin{figure}[h!]
  \centering
      \includegraphics[width=0.9\linewidth]{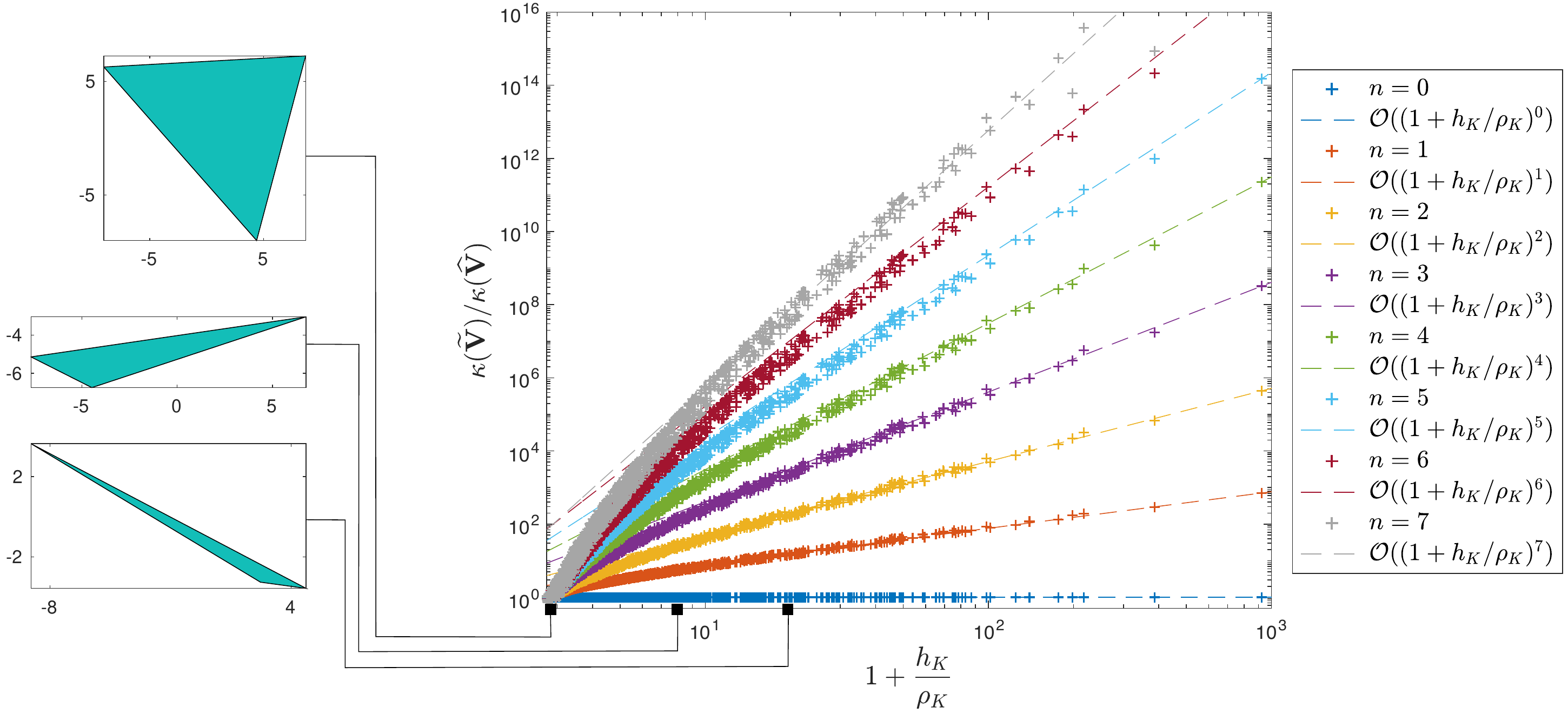}
     \caption{Dependence of the condition number $\kappa(\widetilde{\bf V})$  on the interpolation degree $n$ and the triangle's shape characterized by the ``aspect ratio" $h_K/\rho_K$. Each of the data points marked by a colored plus sign (+)  corresponds to a triangle $K$ with randomly generated vertices contained in $[-10,10]^2$. A total of 1000 triangles were used in this example, each tested for every $n$ in the study. The dashed lines mark the asymptotic behavior of the bound for $\kappa(\widetilde{\bf V})/\kappa(\widehat{\bf V})$  established in~\Cref{eq:conditioning_corr}. A remarkable agreement between the empirical data and theoretical bound is observed.}\label{fig:random_triangles}
\end{figure}
\begin{lemma}\label{lem:Conn}
If ${\bf V}$ is nonsingular and $\left\|{\bf{V}}^{-1}\right\| \leq M$, then the set $\mathcal{Y}$ is $ \Lambda$-poised with $\Lambda = \sqrt{q_n}M$ in the unit disk $B_1(\bol 0)$. Conversely, if the set ${\mathcal{Y}}$ is $\Lambda$-poised in the unit disk $B_1(\bol 0)$, then ${\bf{V}}$ is nonsingular and
\[
\left\|{\bf{V}}^{-1}\right\| \leq \theta_n \Lambda
\]
where $\theta_n>0$ is dependent on $n$ but independent of ${\mathcal{Y}}$ and $\Lambda$.
\end{lemma}

\begin{theorem}\label{eq:conditioning_corr}
    For a given $n\in\N$, let the nodes $\mathcal{I}_n$ be mapped from the Vioreanu-Rokhlin nodes $\widehat{\mathcal{I}}_n$ of interpolation degree $n$ via a non-degenerate affine transformation $\vv T$. We have that the interpolation sets $\mathcal I_n$ and $\widetilde{\mathcal I}_n$ are poised in the normalized monomial basis $\{p_\alpha\}_{|\alpha|\leq n}$. Further, the associated multivariate Vandermonde matrix~$\widetilde{\bf{V}}$ in~\cref{Vandermonde_shiftedscaled_matrix} is invertible and its condition number satisfies the bound
\begin{equation}\label{eq:condition_estimate}
  \kappa(\widetilde{\bf{V}}) \le \chi_n \left(1 + \frac{h_K}{\rho_K}\right)^n \kappa(\widehat{\bf{V}}),
\end{equation}
where $\chi_n$ is computable and depends only on $n$, and where $\widehat{\bf{V}}$ is the multivariate Vandermonde matrix in the normalized basis on $\widehat{\mathcal{I}}_n$.
\end{theorem}

To assess the sharpness of the estimate~\eqref{eq:condition_estimate}, we present \Cref{fig:random_triangles} which displays the quotient  $\kappa(\widetilde{\bf V})/\kappa(\widehat{\bf V})$ plotted against $1+h_K/\rho_K$ for matrices $\widetilde{\bf V}$ corresponding to 1000 triangles $K$ with random vertices uniformly distributed within $[-10,10]^2$ and for various interpolation degrees. These results suggest that our bound~\eqref{eq:condition_estimate} effectively captures the dependence of the condition number $\kappa(\widetilde{\bf V})$ on the triangle's shape characterized here by the ``aspect ratio" $h_K/\rho_K$.

\begin{remark}\label{rem:conditioning_corr}
  Explicit bounds on the ``constant" $\theta_n$ in~\Cref{lem:Conn} can be computed; reference~\cite[Eq.\ (13)]{Conn:08} provides the necessary ingredient for $n=2$, showing a path for bounding $\theta_n$ for larger $n$ as well. Concretely, while the quantity $\kappa(\widehat{\bf{V}})$ can be computed once and for all in the basis $\{p_\alpha\}_{|\alpha|\leq n}$ and is of very small size, the theorem represents a somewhat modest result because $\chi_n$ in~\Cref{eq:conditioning_corr} depends linearly on $\theta_n$ in~\Cref{lem:Conn}, and the latter quantity admittedly, as reference~\cite{Conn:08} notes, likely grows as much as exponentially with $n$.
\end{remark}

\begin{proof}
As a first step, we estimate the ${\Lambda}$-poisedness of $\widehat{\mathcal I}_n$. It follows from~\Cref{lem:Conn}, the identity $\|\widehat{\bf{V}}^{-1}\|=\kappa(\widehat{\bf{V}})\|\widehat{\bf{V}}\|^{-1}$, and the fact that $\|\widehat{\bf {V}}\| \ge 1$, that $\widehat{\mathcal{I}}_n$ is $\widehat\Lambda$-poised with
\begin{equation}\label{eq:ref_lagrange_estimate_condition}
    \widehat{\Lambda} \coloneqq \sqrt{q_n} \kappa(\widehat{\bf{V}}),
\end{equation}
where $\kappa(\widehat{\bf{V}})$ denotes  the condition number of the multivariate Vandermonde matrix for the Vioreanu-Rokhlin nodes in the monomial basis.

We next study the effect of the mapping $\widetilde{\vv T}:\widehat K\to \widetilde K$ on the reference Lagrange-polynomial vector function, denoted $\widehat{\vv{\lambda}}=[\widehat\lambda_1,\ldots,\widehat\lambda_{q_n}]^\top:\R^2\to\R^{q_n}$, which via the concept of $\Lambda$-poisedness introduced above, ultimately provides control on the condition number of the multivariate Vandermonde system. We first note that the polynomials
\begin{equation}\label{lambdaj_tilde_def}
    \widetilde{\lambda}_j(\widetilde\nex) \coloneqq \widehat{\lambda}_j(\widetilde{\vv T}^{-1}(\widetilde{\nex})),
\end{equation}
in fact satisfy the Lagrange interpolation property $\widetilde \lambda_j(\widetilde{\vv{x}}_i) = \delta_{ij}$, $i,j =1,\ldots, q_n$. An analogous statement holds for $\lambda_j(\nex) \coloneqq \widehat{\lambda}_j(\vv T^{-1}(\nex))$. It follows that both $\widetilde{\mathcal{I}}_n$ and $\mathcal{I}_n$ are poised.

From~\cref{lambdaj_tilde_def} it follows that
\[
    \sup_{\widetilde\nex \in B_1(\bol 0)} \left| \widetilde{\lambda}_j(\widetilde\nex) \right| = \sup_{\widetilde\nex \in B_1(\bol 0)} \left| \widehat{\lambda}_j(\widetilde{\vv T}^{-1}(\widetilde\nex)) \right| \le \sup_{\widehat{\nex} \in B_{1+R}(\bol 0)} \left| \widehat{\lambda}_j(\widehat{\nex})\right|, \quad j=1,\ldots, q_n,
\]
where we used the fact that $\widetilde{\vv T}^{-1}(B_1(\bol 0))\subset B_{1+R}(\bol 0)$ where $R\coloneqq \widetilde{h}_K/\widetilde\rho_K = h_K / \rho_K$. This inclusion follows from the fact that $|\widehat\nex|=|\widetilde{\vv T}^{-1}(\widetilde \nex)|=|{\bf A}\widetilde\nex+\vv b|\leq \|{\bf A}\||\widetilde \nex|+|\vv b|\leq R+1$ which we establish next in two parts. Indeed, $\left\|\bf A\right\|\leq \widehat{h}_K/\widetilde\rho_K \leq R$, a fact that itself results from \cite[Lemma 2]{Ciarlet:72} that establishes the first inequality followed by use of the second inequality $\widehat{h}_K \le\widetilde{h}_K$ which in turn follows from the fact that an equilateral triangle $\widehat K$ is the triangle having a minimum diameter among all triangles with a minimum bounding circle of unit radius. The bound $|\vv b|\leq 1$ used above, on the other hand, follows directly from the bijectivity of  $\widetilde{\vv T}^{-1}:\widetilde K\to\widehat K$ together with the fact that, by construction, $\bol 0\in\widetilde K$, which imply $\vv b=\widetilde{\vv T}^{-1}(\bol 0)\in\widehat K\subset\overline{B_1(\bol 0)}$.

It remains to estimate the supremum of $|\widehat\lambda_j|$ over $B_{1+R}(\bol 0)$, for which purpose we build upon the bound, that follows from the definition of $\Lambda$-poisedness~\cref{eq:ref_lagrange_estimate_condition}, $|\widehat\lambda_j(\widehat{\nex})| \le\sqrt{q_n} \kappa(\widehat{\bf{V}})$ for $\widehat{\nex}$ in the more limited region $B_1(\bol 0)$. To extend the estimate to $B_{1+R}(\bol 0)$ we use a Markov-type inequality that provides bounds on homogeneous sub-components of a polynomial relative to the norm of the polynomial on a more limited convex region such as a disk. Specifically, reference~\cite[Thm.\ 1]{ganzburg2002markov} provides for a polynomial $Q(\nex) = \sum_{|\alpha| \le n} d_\alpha \nex^\alpha$, the bound
\begin{equation}\label{eq:homog_polynomial_extension_bound}
    \left|\sum_{|\alpha| = \ell} d_\alpha \nex^\alpha\right| \le B_\ell^{(n)} \left|\nex\right|^\ell \sup_{\ney \in B_1(\bol 0)} \left|Q(\ney)\right|,\quad \ell=0,\ldots, n,
\end{equation}
 where for a given $n$ the $(n+1)$-numbered constants $B_\ell^{(n)}$ are known explicitly (see~\cite{ganzburg2002markov}) and are related to the coefficients of Chebyshev polynomials of degree $n$. In our context, for a given $j\in\{1,\ldots,q_n\}$, choosing $Q (\widehat\nex)= \widehat{\lambda}_j(\widehat\nex)$ and summing over each of the homogeneous polynomials of degree $\ell=0,\ldots, n$ making up $\widehat{\lambda}_j$, we obtain
 $
\widehat\lambda_j(\widehat\nex) =\sum_{\ell=0}^n\sum_{|\alpha|=\ell}d_\alpha \widehat\nex^\alpha
$. Then, employing the triangle inequality, and using~\cref{eq:homog_polynomial_extension_bound}, we readily find
\[
\left|\widehat\lambda_j(\widehat\nex)\right| \leq \sum_{\ell=0}^n\left|\sum_{|\alpha|=\ell}d_\alpha \widehat\nex^\alpha\right|\leq \left(\sup_{\ney \in B_1(\bol 0)} \left|\widehat\lambda_j(\ney)\right|\right) \sum_{\ell=0}^n B_\ell^{(n)}|\widehat\nex|^\ell \leq  \left(\sup_{\ney \in B_1(\bol 0)} \left|\widehat\lambda_j(\ney)\right|\right)\left(\max_{\ell=0,\ldots,n}B^{(n)}_\ell\right)\sum_{\ell=0}^n |\widehat\nex|^\ell.
\]
Letting $B^{(n)} = \max_{\ell=0,\ldots, n} B_\ell^{(n)}$ and taking supremum of the expressions above over $\widehat\nex\in B_{1+R}(\bol 0)$ it hence follows that
\begin{equation*}
    \begin{split}
        \sup_{\widetilde\nex \in B_1(\bol 0)} \left|\widetilde\lambda_j(\widetilde\nex)\right| &\le \sup_{\widehat\nex \in B_{1+R}(\bol 0)} \left|\widehat\lambda_j (\widehat\nex)\right|
        \le 2 B^{(n)} (1+R)^n \sup_{\ney \in B_1(\bol 0)} \left|\widehat\lambda_j(\ney)\right|,\quad  j=1,\ldots, q_n,
    \end{split}
\end{equation*}
where we have used the fact that
$\sum_{\ell=0}^n |\widehat\nex|^\ell\leq\sum_{\ell=0}^n (1+R)^\ell=[(1+R)^{n+1}-1]/R\leq 2(1+R)^n
$ which is in turn obtained from $R=h_K/\rho_K>1$ and $\widehat\nex\in B_{1+R}(\bol 0)$. It follows immediately that the set $\widetilde{\mathcal I}_n$ is $\widetilde\Lambda$-poised with
\begin{equation}\label{eq:tildeLambda_poised_bound}
    \widetilde{\Lambda} \coloneqq  2B^{(n)} (1+R)^n \sqrt{q_n} \kappa(\widehat{\bf{V}})\geq 2 B^{(n)} (1+R)^n  \sup_{\ney \in B_1(\bol 0)} \|\widehat{\vv{\lambda}}(\ney)\|_2\geq \sup_{\widetilde{\nex} \in B_1(\bol 0)} \|\widetilde{\vv{\lambda}}(\widetilde{\nex})\|_2.
\end{equation}
     We have, further, from using~\Cref{lem:Conn} together with~\cref{eq:tildeLambda_poised_bound}, the bound
\begin{equation}\label{eq:tildeV_inv_norm}
    \left\|\widetilde{\bf V}^{-1}\right\| \le \theta_n \widetilde{\Lambda} = 2\theta_n  B^{(n)} (1+R)^n \sqrt{q_n}\kappa(\widehat{\bf{V}}).
\end{equation}
  Here $\theta_n$ denotes the quantity in~\Cref{lem:Conn} (see also \Cref{rem:conditioning_corr}).
  The result follows from this inequality in conjunction with the inequality, see~\cite[Eq.\ (9)]{Conn:08}, $\left\|\widetilde{\bf{V}}\right\| \le q_n^{3/2}$, which eventually yields $\chi_n = 2\theta_nq_n^{2}B^{(n)}$.
\end{proof}
\bibliographystyle{siamplain}
\bibliography{References}
\end{document}

%% file: shared.tex
\usepackage{lipsum}
\usepackage{amsfonts}
\usepackage{graphicx}
\usepackage{epstopdf}
\usepackage{algorithmic}
\ifpdf
  \DeclareGraphicsExtensions{.eps,.pdf,.png,.jpg}
\else
  \DeclareGraphicsExtensions{.eps}
\fi

\newsiamremark{remark}{Remark}
\newsiamremark{hypothesis}{Hypothesis}
\crefname{hypothesis}{Hypothesis}{Hypotheses}
\newsiamthm{claim}{Claim}

\headers{Fast, High-order Accurate Evaluation of Volume Potentials}{T. G. Anderson, M. Bonnet, L. M. Faria, and C. P\'erez-Arancibia}

\title{Fast, high-order numerical evaluation of volume potentials via
polynomial density interpolation\thanks{Submitted to the editors \today.
}}

\author{Thomas G. Anderson\thanks{Department of Computational Applied Mathematics \& Operations Research, Rice University, Houston, TX USA
  (\email{thomas.anderson@rice.edu}).}
\and Marc Bonnet\thanks{POEMS (CNRS, INRIA, ENSTA), ENSTA Paris, 91120 Palaiseau, France
  (\email{marc.bonnet@ensta-paris.fr}, \email{luiz.maltez-faria@inria.fr}).}
\and Luiz M. Faria\footnotemark[3]
\and Carlos P\'erez-Arancibia\thanks{Department of Applied Mathematics, University of Twente, Enschede, The Netherlands (\email{c.a.perezarancibia@utwente.nl}).}
}

\usepackage{amsopn}


%% file: main_v7_arxiv.bbl
\begin{thebibliography}{10}

\bibitem{af2018adaptive}
{\sc L.~af~Klinteberg and A.-K. Tornberg}, {\em Adaptive quadrature by
  expansion for layer potential evaluation in two dimensions}, SIAM Journal on
  Scientific Computing, 40 (2018), pp.~A1225--A1249.

\bibitem{Greengard:16}
{\sc S.~Ambikasaran, C.~Borges, L.-M. Imbert-Gerard, and L.~Greengard}, {\em
  Fast, adaptive, high-order accurate discretization of the
  {L}ippmann--{S}chwinger equation in two dimensions}, SIAM Journal on
  Scientific Computing, 38 (2016), pp.~A1770--A1787.

\bibitem{anderson2023particular}
{\sc T.~G. Anderson, M.~Bonnet, L.~M. Faria, and C.~P{\'e}rez-Arancibia}, {\em
  Construction of polynomial particular solutions of linear
  constant-coefficient partial differential equations}, Computers \&
  Mathematics with Applications, 162 (2024), pp.~94--103.

\bibitem{Inti}
{\sc T.~G. Anderson, L.~M. Faria, and C.~P{\'e}rez-Arancibia}, {\em Inti.jl}.
\newblock \url{https://github.com/IntegralEquations/Inti.jl}, 2024.

\bibitem{Anderson:22a}
{\sc T.~G. Anderson, H.~Zhu, and S.~Veerapaneni}, {\em A fast, high-order
  scheme for evaluating volume potentials on complex 2{D} geometries via
  area-to-line integral conversion and domain mappings}, Journal of
  Computational Physics, 472 (2023), p.~111688.

\bibitem{Atkinson:85}
{\sc K.~E. Atkinson}, {\em The numerical evaluation of particular solutions for
  {P}oisson's equation}, {IMA} Journal of Numerical Analysis, 5 (1985),
  pp.~319--338, \url{https://doi.org/10.1093/imanum/5.3.319}.

\bibitem{Atkinson}
{\sc K.~E. Atkinson}, {\em The Numerical Solution of Integral Equations of the
  Second Kind}, vol.~4, Cambridge University Press, 1997.

\bibitem{Averbuch:00}
{\sc A.~Averbuch, E.~Braverman, and M.~Israeli}, {\em Parallel adaptive
  solution of a {P}oisson equation with multiwavelets}, SIAM Journal on
  Scientific Computing, 22 (2000), pp.~1053--1086.

\bibitem{bao2024singularity}
{\sc G.~Bao, W.~Hua, J.~Lai, and J.~Zhang}, {\em Singularity swapping method
  for nearly singular integrals based on trapezoidal rule}, SIAM Journal on
  Numerical Analysis, 62 (2024), pp.~974--997.

\bibitem{barnett2015spectrally}
{\sc A.~Barnett, B.~Wu, and S.~Veerapaneni}, {\em Spectrally accurate
  quadratures for evaluation of layer potentials close to the boundary for the
  2d {S}tokes and {L}aplace equations}, SIAM Journal on Scientific Computing,
  37 (2015), pp.~B519--B542.

\bibitem{Barnett:14}
{\sc A.~H. Barnett}, {\em Evaluation of layer potentials close to the boundary
  for {L}aplace and {H}elmholtz problems on analytic planar domains}, SIAM
  Journal on Scientific Computing, 36 (2014), pp.~A427--A451.

\bibitem{Bauinger:21}
{\sc C.~Bauinger and O.~P. Bruno}, {\em ``{I}nterpolated factored {G}reen
  function'' method for accelerated solution of scattering problems}, Journal
  of Computational Physics, 430 (2021), p.~110095.

\bibitem{beale2016simple}
{\sc J.~T. Beale, W.~Ying, and J.~R. Wilson}, {\em A simple method for
  computing singular or nearly singular integrals on closed surfaces},
  Communications in Computational Physics, 20 (2016), pp.~733--753.

\bibitem{Borm:03}
{\sc S.~B{\"o}rm, L.~Grasedyck, and W.~Hackbusch}, {\em Introduction to
  hierarchical matrices with applications}, Engineering Analysis with Boundary
  Elements, 27 (2003), pp.~405--422.

\bibitem{Borm:17}
{\sc S.~B{\"o}rm and J.~M. Melenk}, {\em Approximation of the high-frequency
  {H}elmholtz kernel by nested directional interpolation: error analysis},
  Numerische Mathematik, 137 (2017), pp.~1--34.

\bibitem{bremer2015high}
{\sc J.~Bremer, A.~Gillman, and P.-G. Martinsson}, {\em A high-order accurate
  accelerated direct solver for acoustic scattering from surfaces}, BIT
  Numerical Mathematics, 55 (2015), pp.~367--397.

\bibitem{bremer2012nystrom}
{\sc J.~Bremer and Z.~Gimbutas}, {\em A {N}ystr{\"o}m method for weakly
  singular integral operators on surfaces}, Journal of Computational Physics,
  231 (2012), pp.~4885--4903.

\bibitem{bremer2013numerical}
{\sc J.~Bremer and Z.~Gimbutas}, {\em On the numerical evaluation of the
  singular integrals of scattering theory}, Journal of Computational Physics,
  251 (2013), pp.~327--343.

\bibitem{BrunoGarza:20}
{\sc O.~P. Bruno and E.~Garza}, {\em A {C}hebyshev-based rectangular-polar
  integral solver for scattering by geometries described by non-overlapping
  patches}, Journal of Computational Physics, 421 (2020), p.~109740.

\bibitem{Bruno:04}
{\sc O.~P. Bruno and E.~M. Hyde}, {\em An efficient, preconditioned, high-order
  solver for scattering by two-dimensional inhomogeneous media}, Journal of
  Computational Physics, 200 (2004), pp.~670--694.

\bibitem{Bruno:19}
{\sc O.~P. Bruno and A.~Pandey}, {\em Direct/iterative hybrid solver for
  scattering by inhomogeneous media}, SIAM Journal on Scientific Computing, 46
  (2024), pp.~A1298--A1326.

\bibitem{Bruno:22}
{\sc O.~P. Bruno and J.~Paul}, {\em Two-dimensional {F}ourier continuation and
  applications}, SIAM Journal on Scientific Computing, 44 (2022),
  pp.~A964--A992.

\bibitem{Ciarlet:72}
{\sc P.~G. Ciarlet and P.-A. Raviart}, {\em General {Lagrange} and {Hermite}
  interpolation in $\mathbb{R}^n$ with applications to finite element methods},
  Archive for Rational Mechanics and Analysis, 46 (1972), pp.~177--199.

\bibitem{COLTON:2012}
{\sc D.~Colton and R.~Kress}, {\em Inverse Acoustic and Electromagnetic
  Scattering Theory}, vol.~93, Springer, third~ed., 2012.

\bibitem{Conn:08}
{\sc A.~R. Conn, K.~Scheinberg, and L.~N. Vicente}, {\em Geometry of
  interpolation sets in derivative free optimization}, Mathematical
  Programming, 111 (2008), pp.~141--172.

\bibitem{Dangal:17monom}
{\sc T.~Dangal, C.-S. Chen, and J.~Lin}, {\em Polynomial particular solutions
  for solving elliptic partial differential equations}, Computers \&
  Mathematics with Applications, 73 (2017), pp.~60--70.

\bibitem{davis1965ignoring}
{\sc P.~J. Davis and P.~Rabinowitz}, {\em Ignoring the singularity in
  approximate integration}, Journal of the Society for Industrial and Applied
  Mathematics, Series B: Numerical Analysis, 2 (1965), pp.~367--383.

\bibitem{NIST:DLMF}
{\em {\it NIST Digital Library of Mathematical Functions}}.
\newblock \url{https://dlmf.nist.gov/}, Release 1.2.0 of 2024-03-15,
  \url{https://dlmf.nist.gov/}.
\newblock F.~W.~J. Olver, A.~B. {Olde Daalhuis}, D.~W. Lozier, B.~I. Schneider,
  R.~F. Boisvert, C.~W. Clark, B.~R. Miller, B.~V. Saunders, H.~S. Cohl, and
  M.~A. McClain, eds.

\bibitem{driveranalysis}
{\sc B.~Driver}, {\em Analysis Tools with Examples}, 2004,
  \url{https://mathweb.ucsd.edu/~bdriver/DRIVER/Book/anal.pdf}.

\bibitem{Ethridge:01}
{\sc F.~Ethridge and L.~Greengard}, {\em A new fast-multipole accelerated
  {P}oisson solver in two dimensions}, SIAM Journal on Scientific Computing, 23
  (2001), pp.~741--760.

\bibitem{faria2021general}
{\sc L.~M. Faria, C.~P{\'e}rez-Arancibia, and M.~Bonnet}, {\em General-purpose
  kernel regularization of boundary integral equations via density
  interpolation}, Computer Methods in Applied Mechanics and Engineering, 378
  (2021), p.~113703.

\bibitem{folland}
{\sc G.~Folland}, {\em Introduction to Partial Differential Equations},
  vol.~102, Princeton University Press, 1995.

\bibitem{Fryklund:22}
{\sc F.~Fryklund and L.~Greengard}, {\em An {FMM} accelerated {P}oisson solver
  for complicated geometries in the plane using function extension}, arXiv
  preprint arXiv:2211.14537,  (2022).

\bibitem{Fryklund:18}
{\sc F.~Fryklund, E.~Lehto, and A.-K. Tornberg}, {\em Partition of unity
  extension of functions on complex domains}, Journal of Computational Physics,
  375 (2018), pp.~57--79.

\bibitem{ganzburg2002markov}
{\sc M.~I. Ganzburg}, {\em A {Markov}-type inequality for multivariate
  polynomials on a convex body}, Journal of Computational Analysis and
  Applications, 4 (2002), pp.~265--268.

\bibitem{Sauer:00}
{\sc M.~Gasca and T.~Sauer}, {\em Polynomial interpolation in several
  variables}, Advances in Computational Mathematics, 12 (2000), pp.~377--410.

\bibitem{Golberg:03monom}
{\sc M.~Golberg, A.~Muleshkov, C.~Chen, and A.-D. Cheng}, {\em Polynomial
  particular solutions for certain partial differential operators}, Numerical
  Methods for Partial Differential Equations: An International Journal, 19
  (2003), pp.~112--133.

\bibitem{gomez2021regularization}
{\sc V.~G{\'o}mez and C.~P{\'e}rez-Arancibia}, {\em On the regularization of
  {C}auchy-type integral operators via the density interpolation method and
  applications}, Computers \& Mathematics with Applications, 87 (2021),
  pp.~107--119.

\bibitem{Gordon:73a}
{\sc W.~J. Gordon and C.~A. Hall}, {\em Construction of curvilinear co-ordinate
  systems and applications to mesh generation}, International Journal for
  Numerical Methods in Engineering, 7 (1973), pp.~461--477.

\bibitem{Gordon:73b}
{\sc W.~J. Gordon and C.~A. Hall}, {\em Transfinite element methods:
  blending-function interpolation over arbitrary curved element domains},
  Numerische Mathematik, 21 (1973), pp.~109--129.

\bibitem{GreengardLee:96}
{\sc L.~Greengard and J.-Y. Lee}, {\em A direct adaptive {P}oisson solver of
  arbitrary order accuracy}, Journal of Computational Physics, 125 (1996),
  pp.~415--424.

\bibitem{greengard2021fast}
{\sc L.~Greengard, M.~O'Neil, M.~Rachh, and F.~Vico}, {\em Fast multipole
  methods for the evaluation of layer potentials with locally-corrected
  quadratures}, Journal of Computational Physics: X, 10 (2021), p.~100092.

\bibitem{Greengard:87}
{\sc L.~Greengard and V.~Rokhlin}, {\em A fast algorithm for particle
  simulations}, Journal of Computational Physics, 73 (1987), pp.~325--348.

\bibitem{helsing2008evaluation}
{\sc J.~Helsing and R.~Ojala}, {\em On the evaluation of layer potentials close
  to their sources}, Journal of Computational Physics, 227 (2008),
  pp.~2899--2921.

\bibitem{IsaacsonKeller}
{\sc E.~Isaacson and H.~B. Keller}, {\em Analysis of numerical methods}, New
  York: Wiley,  (1966).

\bibitem{Israeli:02}
{\sc M.~Israeli, E.~Braverman, and A.~Averbuch}, {\em A hierarchical {3-D}
  {P}oisson modified {F}ourier solver by domain decomposition}, Journal of
  Scientific Computing, 17 (2002), pp.~471--479.

\bibitem{klockner2013quadrature}
{\sc A.~Kl{\"o}ckner, A.~Barnett, L.~Greengard, and M.~O'Neil}, {\em Quadrature
  by expansion: {A} new method for the evaluation of layer potentials}, Journal
  of Computational Physics, 252 (2013), pp.~332--349.

\bibitem{kress1978transmission}
{\sc R.~Kress and G.~Roach}, {\em Transmission problems for the {H}elmholtz
  equation}, Journal of Mathematical Physics, 19 (1978), pp.~1433--1437.

\bibitem{lubinsky1984rates}
{\sc D.~S. Lubinsky and P.~Rabinowitz}, {\em Rates of convergence of {G}aussian
  quadrature for singular integrands}, {M}athematics of {C}omputation, 43
  (1984), pp.~219--242.

\bibitem{marsden1993elementary}
{\sc J.~E. Marsden and M.~J. Hoffman}, {\em Elementary Classical Analysis},
  Macmillan, 1993.

\bibitem{martin2003acoustic}
{\sc P.~A. Martin}, {\em Acoustic scattering by inhomogeneous obstacles}, SIAM
  Journal on Applied Mathematics, 64 (2003), pp.~297--308.

\bibitem{Matthys:96monom}
{\sc L.~Matthys, H.~Lambert, and G.~De~Mey}, {\em A recursive construction of
  particular solutions to a system of coupled linear partial differential
  equations with polynomial source term}, Journal of Computational and Applied
  Mathematics, 69 (1996), pp.~319--329.

\bibitem{McCorquodale:07}
{\sc P.~McCorquodale, P.~Colella, G.~Balls, and S.~Baden}, {\em A local
  corrections algorithm for solving {P}oisson's equation in three dimensions},
  Communications in Applied Mathematics and Computational Science, 2 (2007),
  pp.~57--81.

\bibitem{Nardini:83}
{\sc D.~Nardini and C.~Brebbia}, {\em A new approach to free vibration analysis
  using boundary elements}, Applied Mathematical Modelling, 7 (1983),
  pp.~157--162.

\bibitem{Olver:06}
{\sc P.~J. Olver}, {\em On multivariate interpolation}, Studies in Applied
  Mathematics, 116 (2006), pp.~201--240.

\bibitem{Partridge}
{\sc P.~W. Partridge and C.~A. Brebbia}, {\em Dual reciprocity boundary element
  method}, Springer Science \& Business Media, 2012.

\bibitem{perez2018plane}
{\sc C.~P{\'e}rez-Arancibia}, {\em A plane-wave singularity subtraction
  technique for the classical {D}irichlet and {N}eumann combined field integral
  equations}, Applied Numerical Mathematics, 123 (2018), pp.~221--240.

\bibitem{perez2019harmonic}
{\sc C.~P{\'e}rez-Arancibia, L.~M. Faria, and C.~Turc}, {\em Harmonic density
  interpolation methods for high-order evaluation of {L}aplace layer potentials
  in {2D} and {3D}}, Journal of Computational Physics, 376 (2019),
  pp.~411--434.

\bibitem{perez2019planewave}
{\sc C.~P{\'e}rez-Arancibia, C.~Turc, and L.~Faria}, {\em Planewave density
  interpolation methods for {3D} {H}elmholtz boundary integral equations}, SIAM
  Journal on Scientific Computing, 41 (2019), pp.~A2088--A2116.

\bibitem{rabinowitz1967gaussian}
{\sc P.~Rabinowitz}, {\em Gaussian integration in the presence of a
  singularity}, SIAM {J}ournal on {N}umerical {A}nalysis, 4 (1967),
  pp.~191--201.

\bibitem{rabinowitz1986rates}
{\sc P.~Rabinowitz}, {\em Rates of convergence of {G}auss, {L}obatto, and
  {R}adau integration rules for singular integrands}, {M}athematics of
  {C}omputation, 47 (1986), pp.~625--638.

\bibitem{saad1986gmres}
{\sc Y.~Saad and M.~H. Schultz}, {\em {GMRES}: {A} generalized minimal residual
  algorithm for solving nonsymmetric linear systems}, SIAM Journal on
  Scientific and Statistical Computing, 7 (1986), pp.~856--869.

\bibitem{saranen2001periodic}
{\sc J.~Saranen and G.~Vainikko}, {\em Periodic integral and pseudodifferential
  equations with numerical approximation}, Springer Science \& Business Media,
  2001.

\bibitem{Sauer:95}
{\sc T.~Sauer and Y.~Xu}, {\em On multivariate {Lagrange} interpolation},
  Mathematics of Computation, 64 (1995), pp.~1147--1170.

\bibitem{Sauter2010}
{\sc S.~A. Sauter and C.~Schwab}, {\em Boundary Element Methods}, Springer,
  2010.

\bibitem{ShenSerkh:22}
{\sc Z.~Shen and K.~Serkh}, {\em Rapid evaluation of {N}ewtonian potentials on
  planar domains}, SIAM Journal on Scientific Computing, 46 (2024),
  p.~A609–A628.

\bibitem{siegel2018local}
{\sc M.~Siegel and A.-K. Tornberg}, {\em A local target specific quadrature by
  expansion method for evaluation of layer potentials in {3D}}, Journal of
  Computational Physics, 364 (2018), pp.~365--392.

\bibitem{Stein:22}
{\sc D.~B. Stein}, {\em Spectrally accurate solutions to inhomogeneous elliptic
  {PDE} in smooth geometries using function intension}, arXiv preprint
  arXiv:2203.01798,  (2022).

\bibitem{Stein:16}
{\sc D.~B. Stein, R.~D. Guy, and B.~Thomases}, {\em Immersed boundary smooth
  extension: a high-order method for solving {PDE} on arbitrary smooth domains
  using {F}ourier spectral methods}, Journal of Computational Physics, 304
  (2016), pp.~252--274.

\bibitem{stein1970singular}
{\sc E.~M. Stein}, {\em Singular Integrals and Differentiability Properties of
  Functions}, Princeton University Press, 1970.

\bibitem{Strang:72}
{\sc G.~Strang}, {\em Approximation in the finite element method}, Numerische
  Mathematik, 19 (1972), pp.~81--98.

\bibitem{sushnikova2023fmm}
{\sc D.~Sushnikova, L.~Greengard, M.~O’Neil, and M.~Rachh}, {\em {FMM-LU}: A
  fast direct solver for multiscale boundary integral equations in three
  dimensions}, Multiscale Modeling \& Simulation, 21 (2023), pp.~1570--1601.

\bibitem{taylor2013partial}
{\sc M.~Taylor}, {\em Partial Differential Equations I: Basic Theory},
  vol.~116, Springer Science \& Business Media, 2013.

\bibitem{vainikko2006multidimensional}
{\sc G.~Vainikko}, {\em Multidimensional Weakly Singular Integral Equations},
  Springer, 2006.

\bibitem{Vioreanu:14}
{\sc B.~Vioreanu and V.~Rokhlin}, {\em Spectra of multiplication operators as a
  numerical tool}, SIAM Journal on Scientific Computing, 36 (2014),
  pp.~A267--A288.

\bibitem{wilhelmsen:74}
{\sc C.~Wilhelmsen}, {\em A {M}arkov inequality in several dimensions}, J.
  Approx. Theory, 11 (1974), pp.~216--220.

\bibitem{xiang2012convergence}
{\sc S.~Xiang and F.~Bornemann}, {\em On the convergence rates of {G}auss and
  {C}lenshaw--{C}urtis quadrature for functions of limited regularity}, SIAM
  journal on numerical analysis, 50 (2012), pp.~2581--2587.

\bibitem{Xiao:10}
{\sc H.~Xiao and Z.~Gimbutas}, {\em A numerical algorithm for the construction
  of efficient quadrature rules in two and higher dimensions}, Computers \&
  mathematics with applications, 59 (2010), pp.~663--676.

\bibitem{ying2006high}
{\sc L.~Ying, G.~Biros, and D.~Zorin}, {\em A high-order 3{D} boundary integral
  equation solver for elliptic {PDE}s in smooth domains}, Journal of
  Computational Physics, 219 (2006), pp.~247--275.

\bibitem{zhu2022high}
{\sc H.~Zhu and S.~Veerapaneni}, {\em High-order close evaluation of {L}aplace
  layer potentials: {A} differential geometric approach}, SIAM Journal on
  Scientific Computing, 44 (2022), pp.~A1381--A1404.

\end{thebibliography}
